\documentclass[12pt,leqno]{article}
\usepackage{amsmath}
\usepackage{amssymb}
\usepackage{amsthm} 
\usepackage{epsf}
\usepackage{physics} 
\usepackage{braket} 
\usepackage{graphicx}
\usepackage{esint} 
\usepackage{dsfont}
\usepackage{hyperref} 
\hypersetup{pdfpagemode=FullScreen,  colorlinks=true} 

\usepackage[breakable]{tcolorbox}
\usepackage{xcolor}
\usepackage{soulutf8}
\soulregister\cite7
\soulregister\underline7
\soulregister\eqref7
\soulregister\ref7

\definecolor{mycolor}{RGB}{229,228,200}
\definecolor{pinkcolor}{RGB}{247,202,193}

\newtcolorbox{changebox}{
    colback=mycolor,    
    colframe=mycolor,   
    boxrule=0pt,
    sharp corners,      
    left=2pt,
    right=2pt,
    top=2pt,
    bottom=2pt,
    breakable}

\numberwithin{equation}{section}

\newtheorem{lem}{Lemma}[section]
\newtheorem{pro}[lem]{Proposition}
\newtheorem{defi}[lem]{Definition}
\newtheorem{thm}[lem]{Theorem}

\theoremstyle{remark}
\newtheorem{rem}[lem]{\bf Remark}

\newcommand{\ms}{\medskip}
\newcommand{\msi}{\par\medskip\noindent}

\newcommand{\R}{\mathbb{R}}
\renewcommand{\H}{\mathcal H}

\newcommand{\bY}{\mathbb {Y}}
\newcommand{\bT}{\mathbb {T}}
\newcommand{\bV}{\mathbb {V}}
\newcommand{\bP}{\mathbb {P}}
\newcommand{\bH}{\mathbb {H}}

\renewcommand{\d}{\partial}
\newcommand{\dist}{\,\mathrm{dist}}

\newcommand{\Angle}{\,\mathrm{Angle}}
\newcommand{\sm}{\setminus}

\newcommand{\wt}{\widetilde}
\newcommand{\wh}{\widehat}
\newcommand{\ol}{\overline}
\newcommand{\ub}{\underbar}

\newcommand{\cL}{{\cal L}}

\newcommand{\cX}{{\cal X}}
\newcommand{\cE}{{\cal E}}

\newcommand{\cP}{{\cal P}}

\newcommand{\cG}{{\cal G}}
\newcommand{\cB}{{\cal B}}

\newcommand{\cC}{{\mathfrak C}}

\begin{document}

\title{An existence theorem for sliding minimal sets}

\author{
    G. David\footnote{ 
    G. David was partially supported by the European Community H2020 grant GHAIA 777822, and the Simons Foundation grant 601941, GD.} 
    \, and C. Labourie
    }

    \date{}
    \maketitle

    \ms\noindent{\bf Abstract.}
    We prove an existence theorem for the sliding  boundary variant of the Plateau problem for $2$-dimensional sets in $\mathbb{R}^n$. The simplest case of sufficient condition is when $n=3$ and the boundary $\Gamma$ is a finite disjoint union of smooth closed curves contained in the boundary of a convex body, but the main point of our sufficient condition is to prevent the limits in measure of a minimizing sequence to have singularities of type $\mathbb{Y}$ along $\Gamma$.

    \ms\noindent{\bf Résumé en Français.}
    On démontre un résultat d'existence pour la variante à frontière glissante du problème de Plateau pour un ensemble de dimension $2$ dans $\mathbb{R}^n$. La condition d'existence la plus simple est quand $n=3$ et on demande que la frontière $\Gamma$ soit une union finie disjointe de courbes lisses fermées, mais le  but principal de notre condition suffisante est d'empêcher que les limites en mesure de suites minimisantes aient des singularités de type $\mathbb{Y}$ le long de $\Gamma$.

    \ms\noindent{\bf Key words/Mots cl\'es.}
    Sliding Plateau problem, Existence, Minimal sets of dimension $2$.

    \ms\noindent

    \tableofcontents

    \section{Introduction}
    \label{S1}

    Let a compact set $\Gamma$ be a finite union of disjoint smooth closed curves (i.e., loops) in $\R^n$; 
    we consider the following variant of the Plateau problem.
    We start with an initial compact set $E_0 \subset \R^n$.
    Then  we consider the  class $\cE = \cE(E_0,\Gamma)$ of sets obtained from $E_0$
    by a continuous deformation that preserves $\Gamma$. That is, we say that $E \in \cE$  
    when we  can find a continuous mapping $\varphi: [0,1]  \times E_0 \to \R^n$, with the properties
    \begin{equation} \label{1a1}
        \varphi(0,x)  = x \ \text{ for } x\in E_0,
    \end{equation}
    \begin{equation} \label{1a2}
        \varphi(t,x)  \in \Gamma \ \text{ when $x\in E_0 \cap \Gamma$ and } 0 \leq t \leq 1
    \end{equation}
    and, denoting $\varphi_t(x)  = \varphi(t,x)$, 
    \begin{equation} \label{1a3}
        \varphi_1 \text{ is Lipschitz on $E_0$,}
    \end{equation}
    and finally
    \begin{equation} \label{1a4}
        E = \varphi_1(E_0).
    \end{equation}
    The condition \eqref{1a3} is traditional. Our result is a tiny bit better with it (we will eventually construct a minimizer $E$ and a $\varphi$ such that $\varphi_1$ is  Lipschitz). 
    Moreover, the limiting set in Theorem \ref{t2a1} is a priori only minimal with respect to Lipschitz deformation.

    We look for minimizers of $\H^2(E)$ (or sometimes a slightly different functional $J$) in the class $\cE$, i.e., sets $E \in \cE$ such that $\H^2(E) = m_0$, where 
    \begin{equation} \label{1a5}
        m_0 = m_0(E_0, \Gamma) = \inf \set{\H^2(E) | E \in \cE}.
    \end{equation}
    Here $\H^2$ denotes the $2$-dimensional Hausdorff measure (recall that $\H^2$ coincides 
    with the surface measure on smooth sets $E$). In general, the sliding Plateau problem would 
    concern more general boundary sets $\Gamma$, and we would minimize the Hausdorff 
    measure $\H^d(E)$, where $d \in \{1,\ldots,n-1\}$, but here our method will force us to 
    restrict to $d=2$.
    We will refer to the minimizers, if they exist, as solutions of the sliding Plateau problem associated
    to $\Gamma$ and $E_0$. The main result of this paper is as follows.

    \begin{thm}\label{t1}
        Let a compact set $\Gamma \subset \R^n$ be the union of a finite family of disjoint closed curves of class 
        $C^{1+\alpha}$, with $\alpha > 0$. Assume in addition that 
        \begin{equation} \label{1a6}
            \text{$\Gamma$ has a good access to the complement of its convex hull;}
        \end{equation}
        we will explain what this means in Section \ref{S3}.
        Let $E_0$ be any compact subset of $\R^n$, and define $\cE = \cE(E_0,\Gamma)$
        as above. Then there exists a solution for the sliding Plateau problem associated
        to $\Gamma$ and $E_0$, i.e., a compact set $E \in \cE$ such that $\H^2(E) = m_0$.
    \end{thm}

    We will explain later in this text that this result can be generalized to some extent, 
    but as far as we know Theorem \ref{t1} is already new as stated, even when we assume 
    that $n=3$ and $\Gamma$ is contained in the boundary of its convex hull. It 
    also gives a very good idea of what we can do. 
    See Theorem \ref{t3a1} for an example that is not too explicit and 
    Section \ref{Slast} concerning extensions of Theorem \ref{t1}.

    When $n=3$, \eqref{1a6} holds as soon as $\Gamma \subset \d H(\Gamma)$, where $H(\Gamma)$ is the convex hull of $\Gamma$.
    When $n >3$, the condition $\Gamma \subset \d H(\Gamma)$ does not mean
    much, because $\Gamma$ (and then everything we do) could be contained in a hyperplane, and 
    then the condition would be void. We will settle for a definition of \eqref{1a6} that looks reasonable
    and allows to work in $\R^n$, $n > 4$. See Definition~\ref{d3a1}.

    Our assumption that the curves that compose $\Gamma$ are $C^{1+\alpha}$ does not seem
    shocking, even though it is probably not optimal; it will allow us to apply some regularity results from \cite{Dvv}. 
    Also, we should probably be able to allow the curves to cross, but the regularity theorem 
    near a cross that we would need is not written yet.

    The main assumption \eqref{1a6} is more problematic. 
    This type of condition is not new; it appears in \cite{Mo1}, in a context similar to ours, 
    but where the sliding minimizers are replaced with size-minimizing currents and
  the minimized surfaces are of co-dimension $1$ instead of $n-2$ 
  (which  makes the definition of \eqref{1a6} simpler).
    In both cases the point is that we can restrict our attention to competitors that are contained in 
    the convex hull of $\Gamma$. We find it interesting that \eqref{1a6} arises almost 
    naturally in our context of sliding minimizers,
    as a way to prevent some annoying blow-up limits from arising.
    At this point the authors still hope that \eqref{1a6} is  not needed, but are unable to prove 
    this because the regularity result in \cite{Dvv} is not good enough. 

    \ms

    Let us add some additional general comments about our sliding Plateau problem before
    we discuss the  proof.
    Typically, our problem will not be trivial when $E_0$ is related to $\Gamma$ with some topological condition, but the pleasant feature of the sliding Plateau problem is that we don't need to say which topological condition.
    We like to keep open the choice of the initial set $E_0$, because we like the fact that
    different initial sets $E_0$ give different classes $\cE$ and often different minimizers $E$.
    Even the value of $m_0$ will depend on $E_0$ in general.

    The reader may worry about bad choices of $E_0$.
    If we start from a stupid choice of $E_0$, it may happen that $m_0 = 0$ and the result is not really interesting.
    This is what happens in the extreme case when $E_0$ is a compact set that does not meet $\Gamma$, and thus there is a retraction
    $\{ \varphi_t \}$ of $E_0$ to a point. In this case the sliding condition \eqref{1a2} is void, 
    and does not help for the non-degeneracy of our Plateau problem. 
    It could also happen that $m_0 = 0$ if $E_0$ contains $\Gamma$, but can still be deformed 
    into a trivial set $E \in \cE$, either reduced to  a point or of dimension~$1$. 
    In any case, we will deal with the case $m_0 = 0$ with a specific argument (at the end of Section \ref{S3})
    
    It could also a priori  happen that $m_0 = +\infty$, i.e., $\H^2(E) = +\infty$ for every $E \in \cE$.
    In this case the theorem would be essentially void (take $E = E_0$), but anyway we claim that this cannot happen due to the fact
    that $\Gamma$ is nice curve and $E_0$ is compact. In this setting, we can construct a Federer-Fleming projection 
    $\varphi$ (on biLipschitz images of faces of cubes) so that $E = \varphi(E_0)$ lies in
    $\cE$ and $\H^2(E) < +\infty$. We leave the details to the reader, because the only point 
    of the remark is to feel a little better about the statement. 

    Our Plateau problem looks a lot like the problem of size minimizers, where one chooses
    a reasonable integral current $S$ on $\Gamma$ (for instance, any sum of integer multiples of the current of 
    integration on the curves that compose it), and then looks for a current $T$ of minimal size such that
    $d T = S$. See for instance \cite{Mo1}; the algebra here is different (in fact, one could argue that
    it is not even visible here), but the problems have strong similarity, because in both cases we minimize
    the Hausdorff measure of the relevant  geometric object. Notice however  that  we seem to have
    more flexibility in the definition of the problem, and in particular we do not care about orientability.
    We refer to \cite{DStein} for a discussion and comparison of the various classical Plateau problems.

    Of course our problem looks even more like the initial Plateau problem treaded by
    Radó \cite{Ra1,Ra2} and Douglas \cite{Do}, for instance, where $E$ is given with a parameterization 
    $f$, and its area is computed from $f$. Here we allow $f$ to be defined on $E_0$, and use slightly
    different rules on the list of allowed competitors; for instance, they only allowed injective mappings,
    or else used a (local) formula that compute the area of the image with multiplicity when $f$ is not injective.
    In both  cases, the difficulty with parameterizations (or our map $\varphi_1$) is that the modulus of continuity 
    of the parameterization may degenerate along a minimizing sequence, thus making it hard to 
    parameterize a limit. 

    In the case when $\Gamma$ is just one closed curve, Radó and Douglas consider sets
    $E = f(D)$, where $D$ is the closed unit disk in $\R^2$ and $f$ is a continuous function
    which coincides on the circle $\d D$ with a nice parameterization of $\Gamma$ by $\d D$.
    Here we could take for $E_0$ such a set $E$ (and taking $f$ injective will be nicer), but our problem 
    is  different because  we allow non injective deformations. 
    Also, some of the minimal sets that we look for are
    more naturally parameterized by some other surfaces than a closed disk, such as a more general
    torus with $m$ holes, minus a small disk whose boundary is sent to $\Gamma$. 
    Thus we want to allow $E_0$ to be one of these objects. As we said before, different choices of $E_0$
    may yield interestingly different minimizers $E$ (and even different values of $m_0$); see \cite{DStein}.

    Our proof of Theorem \ref{t1} relies on two ingredients: a regularity theorem along the boundary for
    sliding almost minimal sets \cite{Dvv}, 
    and a stability theorem concerning limits of such sets \cite{La}.
    The general principle is, as often for the existence of minimizers, to start from a minimizing sequence
    $\{ E_k \}$, extract a sequence with a limit, and show that the limit is a minimizer. 

    There are a few a priori difficulties with this general program, some of them less annoying 
    than one could expect, some of them more delicate. 
    The most obvious attempt would be to take a subsequence $\{ E_k \}$
    that converges for the Hausdorff distance to some limit $E_\infty$, but if we do this, the most likely
    is that $\H^2(E_\infty) = +\infty$, because $E_k$ may have lots of small hairs, with small area, but
    that converge to a large set. The next attempt is to consider a special subsequence with additional 
    regularity properties, so that the measure $\H^2_{\vert E_k}$ has suitable lower semicontinuity
    properties. This was for instance a strategy followed by Reifenberg in \cite{Re}, and in the context 
    of sliding minimal sets as here, this was proposed for instance in \cite{Limits}, and implemented 
    in specific instances in \cite{Li} and \cite{Fv1,Fv2}. The construction of a sequence $\{ E_k \}$
    of ``better competitors'', can be long and  painful, but the general idea is that along a sequence of
    quasiminimal sets, with uniform quasiminimality constants, the amount of Hausdorff measure 
    in an open set behaves in a lower semicontinuous way.

    More recently, the authors of \cite{Ita1} and then \cite{Ita2} found out that for such problems, using
    the weak convergence of the measures $\mu_k = \H^2_{\vert E_k}$ is often much more convenient.
    The idea now is that any weak limit of the $\mu_k$ is of the form $\H^2_{\vert E_\infty}$
    for some minimal set $E_\infty$. Their proof typically uses well known (but difficult) 
    results on minimal sets, plus simpler arguments suited to the precise context. 
    In the present situation, their proof typically yields a set $E_\infty$ which is sliding minimal 
    (in the sense of \cite{Dss}), 
    but this is not yet enough for us because maybe $E_\infty$ is no longer a competitor in 
    our class $\cE$.
    Indeed, we do not know whether we can make the deformations $\varphi_k$ associated to the 
    $E_k$ converge to anything, and also pieces of surface in the $E_k$ may converge to a 
    set of dimension smaller than $d$ (think about thin tubes converging to a wire), 
    which is important in the definition of $\cE$ but disappears from the limit of the $\mu_k$.

    In our case, we will still use weak limits of the Radon measures $\mu_k$, but for the proof of
    sliding minimality for $E_\infty$, we cite \cite{La} which takes into account the sliding boundary and whose proof is more flexible.

    So we need a last piece of information, coming from the fact that $E_\infty$ is a sliding
    minimizer. If none of the blow-up limits of $E_\infty$ at a point of $E_\infty \cap \Gamma$ contains
    a cone of type $\bY$ whose spine is parallel to the tangent of $\Gamma$ 
    (see the definition below), then we can apply a regularity result 
    from \cite{Dvv} 
    and prove that there exist local retractions on $E_\infty$, that preserve $\Gamma$ too.
    Then we can compose such retractions with any of the mappings $\varphi_k$, $k$ large
    enough, to get a deformation $\varphi$ that maps $E_0$ to the limit $E_\infty$
    (in fact, plus a set of vanishing measure), prove that this set lies in $\cE$, and conclude. 
    This last part follows the same route that was used in 
    \cite{Li} and \cite{Fv1,Fv2}, although they worked in different contexts where the retraction 
    was probably harder to find.

    The point of our accessibility assumption (\ref{1a6}) is that 
    the regularity result in \cite{Dvv} is only stated for almost minimal sets that do not have
    blow-up limits of type $\bY$ along $\Gamma$; the extra assumption \eqref{1a6} is precisely a
    simple way to make sure that $E_\infty$ is like that. 
    Of course, the reason why we restricted to $2$-dimensional sets is that \cite{Dvv} only works 
    in this dimension.

    The plan for the rest of this paper is as follows. In Section \ref{S2} we give some of the missing
    definitions, pick any minimizing sequence $\{ E_k \}$, then choose a subsequence
    so that the $\mu_k$ converge weakly to a measure $\mu$, and use the results of
    \cite{La} to show in Theorem \ref{t2a1} that $\mu = \H^2_{\vert E}$ for some sliding minimal 
    set $E$.

    In Section \ref{S3}, we define the notion of good access and state our main practical result,
    Theorem \ref{t3a1}, which says that if we can find a minimizing sequence in a compact set $K$
    (think about the convex hull of $\Gamma$), such that $\Gamma$ has good access to the complement of $K$, then we can find a 
    sliding minimizer in $\cE$. Other functionals $J$ are allowed (as in \eqref{2a5}-\eqref{2a7}),
    but, other than the convex case, the reasons why there would be a compact set $K$
    as above are not discussed before Section \ref{Slast}.

    The proof of Theorem \ref{t3a1} (which implies Theorem \ref{t1}) is done in 
    Sections \ref{S4}-\ref{S10}, where the regularity result in \cite{Dvv} is used to
    control the geometry of the limit set $E_\infty$ (the support of the limit of
    the measures $\H^2_{\vert E_k}$ for a minimizing sequence), and eventually construct a 
    local retraction on $E_\infty$. Section \ref{Slast} contains a discussion of 
    circumstances where we can use Theorem \ref{t3a1}, yielding existence results
    that generalize Theorem \ref{t1}, with essentially the same proof. 

    Finally, we prove in Section \ref{Sbil} that if $E$ is a coral almost minimal set
    of dimension $2$ in $\R^n$, with a sliding boundary composed of disjoint $C^{1+\alpha}$
    closed curves, such that all the blow-up limits of $E$ at $0 \in E$ 
    only have transverse, or half plane, or generic $\bV$ behaviors at the origin, 
    then there is a neighborhood of $0$ where $E$ is biLipschitz-equivalent to one (in fact any) of the tangent cones to $E$ at $x$. 
    This result is not needed for the rest of the paper, but the proof of Theorem \ref{t3a1} 
    almost gives it, and 
    it seems interesting, because this is a situation where we don't know wether there is a 
    unique blow-up limit of $E$ at $0$ (only, they are all biLipschitz-equivalent to each other).

    \section{Almost minimal sets and minimizing sequences}
    \label{S2}

    Since we may also consider slightly more general sliding boundary problems, we shall
    write down the definitions, and some proofs, in more generality than needed for the official
    results of this paper. It is also important to us that we use almost-minimality
    arguments, which are more flexible than the arguments involving true minimality.

    For the moment, we fix an integer $d \in \{1,\ldots,n-1\}$ and a compact boundary set 
    $\Gamma \subset \R^n$ (not necessarily of dimension $d-1$). To each compact set
    $E_0 \subset \R^n$ such that $\H^d(E_0) < \infty$,
    we associate as above the class $\cE(E_0,\Gamma)$ of 
    (results of) deformations of $E_0$ that preserve $\Gamma$.
    That is, $E \in \cE(E_0,\Gamma)$ when $E = \varphi_1(E_0)$ for some family $\{ \varphi_t \}$,
    $0 \leq t \leq 1$, that satisfies \eqref{1a1}-\eqref{1a4}. 

    Let $h: (0,+\infty) \to [0,+\infty]$ be a nondecreasing function such that $\lim_{r \to 0} h(r) = 0$ 
    (we shall call this a \ub{gauge function}); we shall assume in addition that there exist $\alpha > 0$,
    $c_h \geq 0$, and $r_h > 0$ such that
    \begin{equation} \label{2a1}
        h(r) \leq c_h r^{\alpha} \ \text{ for } 0 < r \leq r_h, 
    \end{equation}
    so that various earlier results can be easily applied.

    \begin{defi}\label{d2a1}
        Let $E$ be a compact set of $\R^n$ with $\H^d(E) < \infty$.
        We say that $E \subset \R^n$ is a \ub{sliding almost minimal set} of dimension $d$, 
        with gauge function $h$ and sliding boundary $\Gamma$, when for each choice of ball $B(y,r) \subset \R^n$ with $0 < r \leq r_h$, and each family $\{ \varphi_t \}$ 
        that satisfies the conditions \eqref{1a1}-\eqref{1a4} with $E_0$ replaced by $E$, 
        and in addition the deformation happens entirely in $B(y,r)$, i.e., when for $0 \leq t \leq 1$, 
        \begin{equation} \label{2a2}
            \varphi_t(x,t) = x \text{ for $x \notin B(y,r)$ and }
            \varphi_t(E \cap B(y,r)) \subset B(y,r),
        \end{equation}
        we have that
        \begin{equation} \label{2a3}
            \H^d(E \cap B(y,r)) \leq \H^d(\varphi_1(E \cap B(y,r))) + h(r) r^{d}.
        \end{equation}
    \end{defi}

    A \ub{sliding minimal set} is just a sliding almost minimal set associated to the gauge function $h\equiv 0$ (and without radius constraint, i.e., $r_h = +\infty$); 
    thus for minimal sets, the definition is simpler and boils down to
    \begin{equation} \label{2a4}
        \H^d(E) \leq \H^d(F)  \ \text{ for every } F\in \cE(E,\Gamma)
    \end{equation}
    because $E$ and $F$ are compact sets, and since $h =0$ we may even take $B(y,r)$ very large
    without losing information in \eqref{2a3}. There exist local definitions of sliding minimal and almost
    minimal sets (see for instance \cite{Dss}), but we won't need them here.

    We will use the notion of (sliding) almost minimal sets because it is much more flexible than
    the notion of minimal sets. The simplest way it appears naturally is when you use functionals like
    \begin{equation} \label{2a5}
        J(E) = \int_E f(x) d\H^d(x)
    \end{equation}
    instead of $\H^d(E)$, which corresponds to $f \equiv 1$. For instance, if 
    $\Gamma$ is reasonable, $f$ is such that
    \begin{equation} \label{2a6}
        C_0^{-1} \leq f(x) \leq C_0  \ \text{ for } x \in \R^n
    \end{equation}
    and 
    \begin{equation} \label{2a7}
        |f(x)-f(y)| \leq C_1 |x-y|^\alpha \ \text{ for } x,y \in \R^n
    \end{equation}
    for some constants $C_0, C_1 \geq 1$, and if $E$ is a sliding minimizer for $J$
    (with the same definition \eqref{2a4} as above, but with $\H^2$ replaced by $J$), then we are going to see that $E$ is sliding almost minimal with a gauge function $h$ that satisfies \eqref{2a1} for some choice 
    of $c_h$ and $r_h$ (that depend on $\Gamma$, $n$, $C_0$, and $C_1$).
    A first easy consequence of these assumptions is that $E$ is \ub{quasiminimal}, which means that there exists a constant $M \geq 1$ such that for all ball $B(y,r)$ and all deformation $\varphi_1$ as in Definition \ref{d2a1}, we have
    \begin{equation}
        \H^d(E \cap W) \leq M \H^d(\varphi_1(E \cap W)),
    \end{equation}
    where $W = \set{x \in E \cap B(y,r) | \varphi_1(x) \ne x}$.
    The mild constraint on $\Gamma$ (satisfied for instance if $\Gamma$ is the biLipschitz image
    of a finite union of faces of dyadic cubes of various dimensions, hence also when $\Gamma$
    is as in Theorem \ref{t1}) is then used to verify that $E$ is locally Ahlfors regular.
    This means that we can find a constant $C_d \geq 1$ (which depends only on $n$, $M$ and $\Gamma$) and a radius $r_d > 0$ (which depends only on $n$ and $\Gamma$) such that for all $x \in E^\ast$ (the closed support of $H^d_{\vert E}$) and $0 \leq r \leq r_d$ we have
    \begin{equation} \label{2a8}
        C_d^{-1} r^d \leq \H^d(E \cap B(x,r)) \leq C_d r^d.
    \end{equation}
    See \cite[Proposition 4.74]{Dss} for the verification of \eqref{2a8}, which takes some time but is not surprising.

    Once we know {\eqref{2a8}}, the almost minimality of $E$ follows easily.
    Consider a ball $B(y,r)$ with $0 \leq r \leq r_d$ and a sliding competitor $F = \varphi_1(E) \in \cE(E,\Gamma)$, with $\varphi$ satisfying \eqref{2a2}. By minimality of $E$, we have $J(E) \leq J(F)$ and since $F$ coincides with $E$ in the complement of $B(y,r)$, this simplifies to
    \begin{equation}
        J(E \cap B(y,r)) \leq J(F \cap B(y,r)).
    \end{equation}
    If $\H^d(F \cap B(y,r) \geq C_d r^d$, then we automatically have $\H^d(E \cap B(y,r)) \leq \H^d(F \cap B(y,r))$ by (\ref{2a8}). If $\H^d(F \cap B(y,r)) \leq C_d r^d$, we have
    \begin{eqnarray} \label{2a9}
        f(y) \H^d(E \cap B(y,r)) &=& \int_{E \cap B(y,r)} f(y) d\H^d(x)\cr
        &\,& \hskip-1.2cm \leq \int_{E \cap B(y,r)} f(x) d\H^d(x) + C_1 r^\alpha \H^d(E \cap B(y,r))\cr
        &\,& \hskip-1.2cm \leq \int_{E \cap B(y,r)} f(x) d\H^d(x) + C_1 C_d r^{d+\alpha} 
        \leq \int_{F \cap B(y,r)} f(x) d\H^d(x) + C_1 C_d r^{d+\alpha}\cr
        &\,& \hskip-1.2cm \leq\int_{F \cap B(y,r)} f(y) d\H^d(x) + 2 C_1 C_d r^{d+\alpha} 
        = f(y) \H^d(F \cap B(y,r)) + 2 C_1 C_d r^{d+\alpha}.
    \end{eqnarray}
    We divide by $f(y)$, use \eqref{2a6}, and get \eqref{2a3}.

    We will start our argument searching for a $d$-dimensional set $E \in \cE(E_0,\Gamma)$, that 
    minimizes a functional like $J$ in \eqref{2a5} (with the assumptions \eqref{2a6} and \eqref{2a7}).
    Later in the argument we will make further assumptions as we need them.

    So let $\Gamma$, $E_0$, $J$, be given, set
    \begin{equation} \label{2a10}
        m_0 = m_0(E_0, \Gamma, f) = \inf\set{J(E) | E \in \cE(E_0,\Gamma)}
    \end{equation}
    as above, and consider a minimizing sequence $\{ E_k \}$, $k \geq 0$,
    in $\cE = \cE(E_0,\Gamma)$. That is,
    \begin{equation} \label{2a11}
        \lim_{k \to +\infty} J(E_k) = m_0.
    \end{equation}
    We want to use $\{ E_k \}$ to find a minimizer for $J$, 
    and rather than trying to find a subsequence that converges for the Hausdorff distance,
    we decide to follow \cite{Ita1, Ita2, La} and replace $\{ E_k \}$ with a subsequence for 
    which the measures $\mu_k = \H^{d}_{\vert E_k}$ converge weakly to some limit
    $\mu_\infty$.
    The purpose of this section is to prove that $\mu_\infty = \H^{d}_{\vert E_\infty}$
    for some sliding almost minimal set $E_\infty$.

    We intend to apply Theorem 3.3 or Corollary 4.1 in \cite{La}, so we should assume that our
    boundary set $\Gamma$ is what is referred to in \cite{La} as a Whitney set (Definition 1.8 there). It means a closed set $\Gamma$ which is a Lipschitz neighborhood retract, which is locally diffeomorphic to a cone and which is locally biLipschitz equivalent to an union of faces (of any dimensions) of dyadic cubes. The definition includes $C^1$ compact submanifolds of $\R^n$. Whitney sets are much more general than what we shall need for Theorem \ref{t1} or its variants.
    We want to prove the following.

    \begin{thm}\label{t2a1}
        Let $\Gamma$, $J$, $f$, $E_0$ be as above; in particular $\Gamma$ is a Whitney set and $f$ satisfies 
        \eqref{2a6} and \eqref{2a7}. Also assume that  $0 \leq m_0 < +\infty$.
        Let $\{ E_k \}$ be a minimizing sequence in $\cE(E_0,\Gamma)$ (i.e., assume that 
        $\lim_{k \to +\infty} J(E_k) = m_0$). Then there is a subsequence (which we still denote by $\{ E_k \}$)
        for which the measures $\mu_k = \H^{d}_{\vert E_k}$ converge weakly to some limit
        $\mu_\infty$. Moreover, the support $E_\infty := \mathrm{spt}(\mu)$ is sliding minimal for the functional $J$, i.e.,
        \begin{equation} \label{2a12}
            J(E_\infty) \leq J(F) \text{ for every } F \in \cE(E_\infty,\Gamma).
        \end{equation}
        and $\mu_\infty = \H^{d}_{\vert E_\infty}$. In particular, $E_\infty$ is a sliding almost minimal set, associated to the boundary $\Gamma$ and a gauge function $h$ that satisfies \eqref{2a1}.
    \end{thm}

    We highlight that $E_\infty$ is also \ub{coral}, which means
    that $E_\infty$ is the closed support of $\mu_\infty$, i.e., that $E_\infty$ is closed and $\H^d(E_\infty \cap B(x,r)) > 0$ for $x\in E_\infty$ and $r>0$. 
    Not every minimal or almost minimal set is coral, and in fact minimizers of $J$ above may be forced to have
    a part of vanishing $\H^d$-measure, which is needed because of topological reasons
    (the definition of the class $\cE$), but is forgotten when we take the weak limit of the $\mu_k$.
    Part of the job left for the next sections will be to deal with that part too.

    Theorem \ref{t2a1} is a nice way to use $E_0$ to obtain sliding minimizers, 
    but as far as the sliding Plateau problem of the introduction is concerned, we are not finished 
    yet because probably $E_\infty \notin \cE(E_0,\Gamma)$.

 In view of proving Theorem \ref{t1}, Theorem \ref{t2a1} does not help when $m_0 = 0$ because the $\mu_k$ simply converge to $0$, and we can take $E_\infty = \emptyset$. We will prove Theorem \ref{t1} in this case with a specific argument at the end of Section \ref{S3}.
    The assumption $m_0 < \infty$ does not disturb because Theorem \ref{t1} is trivial in the case $m_0 = +\infty$.

    Let us now prove Theorem \ref{t2a1}, as a direct consequence of Corollary 4.1 in \cite{La}.
    We use the elliptic integrand $J$ defined by $f$; the ellipticity comes from \eqref{2a6} and 
    the continuity of $f$. There is an open set $X$ in the statement of \cite[Corollary 4.1]{La} which plays the role of the ambient domain and which 
    we take equal to $\R^n$.
    We choose the class $\mathcal{C}$ (with the notation of \cite{La}) to be $\cE(E_0,\Gamma)$. 
    The first condition of the corollary is $m_0 < \infty$ and this is satisfied. 
    The second condition is that when $E \in \cE(E_0,\Gamma)$, any sliding deformation of $E$,
    in the sense of \cite{La}, also lies in $\cE(E_0,\Gamma)$. The definition of sliding deformations there 
    (Definition 4.1) looks a tiny bit more complicated because they can be localized in some open set $U$, and it is required that there exists a compact set $C \subset \R^n$ such that $\varphi(t,x) = x$ in $\R^n \setminus C$, for all $t$.
    However, in order to apply \mbox{\cite[Corollary 4.1]{La}}, we don't need to care about localization 
    (i.e., $U = \R^n$),  and the condition $\varphi(t,x) = x$ outside a compact set does not play any constraining role because for us: since $\Gamma$ and $E$ are compact sets, it is always possible to artificially set $\varphi(t,x) = x$ away from $\Gamma$ and $E$.
    Even if we did not keep the constraint \eqref{1a3} in our definition, the sliding deformations 
    of \cite{La} would still be sliding deformations for us.

    So we can apply Corollary 4.1 in \cite{La}, which yields that the measures $J_{\vert E_k}$ 
    converge weakly to the measure $J_{\vert E}$, where $E$ is a coral sliding minimal set for $J$. The convergence of $J_{\vert E_k}$ means that for all test functions
    $\varphi$,
    \begin{equation}
        \lim_{k \to \infty} \int_{E_k} f \varphi d\H^2 = \int_{E} f \varphi d\H^2
    \end{equation}
    but since $f$ is continuous and bounded from above and below, this is equivalent to say that $(\mu_k)$ converge weakly to $\H^d_{\vert E}$.
    We added in our statement that $E_\infty$ is also an almost minimal
    set for $\H^d$, and this follows from the discussion above (near \eqref{2a9}). 
    This concludes the proof of Theorem \ref{t2a1}.
    \qed

    \section{Good access}
    \label{S3}

    We shall continue with the construction of Section \ref{S2}, but add strong new assumptions that
    allow us to apply results from \cite{Dvv}. 
    We now assume that $d=2$, and also that 
    \begin{equation} \label{3a1}
        \Gamma \text{ is a disjoint finite union of closed curves of class $C^{1+\alpha}$}
    \end{equation}
    for some $\alpha > 0$. Since $\Gamma \subset \R^n$ is  compact, this amounts to saying that every
    point of $\Gamma$ has a neighborhood where $\Gamma$ is a $C^{1+\alpha}$ curve.

    We will also require the strange access condition \eqref{1a6}, which will be explained shortly, and whose 
    main purpose is to prevent some annoying blow-up limits of our minimal set $E_\infty$ from arising. 
    Let us first define blow-up limits, say what we want, and then state a sufficient condition 
    on $\Gamma$.

    A \ub{blow-up limit} of a closed set $E$ at a point $x\in E$ is a closed set $F$ which is the limit, in the local Hausdorff
    topology, of a sequence of sets $F_k = r_k^{-1}(E-x)$, where the $r_k$ are positive radii that tend to $0$. 
    The usual way to write the limit is to require that for every $R > 0$, 
    \begin{equation} \label{3a2}
        \lim_{k \to +\infty} d_{0,R}(F, F_k) = 0,
    \end{equation}
    where in general we set 
    for two closed sets $A, B$,
    \begin{equation} \label{3a3}
        d_{x,r}(A,B) = r^{-1} \sup_{y\in A \cap \overline{B}(x,r)} \dist(y,B) + r^{-1} \sup_{y\in B \cap \overline{B}(x,r)} \dist(y,A),
    \end{equation}
    where the supremum is considered in $[0,+\infty]$ (in particular, the supremum of the empty set is $0$). 
    The limit set in (\ref{3a2}) can then be described as
    \begin{eqnarray}
        F   &=& \set{y \in \R^n | \lim_{k \to \infty}  \dist(y, F_k) = 0}\cr
        &=& \set{y \in \R^n | \lim_{k \to \infty} \dist(x + r_k y,E) = 0}.
    \end{eqnarray}
    We observe that $F$ contains $0$ and that for all $t > 0$, the set $t F$ is also a blow-up limit of $E$ at $x$. 
    Therefore, if $E$ has a unique blow-up limit at $x$, it must be a cone.

    We intend to apply the main theorem of \cite{Dvv}, 
    which gives a good local description for $E$ when $E$ is a local coral almost minimal set 
    of dimension $2$ associated to a sliding boundary $\Gamma$ as in \eqref{3a1}. 
    Unfortunately that theorem is not general enough for the result of our dreams; it only works well
    near points $x \in E \setminus \Gamma$ or points $x \in E \cap \Gamma$ for which
    \begin{equation} \label{3a4}
        \begin{gathered}
            \text{no blow-up limit of $E$ at $x$ (here, a point of $E \cap \Gamma$)}\\
            \text{contains a type $\bY$ singularity whose spine is parallel to the tangent of $\Gamma$}. 
        \end{gathered}
    \end{equation}
    As an example, it works well at the points $x \in E \cap \Gamma$ for which
    \begin{equation}
        \begin{gathered} 
            \text{every blow-up limit of $E$ at $x$ is a half plane, a plane,}\\
            \text{a cone of type $\bV$ or a cone of type  $\bY$ or $\bT$}\\
            \text{whose spine is not parallel to the tangent of $\Gamma$ at $x$.}\\
        \end{gathered}
    \end{equation}
    In the case of $\bT$, the spine is composed of four half lines, and we demand that none 
    of these half lines is parallel to the tangent.
    We will discuss all these types later; the point is that there is an important type of
    blow-up limits of $E$ that is excluded here, the cones of type $\bY$ (three half planes
    bounded by a single line $L$, and that make angles of $2\pi/3$ along $L$), with a spine
    $L$ parallel to the tangent of $\Gamma$ at $x$. So our next condition will be designed to avoid this case.

    We shall assume that $\Gamma$ can be wrapped in a compact set $K$, in such a way that
    $\Gamma \subset K$, and all the points of $\Gamma$ have a good access to the complement
    of $K$ in the following sense. 

    \begin{defi}[Good access to the complement]\label{d3a1}
        Let a compact set $\Gamma \subset \R^n$ be a finite disjoint union of $C^{1+\alpha}$ curves. 
        Let $K \subset \R^n$ be a compact set that contains $\Gamma$. 
        We way that $\Gamma$ has a \ub{good access to the complement} of $K$ when for each point $x_0 \in \Gamma$
        and each blow-up limit $K_0$ of $K$ at $x_0$, the following happens. 
        Denote by $L_0$ the vector line parallel to the tangent line to $\Gamma$ at $x_0$ (in particular, $L_0 \subset K_0$).
        Let $e_0$ denote any of the two unit vectors of $L_0$.
        We require that for any cone $Y$ of type $\bY$ with spine $L_0$ and any choice of $c > 0$, $Y \cap B(e_0,c)$ is not contained in $K_0$.
    \end{defi}

   This means that $K_0$ cannot contain a cone showing a $\bY$ singularity with spine $L_0$.
    We allow a general compact set $K$ because we want to allow more general statements 
    than Theorem \ref{t1}, but in the case of Theorem \ref{t1}, $K$ will be chosen to be the 
    convex hull of $\Gamma$.
    As we shall see, this is because when we minimize $\H^2$, it is easy to find minimizing sequences
    that lie in the convex hull $K$, and then $E_\infty \subset K$ too.

    Let us give an example for which Definition \ref{d3a1} holds true. 
    We claim that the condition above is satisfied as soon as there exists linearly independent vectors $e_1,\ldots,e_{n-2}$ in $\R^n$ such that
    \begin{equation}
        K_0 \subset \bigcap_{i=1}^{n-2} \set{y \in \R^n | y \cdot e_i \leq 0}.
    \end{equation}
    In this case, since $L_0 \subset K_0$, we have necessarily $e_0 \cdot e_i = 0$ for all $i = 1,\ldots,n-2$.
    Given a cone $Y$ of spine $L_0$, there exists a linear plane $P$ orthogonal to $L_0$ and three unit vectors $v_1,v_2,v_3 \in P$ such that $v_1 + v_2 + v_3 = 0$ and
    \begin{equation}
        Y = \bigcup_{k=1}^3 \set{t e_0 + s v_k | t \in \R, s \geq 0}.
    \end{equation}
    Let us proceed by contradiction and assume that for some small $c > 0$, $Y \cap B(e_0,c)$ is included in $K_0$.
    This implies that for all $k = 1,2,3$ and $i = 1,\ldots,n-2$, we have $v_k \cdot e_i \leq 0$. The additional condition $\sum_k v_k = 0$ allow us to deduce that for all $k = 1,2,3$ and $i = 1,\ldots,n-2$, we have $v_k \cdot e_i = 0$. This also holds for $i = 0$ by definition of $Y$. Therefore, the three vectors $v_1, v_2, v_3$ belong to a line (the orthogonal complement of $e_0,\ldots,e_{n-2}$ in $\R^n$) and this is not possible.

    The main example of application is that the access condition holds true in dimension $n = 3$ if  $\Gamma$ is contained in the topological boundary of its convex hull $K$. 
    Indeed in this case the blow-up limits of $K$ at points of $\partial K$ are always included in a half-space.
    Thus, Definition~\ref{d3a1} allows at least the
    most classical case, but of course we are also interested in examples in dimensions $n \geq 4$.

    We want to keep the possibility to take wrapping sets $K$ that are different
    from the convex hull of $\Gamma$, but there will be some constraints, because we also want to 
    be able to choose a minimizing sequence $\{ E_k \}$ such that $E_k \subset K$ for all
    $k$, and for this the best option seems to be to  construct a retraction $R$ onto $K$, 
    which diminishes the functional $J$ (and fixes $\Gamma$, since $\Gamma \subset K$), 
    so that initial sets $E_k$ can be replaced with $R(E_k)$. When $K$ is a compact convex set and $J(E) = \H^2(E)$, 
    $R$ will be the shortest distance projection on $K$, which is $1$-Lipschitz and therefore reduces $J$. 
    One can imagine other interesting situations (we discuss this in Section \ref{Slast}), but anyway there will be strong constraints on $\Gamma$, $K$, and $f$.

    Let us not worry about this for the moment, and instead assume that we 
    have $K$ and a minimizing sequence in $K$, and use \cite{Dvv} to get the desired conclusion.
    Our main existence theorem is then as follows and will readily imply Theorem \ref{t1}.

    \begin{thm}\label{t3a1}
        Let $\Gamma$, $J$, $f$, $E_0$ be as in Section \ref{S2}, and assume in addition that
        $d=2$, \eqref{3a1} holds, $0 \leq m_0 < \infty$ and $K$ is a compact set in $\R^n$
        that contains $\Gamma$, such that $\Gamma$ has good access to the complement of
        $K$ (as in Definition \ref{d3a1}). 
        Let $\{ E_k \}$ be a minimizing sequence for $J$ in the class $\cE = \cE(E_0,\Gamma)$.
        Suppose in addition that $E_k \subset K$ for all $k$,
        or more generally that for some subsequence, the measures $\mu_k = \H^d_{\vert E_k}$ converge
        weakly to a measure $\mu_\infty$ whose support is contained in $K$.
        Then we can find a minimizer for $J$ in the class $\cE$, i.e., a set $E \in \cE$ such that 
        $J(E) = m_0$.
    \end{thm}

    As was said before, Theorem \ref{t2a1} already gives us a (sliding) minimizer $E_\infty$ for $J$,
    but $E_\infty$ may not lie in the class $\cE(E_0,\Gamma)$ because the weak limits forgot thin parts
    of the $E_k$ and  did not preserve the topology. We want to obtain our minimizer $E$
    as the image of a cleaner version of $E_k$, $k$ large, projected on $E_\infty$ by a local
    Lipschitz retraction. 

    Let us now explain how to deduce Theorem \ref{t1} from Theorem \ref{t3a1}.
    We take $J = \H^2$ and, in the main case when $0 < m_0 < \infty$, we take for $K$ the convex hull of $\Gamma$.
    Notice that if $E \in \cE(E_0,\Gamma)$ and $\pi$ is the closest point projection on $K$, 
    it is easy to see that $\pi(E) \in \cE(E_0,\Gamma)$, and since $\H^2(\pi(E)) \leq \H^2(E)$, 
    for any minimizing sequence $\{ E_k \}$ in $\cE(E_0,\Gamma)$, we can apply
    Theorem~\ref{t3a1} to the minimizing  sequence $\{ \pi(E_k) \}$.

    We promised to say a few words about the case when $m_0 = 0$. Then there are
    sets $E \in \cE(E_0,\Gamma)$ with $\H^2(E)$ arbitrarily small, and it is easy to see that
    we can use a Federer-Fleming projection with cells adapted to $\Gamma$ (so that the Federer-Fleming projection respects the sliding boundary condition) to project $E$ to a new set $F \in \cE(E_0,\Gamma)$ that does not contain a full $2$-cell, and use this fact to project again on a $1$-grid and find $G\in \cE(E_0,\Gamma)$ such that $\H^2(G) = 0$ (and even $G$ is $1$-dimensional).

    The main step of our proof will be the construction of a retraction that projects on $E_\infty$, and a first step in this direction will be a description of all the blow-up limits of $E_\infty$, which will allow us to
    use the result of \cite{Dvv}, get a nice description of $E_\infty$ near each of its points, and
    build local retractions that will then be put together.

    \section{Blow-up limits of \texorpdfstring{$E_\infty$}{Einfty} at a point} 
    \label{S4}

    In this section and the next ones, we consider the almost minimal set $E_\infty$ 
    obtained by applying Theorem \ref{t2a1} to the sequence of Theorem \ref{t3a1}, 
    and start giving a local description of $E_\infty$. Notice that by assumptions of Theorem \ref{t3a1}, the set
    $E_\infty$ is contained in $K$.

    The task of the current section is to describe the possible blow-up limits of our limit $E_\infty$ at a point $x_0 \in E_\infty$.
    Call $X$ such a blow-up limit,  and let $\{ r_k \}$ denote a sequence of radii such that
    $\lim_{k \to +\infty} r_k = 0$ and 
    \begin{equation}
        X = \lim_{k \to +\infty} r_k^{-1}(E - x)
    \end{equation}
    in local Hausdorff topology.
    We may replace $\{ r_k \}$ with a subsequence for which the $K_k = r_k^{-1} (K-x_0)$ 
    converge to a limit $K_0$, which is of course a blow-up limit of $K$ at $x_0$. 
    In the case $x_0 \in E_\infty \cap \Gamma$, we also let $L_0$ denote the unique blow-up limit of $\Gamma$ at $x_0$ (a vector line). 

    If $x_0 \in E_\infty \setminus \Gamma$, then $X$ is a plain minimal cone (centred at $0$) and if $x_0 \in E_\infty \cap \Gamma$, then $X$ is a sliding minimal cone with respect to the boundary $L_0$ (also centred at the origin).
    This comes from \cite{Dss}, but since this was done there in much more generality than needed here, 
    with some times confusing notation, let us summarize the proof.
    First, the limit $X$ of any blow-up sequence of (sliding) minimal sets is a (sliding) minimal set too; 
    this is done in Part V of \cite{Dss}; see Theorem~23.13 there. 
    In fact, this is one of the main points of the whole book. 
    Then, by the fact that the Hausdorff measure goes to the limit along locally minimizing sequences,
    the density $\theta(r) = r^{-2}\H^{2}(X \cap B(0,r))$, on balls centered at the origin, is constant 
    (because $r^{-2}\H^{2}(E_\infty \cap B(x_0,r))$ was nearly monotone and had a limit), and this forces
    $X$ to be a cone centred at $0$; see Sections 27 and 28 of \cite{Dss}. 

    In the case $x_0 \in E_\infty \cap \Gamma$, we point out that $X$ may be fully transverse to $L_0$ 
    (i.e., $X \cap L_0 = \set{0}$), or $X$ may contain only one half of $L_0$, 
   or the whole line $L_0$.
    If $X$ is not fully transverse to $L_0$, we have an additional constraint coming from the good access condition (Definition \ref{d3a1}).
    Indeed, since $E_\infty \subset K$, we infer that the blow-up limit $X$ satisfies
    \begin{equation} \label{4a1}
        X \subset K_0,
    \end{equation}
    and Definition \ref{d3a1}  prevents $X$ from containing a cone $Y$ whose spine is $L_0$.

 We are now going to present the main (sliding) 2-dimensional minimal cones in $\R^n$. Unfortunately the full classification is unknown as soon as $n \geq 4$, but we still know some rules about the structures of sliding minimal %
 cones.

    We start our description with the simpler case when $x_0 \notin \Gamma$ ($X$ is a plain minimal cone). Then we can use the
    local description of $E_\infty$ that we get from J. Taylor's theorem \cite{Ta}, or its generalization
    in \cite{Dhh, Dcc} 
    to ambient dimensions $n \geq 4$.
    Here already, we have to distinguish between two main cases, where 
    \begin{equation} \label{4a2}
        X \text{ is a minimal cone of type $\bP$, $\bY$, or $\bT$,}
    \end{equation}
    or
    \begin{equation} \label{4a3}
        X \text{ is an exotic minimal cone,}
    \end{equation}
    where by definition an exotic minimal cone is a minimal cone that does not satisfy \eqref{4a2}.
    Recall that a cone of type $\bP$ is just a plane through the origin. A cone of type $\bY$
    is the union of three half planes bounded by a same line $L$ through the origin, and that make 
    $\frac{2\pi}{3}$ angles along $L$; then $L$ is called the spine of $X$. Finally $X$
    is of type $\bT$ when $X$ is the cone over the union of the edges of a regular tetrahedron
    centered at the origin (and contained in some $3$-plane through $0$). In this case the spine
    of $X$ is the union of the four half lines through the vertices of the tetrahedron.

    When \eqref{4a2} holds, \cite{Ta} and \cite{Dcc} give a good description of 
    $E_\infty$ near $x_0$, as the image of the cone $X$ by a $C^{1+\varepsilon}$-diffeomorphism 
    (that sends $0$ to $x_0$).
    We will see later how to use this to find a retraction on $E_\infty$ defined near the origin but let us say a few words about the exotic case. For all plain $2$-dimensional cone $X$ in $\R^n$, \cite{Dhh} shows that $X \cap \d B(0,1)$ is a finite union of closed geodesic arcs (or full great circles) that can only meet at a common endpoint and any endpoint is at the junction of three arcs which make $2\pi/3$ angles. Moreover, there exists a constant $c_n > 0$ which depends only on the dimension $n$ such that the following holds. The full great circles are disjoint from all the other arcs and are even at distance $\geq c_n$ from them. All arcs have a length larger than $c_n$ and the distance between two arcs that do not have a common endpoint is always larger than $c_n$.  We do not have much more information than that. 
    The set $E_\infty$ is still locally equivalent to $X$ through a homeomorphism $\Phi$, 
    but in general we do not know for sure that $\Phi$ can be taken to be a diffeomorphism, 
    as \cite{Dhh} 
    only gives a biH\" older estimate for $\Phi$. So we will have to check that we can still build a 
    Lipschitz retraction on $E_\infty$ near $0$, constructed simply by gluing together retractions 
    defined on annuli where we have a uniform $C^{1+\varepsilon}$ control, and even we will
    use the same tools to prove that all the blow-up limits of $E_\infty$ at $0$ are biLipschitz-equivalent to each other, and locally to $E_\infty$; see Section \ref{Sbil}.

    We will return to this later, but for the moment let us continue our general description of $X$ with the more interesting case when $x_0 \in E_\infty \cap \Gamma$ ($X$ is a sliding minimal cone along $L_0$).
    We follow the same order as in Part V of \cite{Dvv} to simplify the task of the reader.
    A first possibility is that 
    \begin{equation}\label{4a4}
        X \ \text{is a half plane bounded by $L_0$ (we shall also say, a cone of type $\bH$).}
    \end{equation}
    This is the simplest possibility, also with the lowest density, and in this case
    there is a small ball centered at the origin where $E_\infty$ is equivalent to $X$,
    through a $C^{1+\varepsilon}$-diffeomorphism that maps 
    $\Gamma$ to $L_0$. 
    See Section 31 of \cite{Dvv}. 
    This situation is almost as pleasant as when \eqref{4a2} holds, and the desired retraction will be easy to construct.

    Then we consider the case when $X$ is a cone of type $\bV$ bounded by $L_0$, i.e.,
    the union of two half planes bounded by $L_0$ and that make an angle $\alpha \in [2\pi/3, \pi)$ 
    along $L_0$. When this angle is (strictly) larger than $2\pi/3$, we say that 
    \begin{equation} \label{4a5}
        X \text{ is a generic cone of type $\bV$,}
    \end{equation}
    and 
    Section 32 of \cite{Dvv} says that once again there is a small ball centered at the origin 
    where $E_\infty$ is equivalent to $X$, through a $C^{1+\varepsilon}$-diffeomorphism 
    that maps $\Gamma$ to $L_0$. Again local retractions will be easy to construct in this case.

    We excluded the degenerate case when 
    \begin{equation} \label{4a6}
        X \text{ is a plane that contains $L_0$,}
    \end{equation}
    which corresponds to $\alpha = \pi$ in the description above, because the situation
    is somewhat different there. A good $C^{1+\varepsilon}$ description of $E_\infty$
    near the origin is given in 
    Section 33 of \cite{Dvv}; this time $E_\infty$ is no longer
    equivalent to $X$, as it may leave $\Gamma$ and return to it many times (in a tangential way).
    We will have to return to this case carefully, and the construction of the local retraction will be
    unpleasant because of the sliding constraint.

    In the meantime we switch to the other extreme of $\bV$ sets, which is when 
    \begin{equation} \label{4a7}
        X \text{ is a sharp cone of type $\bV$,}
    \end{equation}
    which means that the angle of the two half planes that compose $X$ is equal to $2 \pi /3$.
    The reason why this case is a little special is that the set $E_\infty$ may split at the origin,
    with one part where $E_\infty$ behaves like a $\bV$-set, and another one where it behaves like a
    $\bY$-set with its spine away from $\Gamma$, attached to $\Gamma$ by a very short 
    piece of half plane. This is described in
    Section 34 of \cite{Dvv}; the description is not as 
    straightforward as in the generic case, say, but our retractions will still not be so much harder to construct.

    It may also happen that 
    \begin{equation} \label{4a8}
        X \text{ is fully transverse to $L_0$ at the origin,}
    \end{equation}
    which means that $X \cap \d B(0,1)$ does not meet $L_0$.
    In this case (see Section 37.1 of \cite{Dvv}), $X$ is a plain minimal cone
    (that is, with no sliding condition), and the description of $E_\infty$ near $x_0$
    is just the same as in the cases of \eqref{4a2} and \eqref{4a3} above. In fact, the sliding condition 
    does not play any serious role in the local description of $E_\infty$ or the construction of retractions, because
    $x_0$ is the only point of $E_\infty$ near $x_0$ that lies in $\Gamma$ (so the sliding condition is automatically
    satisfied as long as we keep $x_0$ fixed).

    These are the main sliding minimal cones along a line
 that will show up here, because we excluded the cones $X$ that contain a piece of $\bY$ with a spine that contains $L_0$,
 and as before don't know the full classification. 
    Thanks to Proposition 2.1 in \cite{Dvv}, we know for all sliding $2$-dimensional minimal cone $X$ along $L_0$ the set $X \cap \d B(0,1)$ is a finite union of closed geodesic arcs (or full great circles) that can only meet at a common endpoint and according to specific rules generalizing the ones of plain minimal cones. Following \cite{Dvv}, we will refer to this description as the \underline{general description of minimal cones}.

        By definition, no point of $L_0$ lies in the interior of one of the arcs (otherwise we cut the arc in two).
        With this convention, the full great circles don't meet $L_0$. The full great circles are disjoint from all the other arcs and are even at distance $\geq c_n$ from them, where $c_n$ is a positive constant that depends only on the dimension $n$. Any endpoint away from $L_0$ is at the junction of three arcs which make $2\pi/3$ angles.
        If there is an endpoint $\xi_0$ at $L_0 \cap X$, then $\xi_0$ belongs either to only one arc, or two arcs with an angle $\alpha \in [2\pi/3,\pi]$ or three arcs with $2\pi/3$ angles (a triple junction). The arcs which don't meet $L_0$ have a length larger than $c_n$. If an arc $\gamma_0$ starts from an endpoint $\xi_0 \in L_0 \cap X$ and has length $< c_n$, there either there is no other arc  
        leaving from $\xi_0$ or there is another arc $\gamma$ leaving $\xi_0$ with length $\geq c_n$ and making an angle $\geq 9 \pi/10$ with $\gamma_0$.
        Here, $\gamma_0$ ends at a triple junction $\xi_1$ at distance $< c_n$ from $\xi_0$ and two other arcs $\gamma_1$, $\gamma_2$ leaves $\xi_1$ with $2\pi/3$ angles and with length $\geq c_n$.
        If there is another arc $\gamma$ which leaves $\xi_0$ as above, then we see in particular $\gamma_1$, $\gamma_2$ are at distance $< c_n$ from $\gamma$. This is the only situation where arcs with  no common endpoints are at distance $< c_n$.
        In other words, the distance between two arcs with no common endpoints is always larger than $c_n$, except if they are connected by third arc of length $< c_n$ which has one endpoint in $L_0$.

    Using the general description above and the good access condition (Definition \ref{d3a1}), we will be able to classify the local behavior of $X \cap \partial B(0,1)$ in spherical caps of the form 
    $$
    S(\xi_0,c) := B(\xi_0,c) \cap \partial B(0,1),
    $$ 
    where $\xi_0 \in X \cap \partial B(0,1)$ and $c > 0$ may depend on $n$, $E_\infty$, $x_0$, 
     but not the chosen blow-up limit $X$ (there may a priori be an infinity of blow-up limits $X$ at $x_0$).
    It could be that we can use compactness and get even more uniformity on $c$, but we shall not try.
    We consider directly the case $\xi_0 \in X \cap L_0 \cap \d B(0,1)$, which is more interesting.
    In some cases, we could treat $\xi_0$ and $-\xi_0$ at the same time, but not in general, so we'll need to do a separate discussion near each point $\xi_0$.
    The first case is when
    \begin{equation}\label{4a9}
        \text{$X$ coincides in $S(\xi_0,c)$ with a cone of type $\bP$.}
    \end{equation}
    This is a local version of \eqref{4a6}, but we do not assume that the situation near $-\xi_0$ is the same.
    For instance, $X$ could be a set of type $\bY$ or $\bT$, and in this case $-\xi_0$ is far from $X$.
    But even if $X$ contains a full plane $P$  that contains $L_0$, it could be that it contains other pieces, such as 
    another plane $P'$ nearly orthogonal to $P$ (if $n \geq 4$), as in \cite{Li}. 
    When \eqref{4a9} holds (and we are not in the case of a plane), we will again construct 
    our retractions in intersections of sectors and annuli, and then glue them. 
    This will be a little more complicated than when \eqref{4a6} holds, but
    not fundamentally different.

    Similar to this last case, we will also be able to treat the situations where 
    \begin{equation} \label{4a10}
        \text{$X$ coincides in $S(\xi_0,c)$ with a cone of type $\bH$ or $\bV$}
    \end{equation}
    and finally, our last case is when
    \begin{equation}\label{4a12}
        \text{$X$ coincides in $S(\xi_0,c)$ with a truncated cone of type $\bY$.}
    \end{equation}
    This means that in $S(\xi_0,c)$, the set $X$ is composed of a non-empty arc $\gamma_0$ that goes from $\xi_0$ to a nearby endpoint $\xi_1 \in S(\xi_0,c)$ and then two arcs $\gamma_1, \gamma_2$ that leaves $\xi_1$ and which go all to the way to $\partial S(\xi_0,c)$ (the boundary relative to the unit sphere). The three arcs meet at $\xi_1$ with an angle of $2\pi/3$, and we have assumed that $\gamma_1$ is non-empty so as to distinguish this case from a sharp $\bV$.
    In the definition of (\ref{4a12}), we underline there are no other arcs leaving $\xi_0$ besides $\gamma_0$. The general description of minimal cones would allow such an arc but we are going to see that the good access condition and a compactness argument prevent it (when $c$ is small enough).
    By convention in \eqref{4a10}, \eqref{4a12}, and \eqref{eq_XL0} below, 
 $\bH$ and $\bV$ are bounded by $L_0$ whereas the spine of the truncated $\bY$ does not pass 
 through $\xi_0$ (otherwise, that would be a $\bV$). 

    Finally, we use the good access condition (Definition \ref{d3a1}) to prove that nothing else than the cases mentioned above can happen. This rule out complicated cases, that is, the cases where $X$ could contain 
 a $\bY$ cone of spine $L_0$, or a piece of such a cone with a significant part of $L_0$. 
    For all $x_0 \in E_\infty$, we let $\cX(x_0)$ denote the set of blow-up limits of $E_\infty$ at $x_0$ (we recall that if $x_0 \in E_\infty \setminus \Gamma$, they are plain minimal cones, and if $x_0 \in E_\infty \cap \Gamma$, they are sliding minimal cones along $L_0$, where $L_0$ is the tangent line to $\Gamma$ at $x_0$). It will be convenient to take the convention that $L_0 = \emptyset$ when $x_0 \in E_\infty \setminus \Gamma$ so as to avoid a case distinction. 
    
    \begin{lem}\label{lem_Xstructure}
        Let $x_0 \in E_\infty$, let a blow-up $X \in \cX(x_0)$.
        There exists $0 < c_* < 1$ (that depends on $n$, $E_\infty$, $x_0$ but not $X$) such that the following holds.
        \begin{equation}\label{eq_XL0}
            \begin{gathered}
                \text{For all $\xi_0 \in X \cap L_0 \cap \d B(0,1)$, the cone $X$ coincides}\\
                \text{in $S(\xi_0,c_*)$ with a cone of type $\bP$, $\bH$, $\bV$ or a truncated $\bY$.}
            \end{gathered}
        \end{equation}
        and for all $0 < c \leq c_*$,
        \begin{equation}\label{eq_XY}
            \begin{gathered}
                \text{for all $\xi \in X \cap \d B(0,1)$, if $X \cap S(\xi,10c)$ does not meet $L_0$,}\\
                \text{then $X$ coincides in $S(\xi,c)$ with a cone of type $\bP$ or $\bY$.}
            \end{gathered}
        \end{equation}
    \end{lem}

     Note that in (\ref{eq_XY}), the cone $\bY$ contains $\xi$ but its spine may possibly not pass through $\xi$.
        The important point in the Lemma is that the constant $c_*$ is uniform and works for all blow-up limits at $x_0$ (there may a priori be an infinite number of blow-up limits at $x_0$). Lemma \ref{lem_Xstructure} prevents for example $X$ from being a transverse $\bY$ cone whose spine lies too close to $L_0$.

    \begin{proof}
    We start with \eqref{eq_XL0}.
        We proceed by contradiction and assume that for all constant $\beta_k = 2^{-k}$ (so that $\beta_k \to 0$), 
 we can find a blow-up limit $X_k$ of $E_\infty$ at $x_0$ such that $X_k$ does not satisfy (\ref{eq_XL0}) in $S(\xi_0,\beta_k)$.
        Thus, for all $k$, there exists a point $\xi_0 \in X_k \cap L_0 \cap \d B(0,1)$ such that $X_k$ does not coincide with a cone of type $\bP$, $\bH$, $\bV$ or a truncated $\bY$ in $S(\xi_0,\beta_k)$.
        We may replace $\{ X_k \}$ with a subsequence for which $\xi_0$ is always the same
        (there are only two choices of $\xi_0 \in L_0$). 
        Now let us extract a subsequence (not relabelled) such that $X_k$ converges in local Hausdorff distance to a limit $X$ in $\R^n$.
        It is standard that $X$ is also a blow-up limit of $E_\infty$ at $x_0$,
        and of course $X$ still contains $\xi_0$.
        According to the general description of minimal cones, we can find a sufficiently small radius $r > 0$ such that 
    $X \cap S(\xi_0,10r)$ is composed of one, two or three arcs leaving $\xi_0$ and going all the way to 
    $\partial S(\xi_0,10r)$ (it will be also useful to assume $r < c_n/10$, where $c_n$ is the dimensional constant in the general description of minimal cones).
        If there is only one arc, then $X$ coincides in $S(\xi_0,10r)$ 
        with a cone of type $\bH$. If there are two arcs, then depending on their angle, $X$ coincides in $S(\xi_0,10r)$ with a cone of type $\bP$ or $\bV$ (sharp or generic). Finally, if there are three arcs, then $X$ coincides in $S(\xi_0,10r)$ with a cone of type $\bY$ of spine $L_0$.
        However, the inclusion $E_\infty \subset K$ implies the inclusion of the blow-ups 
        $X \subset K_0$ and the access condition (Definition \ref{d3a1}) means that 
        $K_0$ cannot contain a significant piece of 
        $\bY$ with spine $L_0$, as in this last case. 
        Thus, $X$ must coincide in $S(\xi_0,10r)$ with a cone of type $\bP$, $\bH$, or $\bV$.

        Assuming $k$ large enough, the Hausdorff distance
        \begin{equation}\label{eq_distanceXXk}
            \sup_{x \in S(\xi_0,10r) \cap X} \mathrm{dist}(x,X_k) + \sup_{x \in S(\xi_0,10r) \cap X_k} \mathrm{dist}(x,X)
        \end{equation}
        between $X$ and $X_k$ in $S(\xi_0,10r)$ is as small as we want.
        We then justify that $k$ big enough, $X_k$ must coincide in $S(\xi_0,r)$ with a cone of type $\bP$, $\bH$, $\bV$ or a truncated $\bY$. Taking $k$ also big enough so that $\beta_k = 2^{-k} < r/2$, this will yield a contradiction and prove (\ref{eq_XL0}).

            Our first case is when $X$ coincide with a $\bP$ cone in $S(\xi_0,10r)$. We are going to see that this only leaves for $X_k$ the possibility to coincide with a $\bP$ or a $\bV$ cone (with an angle very close to $\pi$) in $S(\xi_0,r/2)$.
            Taking $k$ big enough, there exists a 2-dimensional plane $P$ such that $X_k \cap S(\xi_0,5r)$ lie in an arbitrarily thin neighborhood of $P$. Remember that $\xi_0 \in X_k \cap L_0$ so, by the general description of minimal cones, $X_k$ may contain one, two or three arcs leaving $\xi_0$. We can already see that three arcs is not possible: as $r < c_n/10$, a triple junction at $\xi_0$ would go all the way to $\partial S(x_0,5r)$ with such angle conditions that it would escape any small neighborhood of $P$. In order to deal with the two other cases, we should rule out the possibility that $X_k$ contains at least one arc $\gamma_0$ leaving $\xi_0$ and an other arc $\gamma_1$ passing through $S(x_0,r)$, but which does not have $\xi_0$ as an endpoint. We justify by contradiction that this cannot happen. 
            First, $\gamma_0$ starts from $\xi_0$ but cannot end in $S(\xi_0,r)$, otherwise that would create a triple junction leaving a small neighborhood of $P$; 
        therefore $\gamma_0$ has length $\geq r$. As $\gamma_1$ does not meet $\xi_0$ and is a geodesic arc, it does not meet $-\xi_0$ either. Thus, it does not meet $L_0$ and the general description states that it has length $\geq c_n \geq r$.
            Now, as $\gamma_0$ and $\gamma_1$ are geodesic arcs passing through $S(\xi_0,r)$, which have length $\geq r$ and whose intersection with $S(\xi_0,5r)$ lie in a small neighborhood of $P$, we deduce that their corresponding great circles are also entirely contained in a small neighborhbood of $P$ (as thin as we want provided we take $k$ big enough).
            Using again $r < c_n$, the arcs $\gamma_0$ and $\gamma$ are at distance $< c_n$ from each other, so they either have a common endpoint and make an angle $2\pi/3$ or they are connected by a third arc $\gamma_3$ with an endpoint in $L_0$,  which makes a $\geq 9\pi/10$ angle with $\gamma$ and a $2\pi/3$ angle with $\gamma_1$. Because of the angle conditions, both cases would contradict the fact that the great circles corresponding to $\gamma_0$ and $\gamma_1$ lie in a thin neighborhood of $P$.
            We conclude that all arcs in $X_k$ which meet $S(\xi_0,r)$ have $\xi_0$ as an endpoint, and our claim follows easily.

            In the case where $X$ coincide with a $\bH$ cone in $S(x_0,10r)$, one can see in the same spirit as before that for $k$ big enough, $X_k$ admits only one arc leaving $\xi_0$ and that no arc can meet $S(\xi_0,r)$ without having $\xi_0$ as endpoint. Therefore, $X_k$ coincides with a $\bH$ cone in $S(\xi_0,r)$.
       
     Similarly, if $X$ coincides with a $\bV$ cone in $S(x_0,r)$, then for $k$ big enough, $X_k$ coincides with  a $\bV$ cone or a truncated $\bY$ cone in $S(\xi_0,r)$ (a truncated $\bY$ with a very short arc looks like a $\bV$ in Hausdorff distance).
 This completes the proof of (\ref{eq_XL0}).
 
        We pass to \eqref{eq_XY}, but we justify a few intermediate properties first.
        Let $\gamma_1, \gamma_2$ be two arcs from the decomposition of $X$ and assume that they have a common endpoint $\zeta \in X \cap \d B(0,1) \setminus L_0$. We observe that for any $c > 0$, if a point $\xi \in \gamma_1$ is at distance $\geq c$ from $\zeta$, then $S(\xi,c/4)$ is disjoint from $\gamma_2$. This comes from the fact that when $\gamma_1$ and $\gamma_2$ meet, they make an angle $2 \pi/3$ which ensures that $\gamma_2$ cannot pass through $S(\xi,c/4)$.

            We are going to deduce that for $0 < c < c_n/2$, if a triple junction $\xi \in X \cap \d B(0,1)$ is such that $X \cap S(\xi,c)$ does not meet $L_0$, then $X$ coincides in $S(\xi,c/4)$ with a $\bY$ cone.
            Let $\gamma_1, \gamma_2, \gamma_3$ be the three arcs leaving $\xi$. For each one of these arcs, either they don't meet $L_0$ and they have a length $\geq c_n$, or they meet $L_0$ but in both cases, they go all the way to $\partial S(\xi,c)$. We then justify by contradiction that no other arc $\gamma$ meet $S(\xi,c/4)$. First, we show that if such an arc $\gamma$ exists, it cannot have a common endpoint with one of the $\gamma_i$. Indeed, say that $\gamma$ has a common endpoint with $\gamma_1$. The endpoint cannot be $\xi$ (as arcs have disjoint interiors, this would make a quadruple junction) and since $\gamma_1$ goes all the way to $\partial S(\xi,c)$, the endpoint is outside of $S(\xi,c)$. It follows from the previous paragraph that $\gamma$ is disjoint from $S(\xi,c/4)$  and we reach a contradiction. 
            Now, $\gamma$ is at distance $< c_n$ from each $\gamma_i$ and has  no common endpoint with them, so it must be connected to each of them by a short arc of length $< c_n$ with endpoint in $L_0$. However, the set $X \cap L_0 \cap \d B(0,1)$ can have at most two points and each point $\xi_0 \in X \cap L_0 \cap \d B(0,1)$ can belong to at most one short arc of length $< c_n$. Again, this situation would create a quadruple junction.

            Now it is easier to deduce (\ref{eq_XY}). Let $\xi \in X \cap \d B(0,1)$ be such that $X \cap S(\xi,10c)$ does not meet $L_0$. If $X$ contains no triple junction in $S(\xi,c)$, then it coincides with a cone of type $\bP$ in $S(\xi,c)$. If $X$ has a triple junction at $\xi_1 \in S(\xi,c)$, then $X \cap S(\xi_1,8c)$ does not meet $L_0$ so $X$ coincides in $S(\xi_1,2c) \supset S(\xi,c)$ with a cone of type $\bY$ whose spine passes through $\xi_1$.
    \end{proof}

    In the next Lemma, we take the convention that $L_0 = \emptyset$ when $x_0 \in E_\infty \setminus \Gamma$ in order to avoid a case distinction.

    \begin{lem}[Covering of $X \cap \partial B(0,1)$]\label{lem_Xcover}
        Let $x_0 \in E_\infty$, let $X \in \cX(x_0)$ and let $0 < \nu \leq 1/2$.
        There exists $0 < c_* < 1$ (that depends on $n$, $E_\infty$, 
        $x_0$ but not $X$ or $\nu$) 
        such that for all $0 < c \leq c_*$, there exists a family of spherical caps $S(\xi_i,c_i)$ where $\xi_i \in X \cap S$, $c_i \in [10^{-7} \nu c,c]$ such that $X \cap S$ is covered by the caps $S(\xi_i,5c_i)$ and for all $i$, and the following description holds true: 
        \begin{equation}
            \begin{gathered}
                \text{if $\xi_i \notin L_0$, then $X$ coincides in $S(\xi_i,10^3 c_i)$}\\
                \text{with a cone of type $\bP$ or $\bY$ whose spine passes through $\xi_i$}
            \end{gathered}
        \end{equation}
        and
        \begin{equation}
            \begin{gathered}
                \text{$\xi_i \in L_0$, then $X$ coincides in $S(\xi_i,10^3 c_i)$}\\
                \text{with a cone of type $\bP$, $\bV$, $\bH$ or}\\
                \text{a truncated $\bY$ whose spines passes through $S(\xi_i, \nu c_i)$.}
            \end{gathered}
        \end{equation}
        We can build the covering in such a way that each radius $c_i$ depends only on $c$, $\nu$ and the type of cone which describes $X$ in $S(\xi_i,10^3 c_i)$. Moreover,
        \begin{enumerate}
            \item the constant $c_i$ is the same for all $\bV$, truncated $\bY$ (near $L_0$) and full $\bY$ (away from $L_0$); 
            \item the spherical caps $S(\xi_i,10^3 c_i)$ where $X$ is a $\bH$, $\bV$, $\bY$ or a truncated $\bY$ are disjoint;
            \item if $X$ is a $\bP$ in $S(\xi_k, 10^3 c_k)$ and a $\bH$, $\bV$, $\bY$ or a truncated $\bY$ in $S(\xi_i, 10^3 c_i)$, then $\xi_k \notin S(\xi_i, 5c_i)$ and $c_k \leq c_i/10$.
        \end{enumerate}
    \end{lem}

    We will use such coverings to build a local Lipschitz retraction onto $E_\infty$ by gluing local formulas. In order for the formulas to glue well, the construction will be intrinsic (independent of the chosen covering).

    \begin{proof}
        Let $c > 0$ be small enough so that $10^3 c \leq c_*$, where $c_*$ is the constant of Lemma \ref{lem_Xstructure}.
        We deal directly with the case $x_0 \in E_\infty \cap \Gamma$, which is more difficult.
        We start by looking at what happens at points $\xi_0 \in X \cap L_0 \cap \d B(0,1)$ 
        (there might be none, one, or two such points).  
        It follows from (\ref{eq_XL0}) in Lemma \ref{lem_Xstructure} that for all $\xi_0 \in X \cap L_0 \cap \d B(0,1)$, either
        \begin{equation}\label{eq_XL0a}
            \begin{gathered}
                \text{the set $X$ coincides in $S(\xi_0, 10^3 c)$ with a cone of type $\bP$, $\bV$ or}\\
                \text{a truncated $\bY$ whose spine passes through $S(\xi_0,\nu c)$}
            \end{gathered}
        \end{equation}
        or
        \begin{equation}\label{eq_XL0b}
            \begin{gathered}
                \text{the set $X$ coincides in $S(\xi_0, \nu c)$ with a cone of type $\bH$.}
            \end{gathered}
        \end{equation}
        When (\ref{eq_XL0a}) holds true, we add $S(\xi_0,c)$ in our family of spherical caps and when (\ref{eq_XL0b}) holds true, we add $S(\xi_0, 10^{-4} \nu c)$. 
        Even if $X \cap L_0 \cap \d B(0,1)$ is composed of two points $\set{\pm \xi_0}$, we can assume $c_*$ sufficiently small so that all spherical caps introduced so far are disjoint.

        Now, we are looking at the points $\xi \in X \cap \d B(0,1) \setminus L_0$. At such a point $\xi$, the set $X$ must be an arc or a triple junction.
        We are first going to add the triple junctions to our family of spherical caps and make sure that the caps $S(\xi_i,10^3 c_i)$ remain disjoint.
        So, let $\xi \in X \cap \d B(0,1) \setminus L_0$ be a triple junction. We distinguish between several cases.

        In the first case, we assume that there exists $\xi_0 \in X \cap L_0 \cap \d B(0,1)$ such that $\abs{\xi - \xi_0} < 10^3 c$. According to (\ref{eq_XL0}), the only possibility is that
        \begin{equation}\label{eq_XLY}
            \text{$X$ coincides with a truncated $\bY$ whose spine passes through $\xi$ in $S(\xi_0,10^3 c)$}.
        \end{equation}
        If $\abs{\xi - \xi_0} < \nu c$, then $S(\xi_0,c)$ is one of the spherical cap introduced above (corresponding to the (\ref{eq_XL0a}) case) and it covers $\xi$ so we don’t need to do anything.
        If on the other hand $\abs{\xi - \xi_0} \geq \nu c$, then $X \cap S(\xi, \nu c)$ is disjoint from $L_0$ so by (\ref{eq_XY}), $X$ coincides in $S(\xi,10^{-1} \nu c)$ with a cone of type $\bY$ whose spine passes though $\xi$.
        We add the spherical cap $S(\xi,10^{-4} \nu c)$ to our family.
        Note that in this case, $X$ coincides in $S(\xi_0, \nu c)$ with a $\bH$ cone. Therefore, $S(\xi_0,10^{-4} \nu c)$ is one of the spherical caps introduced above in our family but since $\abs{\xi - \xi_0} \geq \nu c$, we see that $S(\xi_i, 10^{-1} \nu c)$ and $S(\xi_0,10^{-1} \nu c)$ are disjoint. 
        In case we also have $-\xi_0 \in X$, the spherical cap $S(\xi,10^{-4} \nu c)$ is still disjoint from 
        any spherical caps centered %
        at $-\xi_0$ because we assume $c_*$ small enough.
      Because of (\ref{eq_XLY}), there are no other triple junctions centred at $S(\xi_0,10^3 c)$. There could be a triple junction centred a point in $S(-\xi_0,10^3 c)$ but we can assume $c_*$ sufficiently far so that all spherical caps introduced so far are disjoint.

        Now, we focus on the case where $S(\xi, 10^3 c)$ is disjoint from $L_0$. Then (\ref{eq_XY}) tells us that $X$ coincides in $B(\xi, 100 c)$ with a cone of type $\bY$. In particular, note that $\xi$ is at distance $\geq 100 c$ from all other triple junctions.
        We then add the spherical cap $S(\xi, 10^{-4} \nu c)$ to our family (we want all spherical caps $S(\xi_i,c_i)$ centred at a triple junction in $X \cap \d B(0,1) \setminus L_0$ to have the same radius $c_i$ and the last paragraph forces us to set $c_i = 10^{-4} \nu c$).
        It is clear that $S(\xi, 10^{-4} \nu c)$ is disjoint from all other spherical caps of our family so far, whether they are centred at a point $\xi_0 \in X \cap L_0 \cap \d B(0,1)$ or at a triple junction in $X \cap \d B(0,1) \setminus L_0$.

        We are left with the points $\xi \in X \cap \d B(0,1) \setminus L_0$ where $X$ is not a triple junction.
        If there exists already a spherical cap $S(\xi_i,c_i)$ in our family such that $\abs{\xi - \xi_i} < 5c_i$, there is nothing to do. Otherwise, $\xi$ must be at distance $\geq 10^{-4} \nu c$ from $X \cap L_0 \cap \d B(0,1)$ and from any triple junctions $\xi \in X \cap \d B(0,1) \setminus L_0$.
        Therefore, (\ref{eq_XY}) shows that $X$ must coincide in $S(\xi, 10^{-4} \nu c)$ with a cone on type $\bP$ and we add $S(\xi,10^{-7} \nu c)$ to our family. Finally, we use the compactness of $X \cap \d B(0,1)$ to extract a finite subcover.
    \end{proof}

    In the next Lemma, we choose, for each $x_0 \in E_\infty$, a small radius $r_0 = r_0(x_0) < 1$ so that 
    for $0 < r \leq r_0$, $E_\infty$ is well approximated in $B(x_0, 10r)$ by a blow-up limit $X_r \in \cX(x_0)$.

    \begin{lem}\label{l4a1}
        Let $x_0 \in E_\infty$ be as above. Then for each $\varepsilon > 0$, we can find $r_0 \in (0,1)$
        such that for $0 < r \leq r_0$, there is a blow-up limit $X_r \in \cX(x_0)$ such that
        \begin{equation}\label{4a21}
            d_{0,200r}(E_\infty - x_0, X_r) \leq \varepsilon
        \end{equation}
        with the notation of \eqref{3a3}. 
    \end{lem}

    We are forced to let $r_0$ depend on $x_0$ too, because for instance $x_0$ could be a
    regular point of $E_\infty$ (where $E_\infty$ has a tangent plane) that lies very close 
    to a point of $E_\infty \cap \Gamma$ of type $\bV$.
    In this case, $\cX(x_0)$ contains only planes but $E_\infty$ can be very close to a $\bV$-set at scales 
    larger than $\dist(x_0,\Gamma)$, which can be arbitrarily small.
    Similarly, in the good cases $E_\infty$ has a unique blow up limit $X$ at $x_0$, and then 
    \eqref{4a21} holds with this $X$, but in general we do not know whether this is the case, so we 
    have to allow $X_r$ to depend on $r$.

    \begin{proof}
        For the proof, suppose the lemma fails, so that for some choice of $x_0$ and $\varepsilon$,
        we can find a sequence $(r_k) \to 0$ such that for all $k$ and for all $X \in \cX(x_0)$,
        \begin{equation}
            d_{0,200 r_k}(E_\infty - x_0,X) \geq \varepsilon.
        \end{equation}
        Since $X$ is a cone, this is equivalent to saying
        that $d_{0,200}(F_k,X) \geq \varepsilon$, where $F_k = r_k^{-1}(E_\infty - x_0)$.
        We can find a subsequence (that we still denote by $\{ r_k \}$) such that the sets
        $F_k$ converge to some limit $X$. By definition, $X$ is a 
        blow-up limit of $E_\infty$ at $x_0$, so $X \in \cX(x_0)$.
        But now $X$ satisfies $d_{0,200}(F_k,X) < \varepsilon$ for $k$ big enough because $X$ is the limit of the $F_k$. The lemma follows from this contradiction.
    \end{proof}

    In conclusion, we have the full list of possible behaviors of $X \cap \partial B(0,1)$
    and \cite{Dvv} will now give us enough information on $E_\infty$ to construct local retractions 
    on $E_\infty$.

    \section{A description of \texorpdfstring{$E_\infty$}{Einfty} in annuli centered at \texorpdfstring{$x_0 = 0$}{x0 = 0}} 
    \label{S5}

    In the discussion that follows, the point $x_0 \in E_\infty$ is fixed, and we take $x_0 = 0$ to simplify the notation. 
    For $\varepsilon_0 > 0$, as small as we want, we choose $r_0 \in (0,1)$ as in Lemma \ref{l4a1} 
    so that for all $0 < r \leq r_0$, there exists a blow-up limit $X = X_r \in \cX(0)$ such that
    \begin{equation}\label{5a1}
        d_{0,200 r}(E_\infty,X) \leq \varepsilon_0/200,
    \end{equation}
    i.e.,
    \begin{equation}\label{5a2}
        \begin{gathered}
            E_\infty \cap \overline{B}(0,200 r) \subset \set{x \in \R^n | \mathrm{dist}(x,X) \leq \varepsilon_0 r}\\
            X \cap \overline{B}(0,200 r) \subset \set{x \in \R^n | \mathrm{dist}(x,E_\infty) \leq \varepsilon_0 r}.
        \end{gathered}
    \end{equation}

    This will give us a good enough parametric description of $E_\infty$ in the open annuli $A_r$, where
    \begin{equation}
        A_r := B(0,100 r) \sm \overline{B}(0,r/100), \quad  0 < r \leq r_0.
    \end{equation}
    and more generally for $\lambda \geq 2$, we let
    \begin{equation} \label{5b1}
        A_r(\lambda) := B(0,\lambda r) \sm \overline{B}(0,r/\lambda ), \quad  0 < r \leq r_0
    \end{equation}
    so that $A_r$ is an abbreviation for $A_r(100)$.
    In the good cases, for instance when $X$ is a sharp $\bV$-cone, we could even get such a description directly
    in the whole ball $B(0,100r)$, but $X=X_r$ may also be one of the exotic minimal cones, and then we will only
    be able to get a good control on $A_r$ (with the need to do gluing arguments later). 

    For $\xi_0 \in X \cap \partial B(0,1)$ and $0 < c < 1$, we introduce the conical open domain
    \begin{equation} \label{5b2}
        H(c,\xi_0) := \set{x \in \R^n \setminus \set{0} | \dist(\abs{x}^{-1} x, \xi_0) < c}.
    \end{equation}

    We fix a constant $c \leq 10^{-3} c_*$, where $c_*$ is given by Lemma \ref{lem_Xcover}.
    The constant $c$ is allowed to depend on $n$, $E_\infty$, $x_0$, but neither on $r$ nor on the choice of $X$.
    We fix a small $0 < \nu \leq 10^{-2}$ which will be a universal constant.
    The value $\nu = 10^{-2}$ should be good enough for most of the paper but we prefer to have more flexibility in a few places.
    We also fix a constant $\tau > 0$ which is small and is allowed to depend on $n$, $\nu$ and $c$.
    The constant $\varepsilon_0$ is as small as we want depending on $n$, $E_\infty$, $x_0$, $\nu$, $c$ and $\tau$. And finally, $r_0$ depends on everything else.

    For $0 < r \leq r_0$, we are going to cover $E_\infty \cap A_r(2)$ by boxes of the form $A_r(2) \cap H(10c_i,\xi_i)$, where $\xi_i \in X \cap \partial B(0,1)$, $c_i \in [10^{-3}\nu c,c]$ and such that we have an explicit description of $E_\infty$ in the larger domain $A_r \cap H(100 c_i, \xi_i)$.
    For this purpose, we apply Lemma \ref{lem_Xcover} (with $\nu/10$ instead of $\nu$) to cover $X \cap \partial B(0,1)$ with spherical caps $S(\xi_i, 5c_i)$, where $\xi_i \in X \cap B(0,1)$, $c_i \in [10^{-8} \nu c,c]$ and such that
    \begin{equation}\label{eq_Xcover}
        \begin{gathered}
            \text{$X$ coincides in $S(\xi_i, 10^3 c_i)$,}\\
            \text{with a cone $Z_i$ which is a $\bP$, a $\bV$, a $\bH$, a $\bY$ whose spine passes through $\xi_i$}\\
            \text{or a truncated $\bY$ whose spine passes through $B(\xi_i,\nu c_i/10)$.}
        \end{gathered}
    \end{equation}
    Moreover, the family of spherical caps satisfy the two additional properties at the end of Lemma \ref{lem_Xcover}.
    Let us now justify that $E_\infty \cap A_r(2)$ is covered by the boxes $A_r(2) \cap H(10c_i,\xi_i)$. 
        According to (\ref{5a2}), for all $x \in E_\infty \cap A_r(2)$, there exists $\xi \in X$ such that $\abs{x - \xi} \leq \varepsilon_0 r$ and in particular
    \begin{subequations}
        \begin{align}
            \mathrm{dist}(\abs{x}^{-1}x, \abs{\xi}^{-1} \xi) &\leq \abs{\frac{x}{\abs{x}} - \frac{\xi}{\abs{x}}} + \abs{\frac{\xi}{\abs{x}} - \frac{\xi}{\abs{\xi}}}\label{eq_Ecover}\\
            &\leq 2 \abs{x}^{-1} \abs{x - \xi} \leq 20 \varepsilon_0 r\nonumber.
        \end{align}
    \end{subequations}
    As $\abs{\xi}^{-1} \xi \in X \cap \partial B(0,1)$, it belongs to one of the spherical cap $S(\xi_i,5c_i)$ and thus, for $\varepsilon_0$ small enough (compared to $c$ and $\nu$), the point $x$ belongs to $A_r(2) \cap H(10c_i,\xi)$.

    Observe that given a pair $(c_i,\xi_i)$ and a cone $Z_i$ as in (\ref{eq_Xcover}) and if $\varepsilon_0$ is small enough (compared to $c$ and $\nu$), (\ref{5a2}) implies that
    \begin{equation}
        \begin{aligned}
            \dist(x,Z_i) 
            \leq \varepsilon_0 r \ \text{ for } x \in E_\infty \cap A_r(200) \cap H(200 c_i,\xi_i),\\
            \dist(x,E_\infty) \leq \varepsilon_0 r \ \text{ for } x \in Z_i \cap A_r(200) \cap H(200 c_i,\xi_i).
        \end{aligned}
    \end{equation}
Indeed, for $x \in E_\infty \cap A_r(200) \cap H(200 c_i,\xi_i)$, property (\ref{5a2}) shows there exists $\xi \in X$ such that $\abs{x - \xi} \leq \varepsilon_0 r$ and as $$\mathrm{dist}(\abs{x}^{-1}x, \abs{\xi}^{-1} \xi) \leq 2 \abs{x}^{-1} \abs{x - \xi} \leq 400 \varepsilon_0,$$ we have $\abs{\xi}^{-1} \xi \in S(\xi_i, 10^3 c_i)$ for $\varepsilon_0$ is small enough. In particular, $\xi \in Z_i$ so $\mathrm{dist}(x, Z_i) \leq \varepsilon_0 r$.
      Conversely,  for $x \in Z_i \cap A_r(200) \cap H(200,\xi_i)$, we have in fact $x \in X$ so $\mathrm{dist}(x,E_\infty) \leq \varepsilon_0 r$ holds directly by (\ref{5a2}).

    At this point, we can apply \cite{Dvv} to describe $E_\infty$ in $A_r \cap H(100 c_i,\xi_i)$, depending on the different cases for $Z_i$. To lighten the notation, we will write $Z$ for $Z_i$, $\xi_0$ for $\xi_i$, 
    and $c$ for $c_i$ in all cases.

    Let us go directly to the most interesting case where 
    $x_0 \in E_\infty \cap \Gamma$ and $\xi_0 \in X \cap L_0 \cap \d B(0,1)$.
    Then, $Z$ can be of type $\bP$, $\bV$ (sharp or generic), $\bH$,
    or can be a truncated $\bY$ whose spine passes through $B(\xi_0,\nu c/10)$.

    \msi{\bf Case 1.}
    Our first challenging case is when $Z$ is a plane $P$ containing $\xi_0 \in L_0$; we have
    \begin{equation} \label{5b3}
        \begin{aligned}
            \dist(x,P) \leq \varepsilon_0 r \ \text{ for } x \in E_\infty \cap A_r(200) \cap H(200 c,\xi_0),\\
            \dist(x,E_\infty) \leq \varepsilon_0 r \ \text{ for } x \in P \cap A_r(200) \cap H(200 c,\xi_0).
        \end{aligned}
    \end{equation}
    When \eqref{5b3} holds, with $\varepsilon_0$ small enough, we can apply
    Theorem 33.1 of \cite{Dvv} and get that there is a $\tau$-Lipschitz function $\varphi : P \to P^\perp$ 
    such that $\abs{\varphi} \leq C \tau r$ everywhere and 
    \begin{equation} \label{5b5}
        \text{$E_\infty$ coincides with the graph of $\varphi$ in $A_r \cap H(100 c,\xi_0)$.}
    \end{equation}
    Here $\tau > 0$ is as small as we want (provided that we take $\varepsilon_0$ accordingly small) and it is allowed to depend on $n$ and $c$.
    Theorem~33.1 of \cite{Dvv} actually gives more information on $\varphi$ and the structure
    of $E_\infty \cap \Gamma$ in $A_r$, which we may recall later when needed. This looks good, except for
    that fact that $\Gamma$ is allowed to leave $E_\infty$ (tangentially only), and then return to it, in possibly 
    complicated ways; this will potentially cause trouble with the sliding condition when we try to move points
    along $E_\infty$.

    \msi{\bf Case 2.}
    Our second more challenging case is when $Z$ is a sharp cone $V$ of type $\bV$; we have
    \begin{equation} \label{5b6}
        \begin{gathered}
            \dist(x,V) \leq \varepsilon_0 r \ \text{ for } x \in E_\infty \cap A_r(200) \cap H(200 c,\xi_0),\\
            \dist(x,E_\infty) \leq \varepsilon_0 r \ \text{ for } x \in V \cap A_r(200) \cap H(200 c,\xi_0).
        \end{gathered}
    \end{equation}
    In this case we appeal to Theorem 34.1 of \cite{Dvv} (applied to a finite number of centers that lie on $\Gamma \cap A_r$) to get the following description of $E_\infty$ in $A_r \cap H(100 c,\xi_0)$.
    We start with a bit of notation. Write $x\in \R^n$ as $x = (x_1, x_2, x_3, x_4)$, 
    where $x_1, x_2, x_3 \in \R$ and $x_4 \in \R^{n-3}$ (this coordinate exists only when $n \geq 4$). Assume to simplify that $\xi_0 = (1,0,0,0)$ 
and
    \begin{equation} \label{5b7}
        V =  \set{(x_1, x_2, x_3, x_4) \in \R^n\, ; \,  x_4 = 0 \text{ and } x_3 = \abs{x_2}/{\sqrt3}}.
    \end{equation}
    Since $\Gamma$ is smooth (and if $r_0$ is small enough, depending on $\Gamma$), there exists some function $C^1$ function $\psi^0 : \R \to \R^{n-1}$,
    \begin{equation}
        \psi^0 \colon s \mapsto (\psi^0_2(s), \psi^0_3(s),\psi^0_4(s))
    \end{equation} 
    with $\psi^0_2(s), \psi^0_3(s) \in \R$, and 
    $\psi^0_4(s) \in \R^{n-3}$ such that $\psi^0$ is $\tau$-Lipschitz, $\abs{\psi^0} \leq C \tau r$ and
    \begin{equation} \label{5b8}
        \Gamma \cap B(0, 100 r) =  \set{(s, \psi^0_2(s), \psi^0_3(s), \psi^0_4(s)) | s \in \R} \cap B(0, 100 r).
    \end{equation}
        For $x \in \Gamma \cap A_r$ such that $x_1 \geq 0$, we have $\abs{x - x_1 \xi_0} \leq C \tau r$ so
        \begin{equation}\label{eq_Gammatau}
            \mathrm{dist}(\abs{x}^{-1}x,\xi_0) \leq 2 \abs{x}^{-1} \abs{x - x_1 \xi_0} \leq C \tau.
        \end{equation}
        We can thus assume $\tau$ small enough so that 
    $\Gamma \cap A_r \cap H(100c,\xi_0) \subset H(\nu c,\xi_0)$.

\begin{figure}[!h]
  \centering
    \includegraphics[width=6cm]{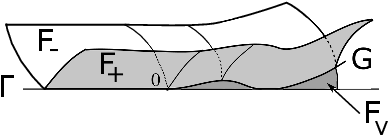}
\vskip-0.5cm
 \caption{The set $E_\infty$ in Case 2}
\label{fig1}
\end{figure}

    Next, if $\varepsilon_0$ is small enough (depending on $\tau$), there is a relatively closed curve 
    $G \subset A_r \cap H(100 c,\xi_0)$, for which there also exists some $\tau$-Lipschitz function 
    $\psi : \R \to \R^{n-1}$ such that $\abs{\psi} \leq C \tau r$ and
    \begin{equation}\label{5b9}
        G = \set{(s, \psi_2(s), \psi_3(s), \psi_4(s)) | s \in \R} \cap A_r \cap H(100c,\xi_0)
    \end{equation}
    and in addition $G$ lies roughly above $\Gamma$ in the sense that
    \begin{equation} \label{5b10}
        \psi_3(x)-\psi^0_3(x)\ \geq 0 \ \ \text{and} \ \
        |\psi_2(x)-\psi_2^0(x)| + |\psi_4(x)-\psi_4^0(x)| \leq \tau \left(\psi_3(x)-\psi^0_3(x)\right).
    \end{equation}
    From the bound $\abs{\psi} \leq C \tau r$, we know that
    \begin{equation}
        \text{$\mathrm{dist}(\abs{x}^{-1}x,\xi_0) \leq C \tau$ for $x \in G$,} 
    \end{equation}
    and we assume $\tau$ small enough so that $G \subset H(\nu c,\xi_0)$.
    The curve $G$ may coincide with $\Gamma$ in some places, in fact it is also $C^{1+\varepsilon}$,
    and when it touches $\Gamma$ this happens tangentially. 
    Finally, there are three relatively closed faces $F_v, F_+, F_- \subset A_r \cap H(100c,\xi_0)$ so that 
    \begin{equation} \label{5b11}
        E_\infty \cap A_r \cap H(100 c,\xi_0) = F_v \cup F_+ \cup F_-,
    \end{equation}
    and with the following rough description (see Figure \ref{fig1}, and \cite{Dvv} for more details),
  
    First, $F_v$ is a piece of graph $\set{z + \varphi_v(z)}$ over the vertical plane $P_v := \set{x_2=0 \text{ and } x_4=0}$ (generated by $\xi_0 = (1,0,0,0)$ and $\xi_v = (0,0,1,0)$)
 of a $\tau$-Lipschitz function $\varphi_v : P_v \cap B(0,100r) \to P_v^\perp$ such that $\abs{\varphi_v} \leq C \tau r$.
This piece of graph is bounded by the two curves $\Gamma$ and $G$, that is,
        \begin{equation}
            F_v = \set{z + \varphi_v(z) | z \in P_v \ \text{such that} \ \psi^0_3(z_1) \leq z_3 \leq \psi_3(z_1)} \cap A_r
        \end{equation}
        and when $z_3 = \psi_3(z_1)$, we have $z + \varphi_v(z) = (z_1,\psi(z_1))$ (resp. the same with $\psi^0$ instead of $\psi$).

    It could be essentially empty (i.e., $F_v = \Gamma \cap A_r \cap H(100c,\xi_0)$), 
    or composed of one or many tiny pieces (depending on how and when $G$ separates 
    from $\Gamma$), and we see it as a thin (roughly) vertical wall that connects $\Gamma$ to $G$ above $\Gamma$.

        Assuming as usual $\tau$ small enough, the vertical piece $F_v$ is entirely contained in $H(\nu c,\xi_0)$. 
        Indeed, it is clear from the above description that for any point $x = z + \varphi_v(z) \in F_v$, we have $\abs{z_3}, \abs{\varphi_v(z)} \leq C \tau r$ so $\abs{x - z_1 \xi_0} \leq C \tau r$, and thus since $\abs{x} \geq r/100$,
        \begin{equation}
            \mathrm{dist}(\abs{x}^{-1} x, \xi_0) \leq 2 \abs{x}^{-1} \abs{x - z_1 \xi_0} \leq C \tau.
        \end{equation}

    The two main pieces $F_+$ and $F_-$ look like wings that are attached to $G$.
    More precisely, $F_\pm$ is a piece of graph $\set{z + \varphi_\pm(z)}$ over the plane
$P_\pm = \set{x_3 = \pm x_2/{\sqrt3}  \text{ and } x_4 = 0}$ (generated by $\xi_0 = (1,0,0,0)$ and $\xi_{\pm} = (0,\pm \sqrt{3}/2, 1/2,0)$) of a $\tau$-Lipschitz function $\varphi_\pm : P_\pm \to P_\pm^\perp$ such that $\abs{\varphi_\pm} \leq C \tau r$. 
        This piece of graph is bounded by $G$ and lies above $G$, i.e.,
        \begin{equation}
            F_{\pm} = \set{z + \varphi_{\pm}(z) | z \in P_{\pm} \, \text{ is such that} 
            \ z \cdot \xi_{\pm} \geq \psi(z_1) \cdot \xi_{\pm}} \cap A_r \cap H(100 c,\xi_0)
        \end{equation}
        and when $z \cdot \xi_{\pm} = \psi(z_1) \cdot \xi_{\pm}$, 
        we have $z + \varphi_{\pm}(z) = (z_1,\psi(z_1))$.
 
    Altogether $E_\infty$ looks like a $Y$-shaped 
    airplane with two wings that detaches itself a little from $\Gamma$.
    When the wings meet the vertical wall on $G \setminus \Gamma$, they make exact angles  
    of $2 \pi /3$.  And  when the two wings meet on $\Gamma \cap G$,
    there is no vertical wall and the angle may vary a little (not too much because of the small $C \tau$ 
    Lipschitz constant). At such a point, the situation is not so different from Case 2, even if the angle is not 
    exactly $2 \pi /3$. 
    This is our most complicated case, and its main feature is that $E_\infty$ does not have exactly 
    the same topology as $V$ locally, so that we cannot parameterize it nicely by $V$. Yet this will not really disturb the construction of retractions below.

    \msi{\bf Case 2 Bis.} 
    Our next case is when $Z$ is a truncated $\bY$ cone $Y$ containing $\xi_0$ and the spine of the triple junction passes through a point $\zeta_0 \in B(\xi_0,\nu c/10)$ with $\zeta_0 \ne \xi_0$; we have
    \begin{equation} \label{5c26}
        \begin{aligned}
            \dist(x,Y) \leq \varepsilon_0 r \ \text{ for } x \in E_\infty \cap A_r(200) \cap H(200 c,\xi_0),\\
            \dist(x,E_\infty) \leq \varepsilon_0 r \ \text{ for } x \in Y \cap A_r(200) \cap H(200 c,\xi_0).
        \end{aligned}
    \end{equation}
    We shall take $\nu$ small, and then his case is almost the same as the previous
    one. However, we prefer to distinguish the two, because this will allow us to construct
    our retraction in a more smooth way. The same description as in Case $2$ applies in 
    $A_r \cap H(100 c,\xi_0)$ but now we prefer to say that 
    $G$ is the graph over the line passing through $\zeta_0$ 
    (instead of $\xi_0$) of a $\tau$-Lipschitz function 
    $\psi : \R \to \R^{n-1}$ such that $\abs{\psi} \leq C \tau r$. We have as usual
    \begin{equation}
            \text{$\mathrm{dist}(\abs{x}^{-1}x,\zeta_0) \leq C \tau$ for $x \in G$}
    \end{equation}
    and, since $\abs{\zeta_0 - \xi_0} \leq \nu c/10$, $G$ and the vertical wall $F_v$ are entirely contained in $H(\nu c,\xi_0)$.

    Recall that we can take $\varepsilon_0$ in \eqref{5c26} as small as we want,
    depending on $c$ and $\nu$ in particular, so we get our description easily from
    multiple applications, at various points of $\Gamma \cap A_r$, of Theorem 31.1 of \cite{Dvv}, 
    for a description near $\Gamma$, and of Theorem 1.15 of \cite{Dcc}, for the other points.

    This will be good enough for the rest of the paper, taking $\nu$ small enough when we need it.

    \msi{\bf Case 3.} Next, we assume that $Z$ is a generic cone $V$ of type $\bV$.
    This time we can use Theorem~32.1 of \cite{Dvv} 
    to obtain a description of $E_\infty \cap A_r \cap H(100 c,\xi_0)$ which is similar to the one above (in \eqref{5b11}), except that now we can take $G = \Gamma$, 
    $F_v$ is essentially empty, i.e., reduced to a subset of $\Gamma$, and 
    $E_\infty$  is reduced to the union of two wings that start directly from $\Gamma$. 
    However, the parameter $\varepsilon_0$ in \cite[Theorem 32.1]{Dvv} should be small enough depending on the angle $\alpha \in (2\pi/3,\pi)$ of $V$.
    The following convention will help us to clarify this.
   If $\varepsilon_0$ is small enough, then depending on the angle of $V$, we are going to see how $E_\infty$ is described either by Case 1, Case 2 or by Case 3.

    In order to get the description of Case 2, we need to have \eqref{5b6} with an 
    $\varepsilon_0$ which is small enough (call the value $\varepsilon_{1}$; 
    it depends on $n$, $\nu$, $\tau$ and $c$).
    Then consider we are in Case 3 only when $V$ makes an angle  $\alpha \geq 2\pi/3 + \varepsilon_1/400$.
    In this case, there is an $\varepsilon_2 >0$ (that depends on $\varepsilon_1$) such that when $\varepsilon_0 \leq \varepsilon_2$ and we assume \eqref{5b6} for the generic cone $V$, we get the desired description, with only two faces, 
    and that make angles in $(\alpha - \varepsilon_1/800, \alpha + \varepsilon_1/800)$. In particular the angles are bigger than $2\pi/3 + \varepsilon_1/800$; there is no sharp angle anywhere in $A_r \cap H(100 c,\xi_0)$.

    Now return to the remaining case when $\alpha < 2\pi/3 + \varepsilon_1/400$.
    In this case, there exists a sharp cone $V'$ which is close to $V$ in the sense that for $x \in V$, 
    there exists $y \in V'$ such that $\abs{y} = \abs{x}$ and 
    $\abs{x - y} \leq \varepsilon_1 \abs{x}/400$ (and reciprocally, any point $x \in V'$ is similarly close to $V$).
    When $\varepsilon_0 < \varepsilon_1/2$, the condition \eqref{5b6} with $V$ and $\varepsilon_0$ imply the same condition \eqref{5b6} with $V'$ and $\varepsilon_1$, so we can pretend we are in Case 2. 

 In fact, we also do the same thing when our generic set $V$ is almost flat. We assume $\varepsilon_1$ small enough so that the description of Case 3 holds whenever $\varepsilon_0 \leq \varepsilon_1$. If $V$ makes an angle $\alpha \geq \pi - \varepsilon_1/400$, we replace $V$  with a plane that contains $L_0$ and pretend that we were in Case 1.
    This way, when we say that we are in Case 3, this also forces the two wings that compose $E_\infty$
    to make angles smaller than $\pi-\varepsilon_1/800$ everywhere along $\Gamma \cap A_r$, for some extremely small $\varepsilon_1 > 0$.

    \msi{\bf Case 4.} This time we assume that $Z$ is a half plane $H$ bounded by $L_0$.
    We use Theorem~31.1 of \cite{Dvv} and get a description as above, but with only one face $F$, 
    which is a half $\tau$-Lipschitz graph (over some plane that contains $L_0$) bounded by $\Gamma$.
    This was our last case with $\xi_0 \in X \cap L_0 \cap \d B(0,1)$.

    \ms
  Now we consider the cases where $\xi_0 \notin L_0$ (allowing also the possibility that $x_0 \in E_\infty \setminus \Gamma$ with the convention $L_0 = \emptyset$ in this case).

    \msi{\bf Case 5.}
    We assume $Z$ to be a plane $P$ through $\xi_0$. 
    In this case, $E_\infty$ coincides with a small Lipschitz graph over $P$ in $A_r \cap H(100 c,\xi_0)$. 
    This is as in \eqref{5b5}, but in the present case we may also assume that $E_\infty \cap A_r \cap H(100 c,\xi_0)$ does not meet $L_0$, because if it did at some point $x_1 \in L_0$, we would have $d(x_1,P) \leq \varepsilon_0 r$ so we could again move the plane $P$ a tiny bit so that it contains $L_0$, and pretend that we are in fact in Case 1.

    \msi{\bf Case 6.} We assume $Z$ to be a cone of type $\bY$ whose spine passes through $\xi_0$.
    In this case, in the region $A_r \cap H(100c,\xi_0)$, $E_\infty$ is a $C^1$ version of $\bY$, and it is composed of three faces that are $\tau$-Lipschitz graphs over the planes that are parallel to the faces of $Y$,
    and meet along the graph $G$ of some $\tau$-Lipschitz function $\psi$ over the line that contains $\xi_0$, with the usual $2\pi/3$ angle along $G$.
    As usual, we have $\abs{\psi} \leq C \tau r$ and thus
    \begin{equation}
        \text{$\mathrm{dist}(\abs{x}^{-1} x,\xi_0) \leq C \tau$ for $x \in G$.}
    \end{equation}
    We assume $\tau$ small enough so that $C \tau < \nu c$, and thus $G$ is contained in $H(\nu c,\xi_0)$.
    We can also assume as before that $E_\infty \cap A_r \cap H(100c,\xi_0)$ does not meet $L_0$ so that we do not need to worry about the sliding condition.
    The proof then comes directly from \cite{Dcc} (or \cite{Ta} when $n=3$).

    \section{Intrinsic projections on small Lipschitz graphs} 
    \label{S6a}

    In the next sections, we will find it pleasant to project points on pieces of 
    $E_\infty$ (or on pieces of curves in $E_\infty$) that are $\tau$-Lipschitz graphs,
    and in a way that only involves $E_\infty$ itself, and not a direction that would be chosen
    arbitrarily. We single out the issue and introduce relevant notation in this section, and
    then we will return to our main subject. 

    Here, we let $\tau$ denotes a constant in $]0,1]$ and $C$ denotes a generic constant $\geq 1$ that may depend on $n$.
    The constant $\tau$ will be assumed small enough, depending only $n$.

    First consider a closed set $E \subset \R^n$, which we assume to be a $\tau$-Lipschitz graph
    over some vector plane $P_0$ of dimension $d$ (for our purposes, $d=1$ and $d=2$ will be enough).
    We let $\pi^0$ denote the orthogonal projection onto $P_0$.
    Notice that $\pi^0$ is a bijective $1$-Lipschitz function from $E$ to $P_0$ and its inverse is $(1 + \tau)$-Lipschitz.
    One can deduce that for all $x \in E$ and for all $r > 0$,
    \begin{equation}
        C^{-1} r^d \leq \H^d(E \cap B(x,r)) \leq C r^d.
    \end{equation}
    The upper estimate can be extended to all $x \in \R^n$, $r > 0$ and the lower estimates can be extended to all $x \in \R^n$, $r > 0$ such that $\mathrm{dist}(x,E) \leq r/2$.

    We are going to explain how to make a Lipschitz projection on $E$ with a varying direction. This will be useful because we  prefer to use projections defined rather intrinsically, and that do not depend on the knowledge of an approximating plane $P_0$, for instance.

    In the Grassmannan space $\mathcal{G}_d = G(d,n)$ of linear planes $P$ of dimension $d$, we can measure the distance between two planes $P_1, P_2$ using the distance
    \begin{equation}
        \mathrm{dist}(P_1,P_2) := \norm{\pi_1 - \pi_2},
    \end{equation}
    where $\pi_i$ is the orthogonal projection onto $P_i$ and $\norm{\cdot}$ is any matrix norm
    (they are all equivalent, but we may need to choose a precise one later).
    As an example, if $L_1$ and $L_2$ are two linear lines generated by respective unit vectors $\nu_1$, 
    $\nu_2$ such that $\nu_1 \cdot \nu_2 \geq 0$, then
    \begin{equation}\label{eq_distL12}
        \mathrm{dist}(L_1,L_2) \sim |\nu_1-\nu_2|.
    \end{equation}
Remember that $E$ is a $\tau$-Lipschitz graph over $P_0$ so if $P \in \mathcal{G}_d$ is such that $\mathrm{dist}(P,P_0) \leq \tau$ with $\tau$ small enough, then $E$ is a $C \tau$-Lipschitz graph over $P$.
    Then, there exists a natural projection $p : \R^n \to E$ onto $E$ defined by
    \begin{equation}\label{eq_definitionP}
        p(z) \in E \ \text{  and  } \  z - p(z) \in P^\perp,
    \end{equation}
    satisfying in particular $p = \mathrm{id}$ on $E$. 
    The map $p$ is $(1 + C \tau)$-Lipschitz and we can assume $\tau$ small enough 
    so that $p$ is $2$-Lipschitz.
    As $(p - \mathrm{id})$ is $10$-Lipschitz and equals $0$ on $E$, we deduce that
    \begin{equation}\label{eq_p}
        |p(z) - z| \leq 10 \mathrm{dist}(z,E) \quad \text{for all $z \in \R^N$.}
    \end{equation}
    For $z \in \R^N$, notice that the projection $p(z) \in E$ is well-defined as soon as there exists $r > 0$ such that $\mathrm{dist}(z,E) \leq r$ and $E$ coincides with a $\tau$-Lipschitz graph above $P_0$ in $B(z,10r)$; we don't need to know $E$ outside $B(z,10r)$.

    Furthermore, one can control how the projections vary when one changes the plane of projections or the graph.
    Let us say that $E_1$ and $E_2$ are two $\tau$-Lipschitz graph over $P_0$ of two functions $\varphi_1$, $\varphi_2$ respectively, and let $P_1$ and $P_2$ are be two linear planes of dimension $d$ such that $\mathrm{dist}(P_i,P_0) \leq \tau$ for $i=1,2$.
    Let $p_1$ be the projection onto $E_1$ in the direction orthogonal to $P_1$, and respectively $p_2$ with respect to $E_2$ and $p_2$. Then, we are are going to show that for all $z,w \in \R^n$,
    \begin{equation}\label{eq_pdist}
        \abs{p_1(z) - p_2(w)} \leq C \mathrm{dist}(P_1,P_2) \, \mathrm{dist}(z, E_2) + C \norm{\varphi_1 - \varphi_2}_{\infty} + (1 + C \tau) \abs{z - w}.
    \end{equation}
    Since the projection $p_2$ onto $E_2$ is $(1 + C \tau)$-Lipschitz, we can estimate 
    \begin{align}
        \abs{p_1(z) - p_2(w)} &\leq \abs{p_1(z) - p_2(z)} + \abs{p_2(z) - p_2(w)}\\
                              &\leq \abs{p_1(z) - p_2(z)} + (1 + C \tau) \abs{z - w}
    \end{align}
    and we are left to control $\abs{p_1(z) - p_2(z)}$.
    Observe that there exists a point $q \in E_1$ such that $\abs{q - p_2(z)} \leq \norm{\varphi_1 - \varphi_2}_{\infty}$.
    Then, we show that $\abs{q - p_1(z)} \leq C \mathrm{dist}(P_1,P_2) \mathrm{dist}(z, E_2) +C \norm{\varphi_1 - \varphi_2}_{\infty}$ and (\ref{eq_pdist}) will follow.
    Let $\pi_1$ and $\pi_2$ be the orthogonal projections onto $P_1$ and $P_2$ respectively.
    As $E_1$ is a $C \tau$-Lipschitz graph above $P_1$ which contains both $q$ and $p_1(z)$, we have $\abs{q - p_1(z)} \leq (1 + C \tau) \abs{\pi_1(q - p_1(z))}$ whence
    \begin{align}
        \abs{q - p_1(z)} &\leq 2 \abs{\pi_1(q - p_2(z))} + 2 \abs{\pi_1(p_1(z) - p_2(z))}\\
                         &\leq 2 \norm{\varphi_1 - \varphi_2}_{\infty} + 2 \abs{\pi_1(p_1(z) - p_2(z))}
    \end{align}
    By definition, $z - p_1(z) \in P_1^\perp$ so $\pi_1(z - p_1(z)) = 0$ and a similar argument shows that $\pi_2(z - p_2(z)) = 0$. We deduce
    \begin{align} 
        \abs{\pi_1(p_1(z) - p_2(z))} &= \abs{\pi_1(z - p_2(z)) - \pi_2(z - p_2(z))}\\
        &\leq C \mathrm{dist}(P_1,P_2) \abs{z - p_2(z)} \leq C \mathrm{dist}(P_1,P_2) \mathrm{dist}(z,E),
    \end{align}
    where for the last line we used the fact that the matrix norm $\norm{\cdot}$ is equivalent to the operator norm induced by the euclidean distance.
    This ends the proof of (\ref{eq_pdist}).

    Now, the goal of this section is to build for each $x \in \R^n$ and $r > 0$ such that $\mathrm{dist}(x,E) \leq r/4$, an ``averaged'' linear plane $P_{x,r}$ of dimension $d$ such that $\mathrm{dist}(P_{x,r},P_0) \leq C \tau$.
    In this way, $E$ is also a $C \tau$-Lipschitz graph above $P_{x,r}$ and there is an associated natural projection $p_{x,r}$ onto $E$.
    The plane $P_{x,r}$ will depend in an intrinsinc %
    way on $E \cap B(x, 10 r)$ and it will depend in a Lipschitz way on $x$ and $r$.
    This will allow us to consider an intrinsic Lipschitz projection $x \mapsto p_{x,r}(x)$ onto $E$, with a varying direction.
    Unfortunately we could not use the closest point projection (which would be intrinsic too), because even for small Lipschitz graphs it is not well-defined (the closest point can jump).

    Set $\mu = \H^d_{\vert E}$, the restriction of Hausdorff measure, and also pick a smooth radial
    non-negative bump function $\varphi$ supported in $B(0,1)$ such that $\abs{\varphi} + \abs{\nabla \varphi} \leq 10$ and $\varphi(x) = 1$ on $B(0,1/2)$. 
    Also set $\varphi_r(y) = r^{-d} \varphi(y/r)$ as usual.
    For almost every $x \in E$, let $P_x$ denote the tangent $d$-plane to $E$ at $x$ (a vector plane), and $\pi_x$ the orthogonal projection on $P_x$ (a linear map). 
    Since $E$ is a $\tau$-Lipschitz graph above $P_0$, we know that for almost-every $x \in E$,
    \begin{equation}
        \norm*{\pi_x - \pi^0} \leq C \tau,
    \end{equation}
    where we recall that $\pi^0$ is the orthogonal projection onto $P_0$.
    For $x \in \R^n$ and $r > 0$ such that $\mathrm{dist}(x,E) \leq r/10$, we define an average projection by
    \begin{equation}\label{6c1}
        \pi^0_{x,r}(w) = \left(\int_E \varphi_r(z-x) \dd{\mu(z)}\right)^{-1} \int_E \varphi_r(z-x) \pi_z(w) \dd{\mu(z)} \quad \text{for all $w \in \R^n$}.
    \end{equation}
    The denominator is always $\leq C$ and we also required that $\dist(x,E) \leq r/10$ to make sure that it is also $\geq C^{-1}$. 
    A priori, $\pi^0_{x,r}$ is a linear map and not necessarily a projection onto a plane.
    However, it stays close to $\pi^0$, that is, $\norm*{\pi^0_{x,r} - \pi^0} \leq C \tau$ because $\norm*{\pi_z - \pi^0} \leq C \tau$ for $z \in E$.

    Since we intend to use some differential geometry in the
    space $\cL$ of linear mappings on $\R^n$, let us represent them by matrices, and assume that $\norm{\cdot}$ is a smooth 
    norm, such as the square root of the sum of the squares of the coefficients.
    Next call $\cP_d \subset \cL$ the set of linear
    mappings that are the orthogonal projection on some vector $d$-plane $P$; this is a smooth
    (and biLipschitz) realization of the Grassman manifold $\cG_d$.
    This is a smooth submanifold of 
    $\cL$, with bounded curvature (by compactness). 
    Thus, there is a constant $\delta_0$ and a smooth $C$-Lipschitz projection $\Pi$ from a $\delta_0$-neighborhood of $\cP_d$ onto $\cP_d$.
    If $\tau$ is small enough compared with $\delta_0$, we can pick
    \begin{equation} \label{6c2}
        \pi_{x,r} = \Pi(\pi^0_{x,r}) \in \cP_d,
    \end{equation}
    and call $P_{x,r}$ the $d$-plane such that $\pi_{x,r}$ is the orthogonal projection onto $P_{x,r}$. 
    The plane $P_{x,r}$ is still close to $P_0$, i.e.,
    \begin{equation}\label{eq_distpi0}
        \norm*{\pi_{x,r} - \pi^0} \leq C \tau
    \end{equation}
    as we can estimate
    \begin{align}
        ||\pi_{x,r} - \pi^0||   &\leq ||\pi_{x,r} - \pi_{x,r}^0|| + ||\pi_{x,r}^0 - \pi^0||\\
        &\leq ||\left(\Pi - \mathrm{id}\right) \left(\pi^0_{x,r} - \pi^0\right)|| + ||\pi_{x,r}^0 - \pi^0||   \nonumber
        \leq C \tau.
    \end{align}
    Provided that $\tau$ is small enough, $E$ is then also a $C\tau$-Lipschitz graph over $P_{x,r}$.
    The main feature of $\pi_{x,r}$ is that it is defined in an intrinsic way, but notice that
    by construction, it depends nicely on $x$ and $r$, i.e., for $x,y \in \R^n$ and $r, s > 0$ such that $\mathrm{dist}(x,E) \leq r/10$ and $\mathrm{dist}(y,E) \leq s/10$,
    \begin{equation}\label{6c3}
        ||\pi_{x,r} - \pi_{y,s}|| \leq C \tau r^{-1}(|x-y| + |r-s|).
    \end{equation}
    For the rather boring proof, notice first that we can assume $\abs{x - y} \leq r$ without loss of generality because we always have $\norm{\pi_{x,r} - \pi_{y,s}} \leq \norm{\pi_{x,r} - \pi^0} + \norm{\pi^0 - \pi_{y,r}} \leq C \tau$.
    Similarly, (\ref{6c3}) is always true when $\abs{r - s} \geq r/2$ so we can assume $r/2 \leq s \leq 2r$ without loss of generality.
    As $\Pi$ is Lipschitz, it suffices to control $\norm{\pi^0_{x,r} - \pi^0_{y,s}}$.
    We write
    \begin{equation}
        \pi^0_{x,r} = \pi^0 + \left(\int_E \varphi_r(z-x) d\mu(z) \right)^{-1} \int_E \varphi_r(z-x) \left(\pi_z - \pi^0\right) d\mu(z),
    \end{equation}
    and we observe that for all $x, y \in \R^n$ and $r, s > 0$ such that 
    $r/2 \leq s \leq 2r$, we have
    \begin{equation}
        \abs{\varphi_r(x) - \varphi_s(y)} \leq C r^{-(d+1)} \left(\abs{x - y} + \abs{r - s}\right).
    \end{equation}
    We deduce that for all $x, y \in \R^n$ such that $\abs{x - y} \leq r$ and $r/2 \leq s \leq 2r$, we have
    \begin{equation}
        \abs{\int_E \varphi_r(z - x)\dd{\mu(z)} - \int_E \varphi_s(z - y) \dd{\mu(z)}} \leq C r^{-1} \left(\abs{x - y} + \abs{r - s}\right).
    \end{equation}
    We can control similarly
    \begin{multline}
        \norm{\int_E \varphi_r(z - x) (\pi_z - \pi^0) \dd{\mu(z)} - \int_E \varphi_s(z - y) (\pi_z - \pi^0) \dd{\mu(z)}} \\\leq C \tau r^{-1} \left(\abs{x - y} + \abs{r - s}\right),
    \end{multline}
    using in addition the fact that $\norm{\pi_z - \pi^0} \leq C \tau$.
    We skip the rest of the proof for the reader's relief.
  Notice as usual that in order for $P_{x,r}$ to be well-defined, we just need that $\mathrm{dist}(x,E) \leq r/10$ and that $E$ coincides with a $\tau$-Lispchitz graph in $B(x,r)$. Similarly for (\ref{6c3}), it suffices that $E$ coincides with a $\tau$-Lipschitz graph in $B(x,r)$ and $B(y,s)$.

\ms
        In the rest of this section, we show how to define a projection onto $E$ which preserves spheres. We will work under the assumption that $E$ coincides with a Lipschitz graph in $A_r \cap H(100c,\xi_0)$, precisely,
        \begin{equation}
            E \cap A_r \cap H(100c,\xi_0) = \set{z + \varphi(z) | z \in P_0} \cap A_r \cap H(100 c,\xi_0)
        \end{equation}
        where $r > 0$ is any radius (the scale at which the description of $E$ holds), $c \in (0,1)$ is a small constant, $\xi_0$ is a unit vector, $P_0$ is a vector plane containing $\xi_0$ and $\varphi : P_0 \to P_0^\perp$ is $\tau$-Lipschitz function satisfying $\abs{\varphi} \leq C \tau r$.

        For $x \in \R^n$, we let $z$ and $z'$ denote the components of $x$ in the sum $\R^n = P_0 + P_0^\perp$, thus $x = z + z'$ or with a slight abuse of notation, $x = (z,z')$. 
        We complete $\xi_0$ with a unit vector $\xi_1 \in P_0$ to form an orthogonal basis $(\xi_0,\xi_1)$ of $P_0$.
        Then for $z \in P_0$, we let $z_0$ and $z_1$ be the coordinates of $z$ in this basis, i.e., $z = z_0 \xi_0 + z_1 \xi_1$ or $z = (z_0,z_1)$.

        We fix a constant $\tau_0 \in (0,1)$ which is small enough (depending on $n$) so that all the properties of this section so far hold true for $\tau_0$-Lipschitz graphs.
        We will assume $c$ small enough depending on $n$ and $\tau_0$. And we will assume $\tau$ small enough compared to $\tau_0$ and $c$.
        The first step in our construction is to justify that for $t \in (r/2,2r)$, the set
        \begin{equation}
            E \cap S_{t} \cap B(t\xi_0,100ct)
        \end{equation}
        is a $\tau_0$-Lipschitz graph above $\R \xi_1$. This will allow us to make a projection onto $E \cap S_t \cap B(t\xi_0,100ct)$ in a well-chosen direction. 
        For this purpose, we assume $c$ small enough depending on $n$ and $\tau_0$ so that for all $t > 0$, the spherical cap $S_t \cap B(t\xi_0,100ct)$ is a graph of some $\tau_0/10$-Lipschitz function $\psi_{t}$ above $\xi_0^\perp$. 
        For $x \in E \cap S_{t} \cap B(t\xi_0,100ct)$, one has both $z' = \varphi(z_0,z_1)$ and $z_0 = \psi_{t}(z_1,z')$, whence 
        \begin{equation}
            z_0 = \psi_{t}(z_1,\varphi(z_0,z_1)).
        \end{equation}
        It easy to deduce, if $\tau$ is small enough depending on $n$ and $\tau_0$, that $z_0$ is uniquely determined by $z_1$ (because the right-hand side is contracting in the variable $z_0$) and even that $z_0$ a $\tau_0/2$-Lipschitz function of $z_1$ (for the same reason). Similarly, $z'$ is a $\tau_0/2$-Lipschitz function of $z_1$. This justifies our claim.

        For each $x \in S_r \cap B(r\xi_0,20cr)$ such that $\mathrm{dist}(x,E) \leq 10^{-3}cr$, we aim to associate a direction $L_{x,r} \in G(1,n)$ which only depends on $E \cap B(x,10^{-1}cr)$ and such that $\mathrm{dist}(L_{x,r},\R\xi_1) \leq C \tau_0$. In this way, $L_{x,r}$ can be used to project onto $E \cap S_r$ with the usual formula: 
        \begin{equation}\label{eq_defpxr}
            p_{x,r}(z) \in E \cap S_r \quad \text{and} \quad z - p_{x,r}(z) \in L_{x,r}^\perp.
        \end{equation}
        We shall say that $p_{x,r}(z)$ is \ul{well-defined} when $x \in S_r \cap B(r\xi_0,20cr)$ is such that $\mathrm{dist}(x,E) \leq 10^{-3}cr$ and $z \in B(r\xi_0,20cr)$ is such that $\mathrm{dist}(z, E \cap S_r) \leq 10^{-3}cr$.

        We could allow $E$ to be bounded by a curve (similarly as the faces $F_\pm$ defined in Section~\ref{S5}); everything works the same as long as we assume additionally that $E$ coincides with a $\tau$-Lipschitz graph above $P_0$ in $B(x,10^{-1}cr)$ and $E \cap S_r$ coincides with a $\tau_0$-Lipschitz graph above $\R \xi_1$ in $B(z,10^{-1}cr)$. Under these assumptions, the bounding curve is sufficiently out of reach of the projection so that it can be ignored.

        We also want $p_{x,r}(z)$ to be well-defined for $z = x$ itself. For this purpose, we will prove at the end of this section that for all $x \in E \cap S_r \cap B(r\xi_0,20cr)$,
        \begin{equation}\label{eq_distES00}
            \mathrm{dist}(x, E \cap S_r) \leq C \mathrm{dist}(x,E),
        \end{equation}
        where $C \geq 1$ only depends on $n$. 
        Letting for $c_1 > 0$,
        \begin{equation}
            W(c_1) := \set{x \in \R^n \setminus \set{0} | \mathrm{dist}(x,E) \leq c_1 \abs{x}},
        \end{equation}
        inequality (\ref{eq_distES00}) shows indeed that if $c_1$ is small enough depending on $n$ and $c$, then $p_{x,r}(z)$ is well-defined for all $x, z \in W(c_1) \cap S_r \cap B(r\xi_0,20cr)$.

        We can define similarly a projection onto $E \cap S_{t}$ for any other radius $t \in (r/2,2r)$ and our last request is that $p_{x,r}$ is Lipschitz in $x$ and $r$.
        Unfortunately, replacing directly $E$ by $E \cap S_r$ in formula (\ref{6c1}) would not provide a direction of projection $L_{x,r}$ which is Lipschitz across spheres.
        We shall consider instead the direction
        \begin{equation}
            L_{x,r} := P_{x,10^{-1}cr} \cap (\R x)^{\perp},
        \end{equation}
  where $P_{x,10^{-1}cr}$ is the intrinsic average plane of $E$ in $B(x,10^{-1}cr)$ as defined earlier in this section. Here, we may assume $c$ small enough so that the condition $x \in S_r \cap B(r\xi_0,20cr)$ implies 
  $x \cdot \xi_0 \geq r/2$ and in particular that $P_{x,r}$ is not a subspace of $(\R x)^{\perp}$. 
  Therefore, $L_{x,r}$ is clearly a line.

        Next  we prove that for all $s, t \in (r/2,2r)$, for all $x, y$ and $z, w$ such that $p_{x,t}(z)$ and $p_{y,s}(w)$ are well-defined, we have
        \begin{equation}\label{eq_pLipschitz}
            \abs{p_{x,t}(z) - p_{y,s}(w)} \leq C \abs{x - y} + C \abs{z - w}.
        \end{equation}
        The reader may want to skip the proof and go directly to Section \ref{S7}.
        Our first step is to show that, under the same assumptions,
        \begin{equation}\label{eq_Lxy0}
            \mathrm{dist}(L_{x,t}, \R \xi_1) \leq \tau_0
        \end{equation}
        and
        \begin{equation}\label{eq_Lxy}
            \mathrm{dist}(L_{x,t},L_{y,s}) \leq C \abs{x - y}.
        \end{equation}
        In particular, property (\ref{eq_Lxy0}) makes sure that $E \cap S_r$ is a graph above $L_{x,r}$ so (\ref{eq_defpxr}) is well-defined.
        The main point is to show that
        \begin{align} 
            \mathrm{dist}(L_{x,t},\R \xi_1) 
            &\leq C \mathrm{dist}(P_{x,10^{-1}ct},P_0) + C \mathrm{dist}(\R x,\R \xi_0), \\ 
            \mathrm{dist}(L_{x,t},L_{y,s})   &\leq C \mathrm{dist}(P_{x,10^{-1}ct},P_{y,10^{-1} c s}) + C \mathrm{dist}(\R x,\R y).
        \end{align} 
        The first estimate implies (\ref{eq_Lxy0}) since $\mathrm{dist}(\R x, \R \xi_0) \leq C c \leq C \tau_0$ (recall that $x \in B(t\xi_0,20cr)$ and $c$ is allowed to be small depending on $\tau_0$).
        The second estimate also implies (\ref{eq_Lxy}) since in view of (\ref{6c3}), we have
        \begin{align} \label{6a35}
            \mathrm{dist}(L_{x,t},L_{y,s})   &\leq C \tau r^{-1} (\abs{x - y} + \abs{t - s}) + C \abs{\frac{x}{\abs{x}} - \frac{y}{\abs{y}}}\\
                                             &\leq C r^{-1} \abs{x - y},
        \end{align}
        and where for the last line, we used the fact that $\abs{t - s} \leq \abs{x - y}$ (recall that $\abs{x} = t$ and $\abs{y} = s$).
        The proof of each estimate is essentially the same, so we only detail the second one. 
        Observe that $P \mapsto P^\perp$ is an 
        isometry between Grassmannians. Thus letting $Q_1$ denote the orthogonal of $P_{x,10^{-1}ct}$ (respectively, $Q_2$ the orthogonal of $P_{y,10^{-1} cs}$), the second estimate amounts to
        \begin{equation}\label{eq_Qxy}
            \mathrm{dist}(Q_1 + \R x, Q_2 + \R y) \leq C \mathrm{dist}(Q_1,Q_2) + C \mathrm{dist}(\R x, \R y),
        \end{equation}

        At this point, we find convenient to set the matrix norm $\norm{\cdot}$ as the operator norm induced by the euclidean distance; this yields in particular that for all $P_1, P_2 \in G(d,n)$,  
        \begin{equation}\label{eq_distVW}
            \mathrm{dist}(P_1,P_2) \leq \sup\set{\mathrm{dist}(z,P_2) | z \in P_1, \abs{z} \leq 1} + \sup \set{\mathrm{dist}(z,P_1) | z \in P_2, \abs{z} \leq 1}.
        \end{equation}
        Let us justify this claim. 
        Let $\pi_1$ and $\pi_1$ be the orthogonal projections onto $P_1$ and $P_2$ respectively.
        For all $x \in \R^n$ with $\abs{x} \leq 1$, the vector $x$ can be decomposed as 
        $x = u + v$ in the sum $\R^n = P_1 + P_1^\perp$ and thus
        \begin{equation}
            \abs{\pi_1(x) - \pi_2(x)} \leq \abs{(\pi_1 - \pi_2)(u)} + \abs{(\pi_1 - \pi_2)(v)},
        \end{equation}
        In order to estimate the first term,  we observe that since $u \in P_1$, 
        we have $\abs{(\pi_1 - \pi_2)(u)} = \abs{u - \pi_2(u)} = \mathrm{dist}(u,P_2)$.
        Next, we estimate $\abs{(\pi_1 - \pi_2)(v)}$ by duality. 
        For all $z \in P_2$ such that $\abs{z} \leq 1$, we have
        \begin{equation}
            \abs{(\pi_1 - \pi_2)(v) \cdot z} = \abs{v \cdot (\pi_1 - \pi_2)(z)} \leq \abs{v} \abs{(\pi_1 - \pi_2)(z)} \leq \mathrm{dist}(z,P_1)
        \end{equation}
        and as $(\pi_1 - \pi_2)(v) \in P_2$ (remember $v \in P_1^\perp$), we deduce that $\abs{(\pi_1 - \pi_2)(v)} \leq \mathrm{dist}(z,P_1)$.
        This proves (\ref{eq_distVW}).

        In  order to show (\ref{eq_Qxy}), it suffices to prove that
        \begin{equation}
            \sup\set{\mathrm{dist}(z,Q_1 + \R x) | z \in Q_2 + \R y, \abs{z} \leq 1} \leq \mathrm{dist}(Q_1,Q_2) + \mathrm{dist}(\R x, \R y),
        \end{equation}
        the other case being symmetric. Let $z \in Q_2 + \R y$ be decomposed as $z = u + v$, where $u \in Q_2$ and $v \in \R y$. Letting $u'$ denote the orthogonal projection of $u$ onto $Q_1$ and $v'$ denote the orthogonal projection of $v$ onto $\R x$, we clearly have
        \begin{equation}
            \mathrm{dist}(z, Q_1 + \R x) \leq \abs{u - u'} + \abs{v - v'} \leq \mathrm{dist}(Q_1,Q_2) \abs{u} + \mathrm{dist}(\R x, \R y) \abs{v}.
        \end{equation}
        All is left to do is to check that $\abs{u} + \abs{v} \leq 2 \abs{z}$.
        Remember that $y \in S_s \cap B(s\xi_0,20cs)$ so $\abs{y/\abs{y} - \xi_0} \leq 100 c$ and this implies that $y$ is nearly orthogonal to $Q_2$ in the sense that for all $u \in Q_2$, since $u \cdot \xi_0 = 0$, we have $\abs{u \cdot y} \leq 100 c \abs{u} \abs{y}$.
        As $v \in \R y$, we deduce that
        \begin{align}
            \abs{u + v}^2 &\geq \abs{u}^2 + \abs{v}^2 - 2 \abs{u \cdot v}\\
            &\geq \abs{u}^2 + \abs{v}^2 - 200 c \abs{u} \abs{v} \geq (1 - 200c) (\abs{u}^2 + \abs{v}^2)
        \end{align}
        and we can assume $c$ sufficiently small so that $\abs{u + v} \geq (\abs{u} + \abs{v})/2$. 
        This concludes the proof of (\ref{eq_Qxy}).

        Our next step in the proof of (\ref{eq_pLipschitz}) is to show that
        \begin{equation}\label{eq_pLipschitz2}
            \abs{\frac{p_{x,t}(z)}{r} - \frac{p_{y,s}(w)}{s}} \leq (1 + C \tau_0) \abs{\frac{z}{r} - \frac{w}{s}} + C \tau_0 r^{-1} \abs{x - y}.
        \end{equation}
        One can easily see that this implies (\ref{eq_pLipschitz}) and we leave the details to the reader.
        The idea behind (\ref{eq_pLipschitz2}) is that $z \mapsto t^{-1} p_{x,t}(z)$ is the projection of $t^{-1}z$ onto the graph $(t^{-1} E) \cap S_1 \cap B(\xi_0,100c)$ in the direction orthogonal to $L_{x,t}$ (the definition (\ref{eq_defpxr}) is well-preserved under rescaling). Thus, we want to deduce (\ref{eq_pLipschitz2}) from (\ref{eq_pdist}). We already know that $\mathrm{dist}(t^{-1} z,(t^{-1} E) \cap S_1) \leq 10^{-3} c \leq \tau_0$ (this is one of the assumptions in order for $p_{x,t}(z)$ to be considered well-defined) and that
        \begin{equation}\label{eq_LLipschitz2}
            \mathrm{dist}(L_{x,t},L_{y,s}) \leq C r^{-1} \abs{x - y}.
        \end{equation}
        Then, considering $(r^{-1} E) \cap S_1 \cap B(\xi_0,100c)$ and $(s^{-1} E) \cap S_1 \cap B(\xi_0,100c)$ as graphs of some function $\varphi_r$ and $\varphi_s$ respectively above $\R \xi_1$, we show that
        \begin{equation}\label{eq_varphits}
            \norm{\varphi_r - \varphi_s}_{\infty} \leq C \tau r^{-1} \abs{t - s}.
        \end{equation}
        We recall that for all $t \in (r/2,2r)$, the spherical cap $S_{t} \cap B(t\xi_0,100ct)$ is a graph of some $\tau_0/10$-Lipschitz function $\psi_r$ above $\xi_0^\perp$. 
        In fact, all these graphs are equivalent to each other under rescaling, i.e., $t^{-1} \psi_{t}(z_1,z') = \psi_1({t}^{-1} z_1, {t}^{-1} z')$. 
        We deduce that any point $(z_0,z_1,z')$ in the graph $(t^{-1} E) \cap B(\xi_0,100c)$ satisfy
        \begin{equation}
            z_0 = \psi_1(z_1, t^{-1}\varphi(t z_0,t z_1))
        \end{equation}
        (and similarly for $s$). 
        Then for all $(z_0,z_1,z') \in (t^{-1} E) \cap B(\xi_0,100c)$ and for all $(w_0,w_1,w') \in (s^{-1} E) \cap B(\xi_0,100c)$ with $z_1 = w_1$, we have
        \begin{equation}
            \abs{z_0 - w_0} \leq C \abs{\frac{\varphi(tz_0,tz_1)}{t} -  \frac{\varphi(sw_0,sw_1)}{s}} \leq C \tau r^{-1} \abs{t - s} + C \tau \abs{z_0 - w_0},
        \end{equation}
        where we used in particular the fact that $\varphi$ is $\tau$-Lipschitz, that $\abs{\varphi} \leq C \tau r$, and that $\abs{z}, \abs{w} \leq 1$. This implies in turn that $\abs{z_0 - w_0} \leq C \tau r^{-1} \abs{t - s}$. One can show similarly that $\abs{z' - w'} \leq C \tau r^{-1} \abs{t - s}$. This ends the proof of (\ref{eq_varphits}). As $\abs{x} = t$, $\abs{y} = s$ and $\tau$ is small depending on $\tau_0$, this implies $\norm{\varphi_t - \varphi_s} \leq C \tau_0 r^{-1} \abs{x - y}$, and in turn (\ref{eq_pLipschitz}) by (\ref{eq_pdist}).

        \begin{rem}\label{rem_pLipschitz3}
            If $x$ and $y$ are collinear with the same orientation, then   \eqref{6a35} yields  
        \begin{equation}
            \mathrm{dist}(L_{x,t},L_{y,s}) \leq C \tau r^{-1} \abs{x - y}
        \end{equation}
        and we can deduce a more precise variant of (\ref{eq_pLipschitz2}), where the error term in 
        $\abs{x - y}$ is controlled by $C \tau r^{-1}$ instead of $C r^{-1}$:  
        \begin{equation}\label{eq_pLipschitz3}
            \abs{\frac{p_{x,t}(z)}{r} - \frac{p_{y,s}(w)}{s}} 
            \leq (1 + C \tau_0) \abs{\frac{z}{t} - \frac{w}{s}} + C \tau r^{-1} \abs{x - y}. 
        \end{equation}
        This will be useful in Section \ref{S8a}.
        \end{rem}

        As promised, we conclude this section by showing that for $x \in S_r \cap B(r\xi_0,20cr)$,
        \begin{equation}\label{eq_distES0}
            \mathrm{dist}(x, E \cap S_r) \leq C \mathrm{dist}(x,E).
        \end{equation}
        The proof also works for any other radius in $(r/2,2r)$.
      
        In the next section, we will need to know that (\ref{eq_distES0}) also holds true when $E$ is bounded by a Lipschitz curve above $\R \xi_0$, that is, when $E \cap A_r \cap H(100c,\xi_0)$ is equal to 
        \begin{equation}\label{eq_Ebounded}
            \set{z + \varphi(z) | z \in P_0 \ \text{such that} \ z \cdot \xi_1 \geq \psi(z_0) \cdot \xi_1} \cap A_r \cap H(100 c,\xi_0),
        \end{equation}
        where $z = z_0 \xi_0 + z_1 \xi_1$, 
        $\varphi : P_0 \to P_0^\perp$ is a $\tau$-Lipschitz function which satisfies $\abs{\varphi} \leq C \tau r$, and similarly, $\psi : \R \xi_0 \to \xi_0^\perp$ is a $\tau$-Lipschitz function which satisfies $\abs{\psi} \leq C \tau r$ (we identify $\R \xi_0$ to $\R$ for notational convenience). 
        Therefore, we detail directly here this more complicated variant. 
        The proof still works if $E$ is bounded by two curves, the strategy is the same. As a consequence, (\ref{eq_distES0}) extends to the case $E$ is given as an union of Lipschitz graphs bounded by curves; the estimate holds independently for each face and therefore for their union.
        The reader is again allowed to skip the following proof.

        Let justify rapidly that (for $x \in S_r \cap B(r\xi_0,20cr)$ as above) we don't need to distinguish between 
        the distance from $x$ to $E$ and the distance from $x$ to $E \cap A_r \cap H(100c,\xi_0)$.
        One can see that $\mathrm{dist}(r\xi_0,E) \leq C \tau r$, since for instance (\ref{eq_Ebounded}) shows that $r\xi_0 + \psi(r\xi_0)$ belong to $E$ (assuming $\tau$ small enough compared to $c$ if necessary). It follows that $\mathrm{dist}(x,E) \leq 20cr + C \tau r \leq 30 cr$ but since $\overline{B}(x,30cr) \subset E \cap A_r \cap H(100c,\xi_0)$, the distance from $x$ to $E$ is attained in $E \cap A_r \cap H(100c,\xi_0)$.

        We decompose $x$ as $x = z + z'$ in the sum $P_0 + P_0^\perp$, with $z = z_0 \xi_0 + z_1 \xi_1$ as usual. 
        Our first idea for (\ref{eq_distES0}) is to look for some $s > 0$ such that $\abs{sz + \varphi(sz)} = r$ but since $sz + \varphi(sz)$ may be on the bad side of   the bounding curve, we need to reproject $sz$ first.
        For $s \in [1/20,20]$, we set 
        \begin{equation}
            z(s) = s z_0 \xi_0 + \max(s z \cdot \xi_1, \psi(sz_0) \cdot \xi_1) \xi_1
        \end{equation}
        so that $z(s)$ is on the good side of the bounding curve, i.e., $z \cdot \xi_1 \geq \psi(sz_0) \cdot \xi_1$.
        Next, let us take note of a few properties.
        Notice that
        \begin{equation}\label{eq_zsz0}
            \abs{\varphi(z) - z'} \leq C \mathrm{dist}(x,E).
        \end{equation}
        since the map $x \mapsto \varphi(z) - z'$ is $C$-Lipschitz on $\R^n$ and is zero on $E \cap A_r \cap H(100c,\xi_0)$.
        Our second property is that for all $x \in S_r \cap B(r\xi_0,20cr)$,
        \begin{equation}\label{eq_zsz}
            \abs{z(s) - sz} \leq C \mathrm{dist}(x,E) + C \tau \abs{1 - s} r.
        \end{equation}
        In order to prove (\ref{eq_zsz}), observe first that for $s$ fixed, $z \mapsto z(s)$ is $C$-Lipschitz on $P_0$ (a maximum of two Lipschitz functions is Lipschitz).
        Moreover, when $z \cdot \xi_1 \geq \psi(z_0) \cdot \xi_1$ (as for instance when $z$ is the projection of a point $x \in E \cap A_1 \cap H(100c,\xi_0)$), we can estimate
        \begin{equation}
            \abs{z(s) - sz} \leq C \tau \abs{1 - s} r.
        \end{equation}
        Indeed, if $s z \cdot \xi_1 \geq \psi(sz_0) \cdot \xi_1$, then $z(s) = sz$ and this estimate is trivial. Otherwise, $z(s) - sz$ only has a $\R \xi_1$ coordinate, which we can control as follows:
        \begin{align}
            0 \leq \psi(sz_0) \cdot \xi_1 - s z \cdot \xi_1 &= \left(\psi(sz_0) - s\psi(z_0)\right) \cdot \xi_1 + s \left(\psi(z_0) \cdot \xi_1 - z \cdot \xi_1\right)\\
            &\leq \left(\psi(sz_0) - s\psi(z_0)\right) \cdot \xi_1\\
            &\leq \abs{\psi(sz_0) - \psi(z_0)} + \abs{\psi(z_0) - s\psi(z_0)} \leq  C \tau \abs{1 - s} r,
        \end{align}
        where for the last line, we used the fact that $\psi$ is $\tau$-Lipschitz, $\abs{\psi} \leq C \tau r$ and that $\abs{z} \leq \abs{x} \leq 10$.
        Then , we deduce (\ref{eq_zsz}) from the fact that $x \mapsto z(s) - sz$ is $C$ Lipschitz on $\R^n$ and is bounded by $C \tau \abs{1 - s} r$ when $x \in E$

        For $s \in [1/20,20]$, we set $x(s) = z(s) + \varphi(z(s))$. 
        We recall that since $\abs{x} = r$, we have $\abs{z} \leq r$. 
        Moreover, $c$ is small enough so that the condition $x \in S_r \cap B(r\xi_0,20cr)$ implies $x \cdot \xi_0 \geq r/2$ and thus $\abs{z} \geq r/2$. 
        In view of (\ref{eq_zsz}) and the fact $\mathrm{dist}(x,E) \leq 20 cr$ and that $r/2 \leq \abs{z} \leq r$, it is clear that $\abs{x(s)} \in (r/100,100r)$ and that
        \begin{equation}
            \mathrm{dist}(\abs{x}^{-1}x, \abs{x(s)}^{-1} x(s)) \leq 2 \abs{sx}^{-1} \abs{sx - x(s)} \leq C \tau
        \end{equation}
        so $x(s) \in A_r \cap H(100c,\xi_0)$ (assuming possibly $\tau$ small enough depending on $c$) and in turn that $x(s) \in E$ by (\ref{eq_Ebounded}).
        Now, in order to show (\ref{eq_distES0}), we want to prove that there exists some $s \in [1/20,20]$ such that $\abs{x(s)} = r$ and
        \begin{equation}\label{eq_xse}
            \mathrm{dist}(x, x(s)) \leq C \mathrm{dist}(x,E).
        \end{equation}
        The map $s \mapsto \abs{x(s)}$ goes continuously from a value $\leq r/2$ at $s = 1/20$ to a value $\geq 2r$ at $s=20$, so there exists $s \in [1/20,20]$ such that $\abs{x(s)} = r$. 
 As $x = z + z'$ and $z(s) = z(s) + \varphi_i(z(s(s))$, (\ref{eq_xse}) amounts to showing that 
        $\abs{z - z(s)} + \abs{z' - \varphi(z(s))} \leq C \mathrm{dist}(x,E)$.
        In fact, using (\ref{eq_zsz0}), we have 
        \begin{equation}
            \abs{z' - \varphi(z(s))} \leq \abs{z' - \varphi(z)} + \abs{\varphi(z) - \varphi(z(s))} \leq C \mathrm{dist}(x,E) + C \abs{z - z(s)},
        \end{equation}
        so all is left to do is to check that $\abs{z - z(s)} \leq C \mathrm{dist}(x,E)$.
        Observe that if two unit vectors $x = z + z'$ and $y = w + w'$ in $P_0 + P_0^\perp$ are such that $\min(\abs{z'},\abs{w'}) \leq 1/2$, then $\abs{\abs{z} - \abs{w}} \leq C \abs{z' - w'}$.
        Applying this observation to $x = z + z'$ and $x(s) := z(s) + \varphi(z(s))$ which have the same norm $r$, we obtain
        \begin{align}
            \big|\abs{z(s)} - \abs{z}\big| &\leq C \abs{\varphi(z(s)) - z'}\\
            &\leq C \abs{\varphi(z(s)) - \varphi(z)} + C \abs{\varphi(z) - z'}\\
            &\leq C \tau \abs{z(s) - z} + C \mathrm{dist}(x,E).
        \end{align}
        Using $\abs{z(s) - z} \leq \abs{z - sz} + \abs{z(s) - sz}$ and (\ref{eq_zsz}), we have
        \begin{equation}\label{eq_zzz}
            \abs{z(s) - z} \leq C \abs{1 - s} r + C \mathrm{dist}(x,E) 
        \end{equation}
        so in particular
        \begin{equation}\label{eq_bigzsz}
            \big|\abs{z(s)} - \abs{z}\big| \leq C \tau \abs{1 - s} r + C \mathrm{dist}(x,E).
        \end{equation}
        By (\ref{eq_zsz}) and the fact that the map $u \mapsto \big|\abs{z} - \abs{u}\big|$ is $1$-Lipschitz, we can also bound from below
        \begin{align}
            \big|\abs{z} - \abs{z(s)}\big|  &\geq \big|\abs{z} - \abs{sz}\big| - \abs{z(s) - sz}\\
                                            &\geq C^{-1} \abs{1 - s}r  - C \mathrm{dist}(x,E)\label{eq_bigzsz2}
        \end{align}
        and it follows that $\abs{1 - s}r \leq C \mathrm{dist}(x,E)$. Finally, (\ref{eq_zzz}) shows that $\abs{z - z(s)} \leq C \mathrm{dist}(x,E)$ as wanted.
  
    \qed

\section{Projections along the spheres \texorpdfstring{$S_r$}{Sr}} 
\label{S7}

We keep the notation and assumptions 
of Section \ref{S5},
and we shall use the description of $E_\infty$ in the $A_r$ that was obtained there
to construct a first projection on $E_\infty$ that preserves the spheres. 

Denote by $S_r$ the sphere centered at the origin and with radius $r$
(so that $S_1 = \d B(0,1)$, for instance), and set $r(x) = |x|$ for $x\in \R^n$.

\begin{pro}\label{l6a1}
    Suppose as above that $0 \in E_\infty$, and that $r_0$ and $\varepsilon_0$ were chosen as in Section \ref{S5}.
    Then  there is a projection $\pi$, defined on the region 
    \begin{equation} \label{7a1}
        W(c_1) = \set{x \in B_0 = B(0,r_0) | \dist(x, E_\infty) \leq c_1 r(x)},
    \end{equation}
    where we can choose the small constant $c_1 > 0$ depending only on $n$ and the constants $\nu$, $c$ of the last section,  such that for $0 < r \leq r_0$, 
    \begin{equation} \label{7a2}
        \pi(x) \in S_r \cap E_\infty \ \text{ for } x\in S_r \cap W(c_1),
    \end{equation}
    \begin{equation} \label{7a3}
        \pi : W(c_1) \to E_\infty \ \text{ is $C$-Lipschitz,}
    \end{equation}
    and naturally
    \begin{equation} \label{7a4}
        \pi(x) = x  \ \text{ when } x \in E_\infty \cap B(0,r_0).
    \end{equation}
    The letter $C$ denotes a constant $\geq 1$ that depends only on $n$.
\end{pro}

So we decided to let our first projection act on spheres
separately (but notice that we require $\pi$ to be globally Lipschitz).
This is convenient because we don't need to know $E_\infty$ at scales much smaller
than $r$ to define $\pi$ on $S_r$.
Notice however that $W(c_1)$ is a small conical neighborhood of $0$, but does not contain 
a neighborhood of the origin. It would have been nice to have a projection $\pi$ defined on 
a whole neighborhood of $E_\infty$, but this will not happen, for simple topological reasons.
Even if $E_\infty$ is a plane through the origin, we cannot map the unit sphere to 
$E_\infty \cap S_1$ continuously (where do we send the poles?). This means that we will need something else than $\pi$ in later sections.

Let us use the description of Section \ref{S5}, that is, we cover $E_\infty \cap B(0,r_0)$ by boxes $A_r(2) \cap H(10c_i,\xi_i)$ (where $0 < r \leq r_0$) such that $E$ has an explicit description in $A_r \cap H(100 c_i,\xi_i)$, and we use each such box to define $\pi$ on $W(c_1) \cap A_r(2) \cap H(10c_i,\xi_i)$.
Our construction of $\pi$ will be intrinsic and will not depend on the choice of the covering.

We fix a radius $0 < r \leq r_0$ and a box $A_r \cap H(100 c_i,\xi_i)$ where $E_\infty$ is described 
by one of the cases of Section \ref{S5}.
To simplify the notation, 
we write $(c,\xi_0)$ for $(c_i,\xi_i)$.
We only detail how to define $\pi$ on $W(c_1) \cap S_r \cap H(10c,\xi_0)$, but the construction could be easily applied to the thicker domain $W(c_1) \cap A_r(2) \cap H(10c,\xi_0)$.
Similarly, the definition of $\pi$ on $W(c_1) \cap S_r$ could be done by using boxes associated to another radius in $(r/2,2r)$; it does not matter since our construction is intrinsic.

We let $\tau_0 > 0$ be a parameter that is small enough, depending on $n$.
Typically, we will work with $C \tau_0$-Lipschitz curves, 
and we will need $\tau_0$ small enough so that all the properties of Section~\ref{S6a} holds true.
We let $\tau > 0$ be the parameter used in Section~\ref{S5} in the definition of Cases 1--6. 
It is assumed to be very small compared to $\tau_0$ and $c$. Moreover, $c$ will also be 
assumed to be small compared to $\tau_0$.

We start with Case 2, if it exists (i.e., there exists a point in $S_r \cap E_\infty \cap \Gamma$ of type $\bV$, with a sharp enough angle). The following construction also applies directly to Case~2~Bis. 
Let us use the description \eqref{5b9}-\eqref{5b11} that we found in this region, which works in $A_r \cap H(100c,\xi_0)$.
Recall that (in this region) $E_\infty$ is composed of three relatively closed faces $F_i \subset A_r \cap H(100c,\xi_0)$, $i \in I = \set{v, +, -}$, bounded by a curve $G$.
Since $\Gamma$ is a $C \tau$-Lipschitz graph over $L_0$ which stays at distance $\leq C \tau r$ from $L_0$, the piece $\Gamma \cap A_r \cap H(100c,\xi_0)$ is transverse to the spheres.\footnote{
        We say that a continuous path $f : [t_0,t_1] \to \R^n$ is transverse to the spheres if $t \mapsto \abs{f(t)}$ is strictly increasing. 
        It suffices for instance that there exists $c > 0$ such that $f(s) \cdot (f(t) - f(s)) \geq c^{-1} \abs{f(s)} \abs{f(t) - f(s)}$ for all $s < t$.
        Indeed, this condition implies $\abs{f(t)} - \abs{f(s)} \geq c^{-1} \abs{f(t) - f(s)}$, which can be checked by computing $\tfrac{\dd}{\dd{x}} \abs{f(s) + x(f(t) - f(s))} \geq c^{-1} \abs{f(t) - f(s)}$ for $x \in [0,1]$. 
        This further shows that the map $\rho \mapsto f(s(\rho))$, where $s = s(\rho)$ is such that $\abs{f(s))} = \rho$, is $c$-Lipschitz.
        In case of a piece of graph $f(t) = (t,\psi(t))$ parametrized by $t \in [t_0,t_1]$, this sufficient condition follows easily if $\psi$ is $\tau$-Lipschitz with $\tau \leq 1/100$ and if $\abs{\psi(t_0)} \leq 10t_0$ (i.e., the vector $f(t_0)$ is not too orthogonal to the line of parametrization, unless $f(t_0) = 0$). Moreover in this case, the map $\rho \mapsto f(s(\rho))$ is $(1+C\tau)$-Lipschitz.
 }
 By transversality, there exists a unique point $\xi(r) \in \Gamma \cap S_r$ that lies on the same side as $\xi_0$ and the map $r \mapsto \xi(r)$ is $(1 + C\tau)$-Lipschitz.
The curve $\Gamma \cap A_r \cap H(100c,\xi_0)$, $G$ and the vertical face $F_v$ are entirely contained in $H(\nu c, \xi_0)$ and, in particular, we have $\abs{\xi(r) - r\xi_0} \leq \nu r c$.
Similarly, we call $\zeta(r)$ the unique point of $S_r \cap G$, and it satisfies $\abs{\zeta(r) - r \xi_0} \leq \nu r c$.
We define 
\begin{equation}
    \gamma_i := F_i \cap S_r
\end{equation}
so that $\gamma_i$ is a relatively closed subset of $S_r \cap B(r \xi_0,100 rc)$, and we have the decomposition
\begin{equation}
    E_\infty \cap S_r \cap B(r \xi_0, 100 rc) = \gamma_+ \cup \gamma_- \cup \gamma_v.
\end{equation}

\vskip-0.4cm
  \begin{figure}[!h]
  \centering
    \includegraphics[width=6cm]{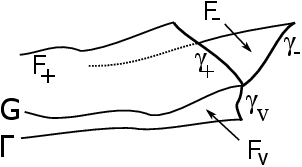}
\vskip-0.4cm
 \caption{Case 2, the intersection with $S_r = \d B(0,r)$}
\label{fig2b}
\end{figure}

We assume $c$ and $\tau$ small enough (depending on $n$, $\tau_0$) so that each $\gamma_i$ is a piece of $1$-dimensional $\tau_0$-Lipschitz graph above $\R \xi_i$, where we recall that $\xi_{\pm}$ and $\xi_v$ are defined in Case 2 of Section \ref{S5}.
For $i \in \set{\pm}$, the curve $\gamma_i$ starts from $\zeta(r)$ and goes all the way so $S_r \cap \partial B(r \xi_0,100 cr)$. 
On the other hand, the vertical curve $\gamma_v$ goes from $\xi(r)$ to $\zeta(r)$ and is very short; it is entirely contained in $S_r \cap B(r \xi_0, \nu cr) \subset S_r \cap B(\zeta(r), 2 \nu c r)$.
All we need to know about $\nu$ in this section is that $\nu \leq 1/100$ so that $\gamma_v$ is entirely contained in $S_r \cap B(\zeta(r),c r /10)$.
We also define the direction 
\begin{equation} \label{7a5}
    \theta(x) = \frac{x-\zeta(r)}{|x-\zeta(r)|} \in S_1
    \ \text{ for }x \in \R^n \sm \{ \zeta(r) \}.
\end{equation}

We will need to project on the $\gamma_i$.
We start by picking the unique point $\zeta_\pm(r)$
of $\gamma_\pm$ that lies at distance $c r$ from $\zeta(r)$, and then set
$e_\pm = \theta(\zeta_\pm(r))$; this will be the direction that we choose to define 
the projection onto $\gamma_\pm$. For the short $\gamma_v$, we simply take
$e_v = - (e_+ + e_-)/|e_+ + e_-|$ (this is more stable
because $\gamma_v$ may be very short).
For all $i \ne j$, we have at least $\abs{e_i - e_j} \geq 1/2$ (we can take $c$ and $\tau$ smaller if necessary).

    Let us remark that if $s \mapsto (s,\psi(s))$ is any $\tau_0$-Lipschitz graph above a line $\R e$, where $e = (1,0)$, then for all $x = (s,\psi(s))$ and $y = (t,\psi(t))$ with $t > s$, we have
    \begin{equation}
        \abs{\frac{y - x}{\abs{y - x}} - e} \leq 10 \tau_0.
    \end{equation}
    This uses the fact that $\abs{(y - x) - (t - s) e} \leq C \tau \abs{x - y}$ and the usual formula
    \begin{equation}\label{eq_theta_lip}
        \abs{\frac{u}{\abs{u}} - \frac{v}{\abs{v}}} \leq 2 \abs{u}^{-1} \abs{u - v} \quad \text{for all $u,v \in \R^n \setminus \set{0}$},
    \end{equation}
    applied to $u = (y - x)$ and $v = (t - s)e$.

    In our situation, $\gamma_{\pm}$ is a piece of $\tau_0$-Lipschitz graph above $\R \xi_{\pm}$ which starts from $\zeta(r)$, so $e_{\pm} \cdot \xi_{\pm} \geq 0$ and the above reasoning show 
     that $\abs{e_\pm - \xi_\pm} \leq 10 \tau_0$. As a consequence, $\gamma_{\pm}$ is a $C \tau_0$ Lipschitz graph above $\R e_{\pm}$. Using $\xi_v = -(\xi_+ + \xi_-)/\abs{\xi_+ + \xi_-}$ (and similarly for $e_v$), 
     one can also check that $e_{v} \cdot \xi_v \geq 0$ and $\abs{e_{v} - \xi_v} \leq C \tau_0$
    using (\ref{eq_theta_lip}), thus $\gamma_v$ is a $C \tau_0$-Lipschitz graph above $\R e_v$.

    The point of replacing $\xi_i$ by $e_i$ is to have a direction of projection onto $\gamma_i$ which is intrinsic.
    (notice that even though $e_i$ depends on $c = c_i$, Lemma \ref{lem_Xcover} guarantees that $c_i$ is always the same when Case 2 or Case 2 Bis applies).
    Indeed, there could be different choices for $X$ in (\ref{5a1}), which are $\varepsilon_0$ close to each other and which could present a slightly rotated $\bV$ at $\xi_0$. If the construction of $\pi$ on $S_r$ depends on the choice of $X$, it might not glue well when $r$ varies.

Now, we define open regions where we can easily project on $\gamma_i$, $i\in I$, by  
\begin{equation} \label{7a6}
    H_i = \Set{x \in S_r \cap B(r\xi_0,20cr) | x \ne \zeta(r),
    \ \abs{\theta(x) - e_i} < 10^{-7}}. 
\end{equation}
We will also need to define a projection in the larger convex domain (not a subset of $S_r$)
\begin{equation}
    \widehat{H}_i = \Set{x \in B(r\xi_0,20cr) | x \ne \zeta(r),\ \abs{\theta(x) - e_i} \leq 10^{-6}} \cup \set{\zeta(r)}.
\end{equation}
As $\gamma_i$ is a piece of $\tau_0$-Lipschitz graph above $\R e_i$ which starts from $\zeta(r)$, we can assume $\tau_0$ small enough so that for all $x \in \gamma_i \setminus \set{\zeta(r)}$, we have $\abs{\theta(x) - e_i} \leq 10^{-8}$, in particular $\gamma_i \cap B(\zeta(r),20cr) \setminus \set{\zeta(r)} \subset H_i$.

We justify that the regions $\widehat{H}_i$ are far from each other in the sense that for all $i \ne j$,
\begin{equation}\label{eq_distHi}
    \text{for all $x \in \widehat{H}_i$, \quad $\mathrm{dist}(x,\widehat{H}_j) \geq \abs{x - \zeta(r)}/2$.}
\end{equation}
For $i \ne j$, one can compute $\abs{\xi_i - \xi_j} \geq 3/2$ 
(they are unit vectors making an angle of 
 $2\pi/3$) and, since $\gamma_i$, $\gamma_j$ are $\tau_0$-Lipschitz graphs above $\R \xi_i$ and $\R \xi_j$ respectively, we deduce that $\abs{e_i - e_j} \geq 5/4$ (upon assuming $\tau_0$ small enough). The definition of $\widehat{H}_i$ and $\widehat{H}_j$ and the triangular inequality imply in turn that for all $x \in \widehat{H}_i$ and $y \in \widehat{H}_j$, we have $\abs{\theta(x) - \theta(y)} \geq 1$.

Then we use the usual formula (\ref{eq_theta_lip}) to deduce that for all $x \in \widehat{H}_i$ and for all $y \in \widehat{H}_j$ (with $i \ne j$)
\begin{equation}\label{eq_thetaxy}
    1 \leq \abs{\theta(x) - \theta(y)} \leq \frac{2 \abs{x - y}}{\abs{x - \zeta(r)}}
\end{equation}
and thus $\abs{x - y} \geq \abs{x - \zeta(r)}/2$. This proves our claim.

Next, we show that for all $x \in \widehat{H}_\pm$,
\begin{equation}\label{eq_H1}
    \mathrm{dist}(x, \gamma_\pm) \leq 10^{-5} \abs{x - \zeta(r)}
\end{equation}
and
\begin{equation}\label{eq_H2}
    E_\infty \cap S_r \cap B(x, 10^{-1} \abs{x - \zeta(r)}) = \gamma_\pm \cap B(x, 10^{-1} \abs{x - \zeta(r)}).
\end{equation}

    Before passing to the proof, let us explain the goal of these properties. For $x \in \widehat{H}_{\pm}$, the endpoints of $\gamma_\pm$ lie outside of the ball $B(x,10^{-1} \abs{x - \zeta(r)})$ so $\gamma_{\pm}$ coincides in $S_r \cap B(x,10^{-1} \abs{x - \zeta(r)})$ with a $C \tau_0$-Lipschitz graph above the line $\R e_{\pm}$ and so does $E_\infty \cap S_r$ by (\ref{eq_H2}). This will allow us to apply Section \ref{S6a} in order to make an intrinsic projection of a point $x \in \widehat{H}_{\pm} \setminus B(\zeta(r),cr)$ onto $\gamma_{\pm}$. Precisely, for such a point $x$, we have $\mathrm{dist}(x,\gamma_{\pm}) \leq 10^{-3} cr$ and $E \cap S_r$ coincides in $B(x,10^{-1}cr)$ with a $C \tau_0$-Lipschitz graph. This shows $p_{x,r}(z)$ is well-defined for all $x \in S_r \cap \widehat{H}_{\pm}$ and $z \in \widehat{H}_{\pm}$ (we refer to Section \ref{S6a} for the definition of $p_{x,r}$).

    The first point is based on the observation that any point $x \in \R^n \setminus \set{\zeta(r)}$ such that $\abs{\theta(x) - e_{\pm}} \leq 10^{-6}$ is relatively close to its projection on $\zeta(r) + \R^+ e_{\pm}$ in the sense that
    \begin{align}
        \mathrm{dist}(x, \zeta(r) + \R^+ e_{\pm}) &\leq \mathrm{dist}(x, \zeta(r) + \abs{x - \zeta(r)} e_{\pm})\\
        &= \abs{x - \zeta(r)} \abs{\theta(x) - e_{\pm}} \leq 10^{-6} \abs{x - \zeta(r)}\label{eq_gammaP}.
    \end{align}
    (The final estimate also holds true for $x = \zeta(r)$.)
    This applies in particular to the points of $\gamma_{\pm}$ and since $\gamma_{\pm}$ goes continuously from $\zeta(r)$ to $S_r \cap \partial B(\zeta(r),40cr)$, we deduce that its projection on $\zeta(r) + \R^+ e_{\pm}$ contains at least the interval $\zeta(r) + [0,30cr] e_{\pm}$.
    For $x \in \widehat{H}_{\pm}$, one can see that its orthogonal projection $z$ on $\zeta(r) + \R^+ e_{\pm}$ belongs to $\zeta(r) + [0,30cr] e_{\pm}$ (as $\abs{\theta(x) - e_{\pm}} < 1$, we have $\theta(x) \cdot e_{\pm} \geq 0$ so $(x - \zeta(r)) \cdot e_{\pm} \geq 0$, and it is also clear that $\abs{z - \zeta(r)} \leq \abs{x - \zeta(r)} \leq 30cr$) so $z$ coincides with the projection of a point $y \in \gamma_{\pm}$. 
    We finish the proof of (\ref{eq_H1}) by estimating $\abs{x - y} \leq 10^{-5} \abs{x - \zeta(r)}$. According to the triangular inequality, $\abs{x - y} \leq \abs{x - z} + \abs{y - z}$ and we already know that $\abs{x - z} \leq 10^{-6} \abs{x - \zeta(r)}$. In order to estimate $\abs{y - z}$, we observe first using (\ref{eq_gammaP}) that
    \begin{equation}
        \abs{y - \zeta(r)} \leq \abs{y - z} + \abs{z - \zeta(r)} \leq 10^{-6}\abs{y - \zeta(r)} + \abs{x - \zeta(r)}
    \end{equation}
    so $\abs{y - \zeta(r)} \leq 2 \abs{x - \zeta(r)}$ and then using (\ref{eq_gammaP}) again
    \begin{equation}
        \abs{y - z} \leq 10^{-6} \abs{y - \zeta(r)} \leq 2 \cdot 10^{-6} \abs{x - \zeta(r)}.
    \end{equation}
    This concludes the proof of (\ref{eq_H1}).

Next, (\ref{eq_H2}) comes from the fact that
\begin{equation} \label{7c22}
    E_\infty \cap S_r \cap B(r \xi_0,100 cr) = \gamma_+ \cup \gamma_- \cup \gamma_v 
\end{equation}
and that for $x \in \widehat{H}_\pm$, the ball $B(x, 10^{-1} \abs{x - \zeta(r)})$ is disjoint from the other 
$\widehat{H}_i$ (see \ref{eq_distHi}).

\ms
In order to project onto $\gamma_{\pm}$, we define a projection $p_\pm : \widehat{H}_\pm \to \gamma_\pm$ perpendicular to the direction of $e_{\pm}$ :
\begin{equation} \label{7a8}
    p_\pm(x) \in \gamma_\pm \ \text{ and } \ (p_\pm(x) - x) \perp e_\pm.
\end{equation}
This is well-defined by Section \ref{S6a}, (\ref{eq_H1}) and (\ref{eq_H2}) as usual.
Furthermore, the map $p_{\pm}$ is the identity on $\gamma_{\pm}$ and, 
by (\ref{eq_p}) and (\ref{eq_H1}),
\begin{equation} \label{7c24}
    \text{$\abs{p_{\pm} - \mathrm{id}} \leq 10^{-2} c r$ on $\widehat{H}_\pm$ \quad \text{and} \quad $p_{\pm}$ is $2$-Lipschitz.} 
\end{equation}
The projection $p_\pm$ is defined in particular at $\zeta(r)$ where it satisfies $p_\pm(\zeta(r)) = \zeta(r)$. As $p_\pm$ is $2$-Lipschitz (see Section \ref{S6a}, just above (\ref{eq_p})), we deduce that $\abs{p_\pm(x) - \zeta(r)} \leq 2 \abs{x - \zeta(r)}$ for all $x \in \widehat{H}_\pm$.
In order to project onto $\gamma_v$, we proceed slightly differently because $\gamma_v$ is short.
We extend $\gamma_v$ by a half line that starts from $\xi(r)$ and goes in the direction of $e_v$;
this gives an extended graph $\wh \gamma_v$ above $\R e_{\pm}$, we can define a projection
$\wh p_v: \widehat{H}_v \to \wh \gamma_v$ perpendicular to $e_v$ as above, 
and finally we define $p_v = \kappa \circ \wh p_v$, where  
$\kappa$ is the projection $\wh \gamma_v \to \gamma_v$ defined by $\kappa(x) = x$ when 
$x \in \gamma_v$ and $\kappa(x) = \xi(r)$ otherwise.   
Again $p_v : \widehat{H}_v \to \gamma_v$ is 
the identity on $\gamma_v$, is $2$-Lipschitz and satisfies $\abs{p_v(x) - \zeta(r)} \leq 2 \abs{x - \zeta(r)}$ for $x \in \widehat{H}_v$.

We will use the $p_i$ to define a projection
\begin{equation} \label{7c25}
    \pi : W(c_1) \cap S_r \cap B(r\xi_0,10cr) \to E_\infty \cap S_r,
\end{equation}
where $c_1 > 0$ depends on only on $c$.
Actually, we will define $\pi$ on 
$S_r \cap \left[H_+ \cup H_- \cup B(\zeta(r),cr)\right]$.

We are going to check first that, if $c_1$ is small enough compared to $c$, this domain contains $W(c_1) \cap S_r \cap B(r\xi_0,10cr)$, i.e,
\begin{equation} \label{7a10}
    W(c_1) \cap S_r \cap B(r\xi_0,10cr) \setminus B(\zeta(r),c r) \subset H_+ \cup H_-.
\end{equation}
We recall from (\ref{eq_distES0}) that for all $x \in S_r \cap B(r\xi_0, 20cr)$,
\begin{equation}\label{eq_distES}
    \mathrm{dist}(x, E_\infty \cap S_r) \leq C \mathrm{dist}(x,E_\infty),
\end{equation}
where $C \geq 1$ depends only on $n$.
Let us now see how to deduce (\ref{7a10}).
As $E_\infty \cap S_r \cap B(r\xi_0,100cr) = \gamma_+ \cup \gamma_{-} \cup \gamma_v$, (\ref{eq_distES}) implies that all the points of $W(c_1) \cap S_r \cap B(r\xi_0,10cr)$
are within distance $\leq C c_1 r$ of $\gamma_+ \cup \gamma_- \cup \gamma_v$. 
Remember that as $\gamma_v$ is short, it stays inside $B(\zeta(r), c r/10)$.
Then if $c_1$ is small enough compared to $c$, the points of $W(c_1) \cap S_r \cap B(r\xi_0, 10cr) \setminus B(\zeta(r), c r)$ are actually within $C c_1 r$ of $\gamma_+ \cup \gamma_- \setminus B(\zeta(r), c r/2)$.
We recall that as $\gamma_i$ is $C\tau_0$-Lipschitz graph above $\R e_i$, we have for all $x \in \gamma_i$, $\abs{\theta(x) - e_i} \leq C \tau_0 \leq 10^{-8}$ and that by (\ref{eq_theta_lip}), the function $x \mapsto \theta(x)$ is $C / (c r)$-Lipschitz outside $B(\zeta(r), c r/2)$. Hence, we can assume one more $c_1$ small enough so that if a point $x \in \R^n \setminus B(\zeta(r),cr)$ is at distance $\leq C c_1 r$ from $\gamma_i \setminus B(\zeta(r),c r/ 2)$, then $\abs{\theta(x) - e_i} < 10^{-7}$. This concludes the proof of (\ref{7a10}).

We now start the construction of $\pi$. 
We use a cut-off to single out the regions $H_i$, $i \in I$. 
Pick a Lipschitz cut-off function $\eta$ on $[0,1]$
such that $\eta(t) = 1$ when $0 \leq t \leq 10^{-7}$,
$\eta(t) = 0$ when $t \geq 10^{-6}$, and which interpolates linearly in-between.
Also set $\eta_i(x) = \eta(|\theta(x) - e_i|)$; this map is only locally Lipschitz away
from $\zeta(r)$, with $|\nabla \eta_i| \leq C |x - \zeta(r)|^{-1}$, but this will be compensated
in Lipschitz estimates because we will always multiply by terms like $x-\zeta(r)$.
Then set
\begin{equation} \label{7a9}
    \pi(x) = \zeta(r) \ \text{ for } x \in 
    S_r \cap B(\zeta(r), c r) \sm \cup_{i \in I} \widehat{H}_i
\end{equation}
and 
\begin{equation} \label{7a9bis}
    \pi(x) =  p_i \big[\eta_i(x) p_i(x) + (1-\eta_i(x)) \zeta(r) \big], 
    \ \text{ for } \ x \in S_r \cap  B(\zeta(r), c r) \cap \widehat{H}_i.
\end{equation}
The second formula is well-defined because for $x \in \widehat{H}_i \cap B(\zeta(r),cr)$, we have
\begin{equation}
    p_i(x) \in \gamma_i   \cap B(\zeta(r),2cr) \subset H_i \subset \widehat{H}_i
\end{equation}
and $\widehat{H}_i$ is convex.
This completes our definition of $\pi$ in $S_r \cap  B(\zeta(r), c r)$.
Observe that in each $H_{i}$,  $\eta_i = 1$  
so $\pi$ coincides there with $p_i$. The function $\pi$ built so far is $C$-Lipschitz
because of the remark on the Lipschitz character of $\eta_i(x)$ and the observation that 
$\abs{p_i(x) - \zeta(r)} \leq 2 \abs{x - \zeta(r)}$ in $\widehat{H}_i$.

Let us also note that for all $x \in S_r \cap B(\zeta(r), c r)$,
    \begin{equation}\label{eq_PIdistES}
    \abs{\pi(x) - x} \leq C \mathrm{dist}(x, E_{\infty} \cap S_r).
\end{equation}
This is clear if $x = \zeta(r)$ or if there exists $i$ such that $x \in H_i$ by properties of projections. Otherwise, $\abs{\theta(x) - e_i} \geq 10^{-7}$ for all $i$, and in this case, the main argument behind (\ref{eq_PIdistES}) is the fact that $\abs{x - \zeta(r)} \leq C \mathrm{dist}(x, E_{\infty} \cap S_r)$. Precisely, the distance from $x$ to $E_{\infty} \cap S_r$ is attained at a point $y \in E_{\infty} \cap S_r \cap B(r\xi_0,100cr)$ and there exists $i$ such that $\abs{\theta(y) - e_i} \leq 10^{-8}$ so $\abs{\theta(x) - \theta(y)} \geq 10^{-7}$. Reasoning as in (\ref{eq_thetaxy}), we get that 
$\abs{x - y} \geq 10^{-8} \abs{x - \zeta(r)}$.
    In order to deduce (\ref{eq_PIdistES}), observe finally that $\abs{\pi(x) - x} \leq \abs{\pi(x) - \zeta(r)} + \abs{\zeta(r) - x} \leq C \abs{x - \zeta(r)}$.

Our next step is to define $\pi$ in $H_\pm \sm B(\zeta(r), c r)$.
We start with the region $H_\pm \sm B(\zeta(r), 2c r)$.
According to Section \ref{S6a} and the comments below (\ref{eq_H1})-(\ref{eq_H2}), 
we can set in $H_\pm \sm B(\zeta(r),2 c r)$,
\begin{equation}\label{7a11}
    \pi(x) = p_{x,r}(x) \in \gamma_{\pm}.
\end{equation}
Moreover, in each domain $H_\pm \sm B(\zeta(r), 2 c r)$, $\pi$ is $C$-Lipschitz and, by (\ref{eq_p}) we have
\begin{equation}
    \abs{\pi(x) - x} \leq 10 \mathrm{dist}(x,E \cap S_r) \leq 10^{-2} c r.
\end{equation}

We then define $\pi$ in $H_\pm \cap B(\zeta(r), 2 c r) \setminus B(\zeta(r), c r)$.
There, we interpolate linearly the unit vector that we use to define the projection so that 
the various definitions match on the interfaces $\d B(\zeta(r), 2c r)$ and $\d B(\zeta(r), c r)$.
For $x \in H_\pm \setminus B(\zeta(r), c r)$, we let $e_\pm(x,r)$ be the unit vector generating $L_{x,r}$ 
(the direction of $p_{x,r}$), oriented in such a way that $e_{\pm}(x,r) \cdot e_{\pm} \geq 0$. 
In $H_\pm \cap B(\zeta(r), c r)$, we used the vector $e_\pm$ to project onto $E_\infty \cap S_r$
and in $H_{\pm} \setminus B(\zeta(r), 2c r)$, we used the vector $e_\pm(x,r)$.
In the intermediate region, we use
\begin{equation} \label{7a12}
    e_{\pm}(x) := \frac{e'_{\pm}(x)}{|e'_{\pm}(x)|},
\end{equation}
where
\begin{equation}
    e'_\pm(x) := \frac{2c r - \abs{x - \zeta(r)}}{c r} \, e_{\pm} + \frac{\abs{x - \zeta(r)} - c r}{c r} \, e_{\pm}(x,r).
\end{equation}
Notice that even though the direction of the projection depends (slowly) on $x$,
this does not prevent us from using it.
We do not need any injectivity property of the projections anyway.
Writing the definition \eqref{7a11} explicitly in terms of $E_\infty \cap S_r$ has the advantage of gluing well across the $S_r$.

Let us check that in the domain $H_{\pm} \cap B(\zeta(r),2 c r) \setminus B(\zeta(r), c r)$,
\begin{equation}\label{eq_ei}
    \text{$\abs{e_{\pm}(x) - e_{\pm}} \leq C \tau_0$ \quad and \quad $x \mapsto e_{\pm}(x)$ is $C/(cr)$-Lipschitz.}
\end{equation}
Since each point $x \in H_{\pm} \cap B(\zeta(r),2 c r) \setminus B(\zeta(r), c r)$ is at 
distance $\leq 10^{-3} c r$ from $\gamma_{\pm}$, 
and since $\gamma_{\pm}$ coincides with a $C\tau_0$-Lipschitz graph above $\R e_{\pm}$ in $B(x, 10^{-1} c r)$, (\ref{eq_ei}) allows us to define $\pi(x)$ in $H_{\pm} \cap B(\zeta(r),2 c r) \setminus B(\zeta(r), c r)$ as the projection onto $\gamma_{\pm}$ in the direction orthogonal to $e_{\pm}(x)$. Moreover, (\ref{eq_ei}) shows that $\pi$ is $C$-Lipschitz in $H_{\pm} \cap B(\zeta(r),2 c r) \setminus B(\zeta(r), c r)$; this uses also (\ref{eq_pdist}) and the fact that we have $\mathrm{dist}(x, \gamma_{\pm}) \leq 10^{-3} cr$ in this domain.
    
So let's justify (\ref{eq_ei}).
We know from the definition of $e_{\pm}(x,r)$ and (\ref{eq_Lxy0})-(\ref{eq_Lxy}) that
that 
\begin{equation}
    \text{$\abs{e_{\pm}(x,r) - e_{\pm}} \leq C \tau_0$ \quad and \quad $\abs{e_{\pm}(x,r) - e_{\pm}(y,r)} \leq \frac{C}{r} \abs{x - y}$.}
\end{equation}
Then it follows immediately from the formula of $e'_i$ that 
\begin{equation}
    \text{$\abs{e_{\pm}'(x) - e_{\pm}} \leq C \tau_0$ \quad and \quad $\abs{e'_{\pm}(x) - e'_{\pm}(y)} \leq \frac{C}{cr} \abs{x - y}$}.
\end{equation}
We then deduce (\ref{eq_ei}) using the formula (\ref{eq_theta_lip}).

We are done with the definition of $\pi$ in $H_{\pm} \cap B(\zeta(r),2 c r) \setminus B(\zeta(r), c r)$. By properties of projections, we have as usual
\begin{equation}\label{eq_pi12}
    \abs{\pi(x) - x} \leq 10 \mathrm{dist}(x,E \cap S_r) \leq 10^{-2} c r
\end{equation}
in this domain.

    Finally, we check that 
    \begin{equation}
        \text{$\pi$ is $C$-Lipschitz on $S_r \cap \left[H_+ \cup H_- \cup B(\zeta(r), c r)\right]$.}
    \end{equation}
    By construction, the map $\pi$ is independently $C$-Lipschitz on  
    $H_+ \setminus B(\zeta(r),2cr)$ and $H_+ \cap B(\zeta(r),2cr) \setminus B(\zeta(r),cr)$. 
    In order to show  that $\pi$ is $C$-Lipschitz in their union $H_+ \setminus B(\zeta(r),cr)$, we observe that for any $x \in H_+ \setminus  B(\zeta(r),2cr)$ and $y \in H_+ \cap B(\zeta(r),2cr) \setminus B(\zeta(r),cr)$, the geometry of these domains allows to find $z \in H_+ \cap \partial B(\zeta(r),2cr)$ such that $\abs{x - z} + \abs{z - y} \leq C \abs{x - y}$ and thus
    \begin{align}
        \abs{\pi(x) - \pi(y)} &\leq \abs{\pi(x) - \pi(z)} + \abs{\pi(z) - \pi(y)}\\
        &\leq C \abs{x - z} + C \abs{z - y} \leq C \abs{x - y}.
    \end{align}
    As $\abs{\pi - \mathrm{id}} \leq 10^{-2} c r$ in both $H_+ \setminus B(\zeta(r), c r)$ and $H_- \setminus B(\zeta(r),cr)$, and since these two domains are at distance $\geq cr/10$ from each other, the map $\pi$ is again $C$-Lipschitz in their union. The only case left is to check whether $\pi$ is Lipschitz in the union of $H_{\pm} \setminus B(\zeta(r),cr)$ and $S_r \cap B(\zeta(r),cr)$. The argument is the same as before: $\pi$ glues continuously across their interface and for any $x \in S_r \cap H_\pm \setminus B(\zeta(r), cr)$ and $y \in S_r \cap B(\zeta(r),c r)$, one can find a point $z \in S_r \cap H_\pm \cap \partial B(\zeta(r),cr)$ such that $\abs{x - z} + \abs{z - y} \leq C \abs{x - y}$.

As mentioned at the beginning of the proof, the construction is the same for every other radius in $(r/2,2r)$. 
Thus, $\pi$ is $C$-Lipschitz along each piece of sphere $W(c_1) \cap S_r(10) \cap B(t\xi_i,10c_it)$ for $t \in (r/2,2r)$. 
One can also check that $\pi$ is $C$-Lipschitz across the spheres, that is, in the whole box $W(c_1) \cap A_r(2) \cap H(10c_i,\xi_i)$; this relies in particular on the Lipschitzness of $\rho \mapsto \xi(\rho), \zeta(\rho)$ (see the footnote at the beginning of this section), the Lipschitzness $r \mapsto e_i$ for each $i \in I$ (this type of property will be detailed in the next section, so we postpone the details until that point) and the Lipschitzness of $p_{x,r}$ (see (\ref{eq_pLipschitz})).

\ms

This ends the definition of $\pi$ near $r\xi_0$ when $\xi_0$ belongs to Case 2 (sharp $\bV$) or Case 2 Bis (truncated $\bY$). 
Fortunately for the other cases, the definitions will be the same, but simpler. 
Let us recall that the family of spherical caps $S(\xi_i,c_i)$ satisfy the two additional assumptions of Lemma \ref{lem_Xcover}. This makes sure that the boxes $A_r \cap H(100c_i,\xi_i)$, where $E_\infty$ is close to a $\bH$, $\bV$, $\bY$ or a truncated $\bY$ are disjoint. Therefore, the constructions of $\pi$ in these cases will be independent from each other.
However, a box $A_r \cap H(100c_i,\xi_i)$ can meet another box $A_r \cap H(100c_k, \xi_k)$, where $E_\infty$ is close to a $\bP$, but in this case $\xi_k \notin S(\xi_i,5c_i)$ and $c_k \leq c_i/10$ so the overlapping regions $H(10c_i,\xi_i) \cap H(10c_k,\xi_k)$ is outside $H(4c_i,\xi_i)$. This observation will allow the reader to check easily that all definitions coincide when two box meet.

We start with Case 3, with a generic cone $V$ of type $\bV$ with angle $\alpha$ in $(2\pi/3 + \varepsilon_1/400, \pi - \varepsilon_1/400)$, where $\varepsilon_1$ is defined in Section \ref{S5}.
There we can just use the same formula as in Case~2; we just have to consider that $F_v$ was essentially empty, $\xi(r) = \zeta(r)$ and we ignore the region $H_v$.

    When $V$ in Case 3 is nearly flat, i.e when $\alpha' := \pi - \alpha$ is small, we modify the construction a little bit, because we want to organize a smooth transition to the flat case where we proceed slightly differently.

    Let $\tau_1 \in (0,1)$ be such that previous constructions work for all $\tau \leq \tau_1$ (it only depends on $n$, $\nu$, $c$).
    We assume the constant $\varepsilon_1$ in Section \ref{S5} small enough so that $\varepsilon_1/400 \leq \tau_1$ and in particular, we have $\alpha' \leq \tau_1$ in Case 1.
    For both Case 1 (a plane $\bP$) and in Case 3 with $\alpha' \leq \tau_1$, we set $\pi$ with the same formula
    \begin{equation}\label{7a11bis}
        \pi(x) =  p_{x,r}(x)  \quad \text{in $W(c_1) \cap S_r \cap B(r\xi_0,10cr)$,}
    \end{equation}
    as in \eqref{7a11}, where now we can use the fact that $E_\infty \cap B(r\xi_0,100 rc)$ is a single $1$-dimensional $C \tau_1$-Lipschitz graph to define $p_{x,r}(x)$.
    Here and hereafter, we assume $c_1$ small enough compared to $c$ so that every point in $W(c_1) \cap S_r \cap B(r\xi_0,10cr)$ is at distance $\leq 10^{-3} cr$ from $E_\infty \cap S_r$ (see (\ref{eq_distES0})); this makes sure that $p_{x,r}(x)$ is well-defined.


    In Case 3 with $\tau_1 \leq \alpha' \leq 2 \tau_1$, we interpolate: the projection $\pi$ built in Case 3 is replaced by
    \begin{equation}\label{7a13}
        p_{x,r}((2 - \alpha'/\tau_1) x + (\alpha'/\tau_1 - 1) \pi(x)).
    \end{equation}
    This is well-defined for all $x \in W(c_1) \cap S_r \cap B(r\xi_0,10cr)$ because the point $z = (2 - \alpha'/\tau_1) x + (\alpha'/\tau_1 - 1)$ satisfies
    \begin{align}
        \mathrm{dist}(z, E \cap S_r) &\leq \mathrm{dist}(z, \pi(x))\\
                                     &\leq \abs{x - \pi(x)}\\
                                     &\leq C \mathrm{dist}(x,E_{\infty} \cap S_r)\leq C \mathrm{dist}(x,E)
    \end{align}
    and $\mathrm{dist}(x,E) \leq c_1 r \leq 10^{-3} cr$.

    Notice that as $\alpha$ goes from $2\tau_1$ to $\tau_1$, the role of $\zeta(r) = \xi(r)$ disappears in the construction. This is needed since in Case 1, the curve $\Gamma$ is allowed to leave $E_\infty$ so the point $\zeta(r) \in \Gamma \cap S_r$ may not belong to $E_\infty$ and cannot be used to project anymore.

One could object that $\alpha = \alpha(r)$ was not defined intrinsically (it depends on the chosen blow-up $X$ in (\ref{5a1})), but for the purpose of this construction, we can replace $\alpha$ with the computable number $\Angle(e_+,e_-) \in (0,\pi)$.
The fact that $E_\infty$ is extremely close to a $\bV$-cone
in $A_r \cap H(100 c,\xi_0)$ is useful to get the good description of $E_\infty$, but 
we do not need all that precision for the definition of the projection $\pi$ and the verification that it is Lipschitz.

    Next we consider Case 4, where $E_\infty$ is approximated by a half plane
    $H \in \bH$. In this case, $E_\infty \cap S_r \cap B(r \xi_0, 100 cr)$ is a single curve
    $\gamma = \gamma_+$ (compare with \eqref{7c22}). We can proceed as in Case 2, with $\xi(r) = \zeta(r)$ and ignoring $H_-$ and $H_v$.

Now consider Case 6, that is, the case of $\bY$-points in $A_r \cap H(100c,\xi_0)$.
We use the $\bY$-points of $S_r \cap E_\infty$ as new points $\xi(r) = \zeta(r)$, and perform the same construction as in Case 2 (this time, with three long curves $\gamma_i$).
This gives a definition of $\pi$ on the $W(c_1) \cap S_r \cap B(r\xi_0, 10c r)$, and as above we make sure to use \eqref{7a11} outside $B(\zeta(r), 2 c)$.
Now although on $S_r$ the points $\xi_0$ with a singularity of type 
$\bY$,  and Case 2 bis are far from each other, as $r$ varies they can be transformed in 
each other, and this  is why we try to use intrinsic formulae. For the present case, the reader may worry that we start with a $\bY$ cone, which slowly evolves into the previous case 
of a truncated $\bY$, with a third leg that becomes shorter and shorter. In the first case,
we used an intrinsic formula based on the third curve to define the direction of the associated projection, while in the case of a short leg, we used the formula $e_v = - (e_+ + e_-)$ which was more stable. Let us simply decide, depending on the (intrinsically evaluable) length
of the short leg, to use a formula that interpolates nicely between the two cases; we skip the
formula because it may only add to the confusion; observe that we did something similar
for the transition between a plane and a flat $\bV$.

We are left with regions where $E_\infty$ is well approximated by planes (Case 5), and there, as promised above,
we simply use \eqref{7a11bis}. Fortunately, we always made sure to use \eqref{7a11} away from $\xi_i$. This way, all our definitions glue nicely with each other.

    At this point, we have defined a $C$-Lipschitz projection 
 $\pi : W(c_1) \cap A_r(2) \cap H(10c_i,\xi_i) \to E_\infty$ for each radius $0 < r \leq r_0$ and for each box $A_r \cap H(100c_i,\xi_i)$ associated to scale $r$ in Section \ref{S5}. Moreover, all definitions of $\pi$ coincide in overlapping regions. It is left to check that gluing everything induces a $C$-Lipschitz projection $\pi : W(c_1) \to E_\infty$.

    We check first that for any fixed radius $0 < r \leq r_0$, $W(c_1) \cap A_r(2)$ is covered by the charts $A_r(2) \cap H(10c_i,\xi_i)$ associated to scale $r$. This ensures that gluing everything provides a mapping defined on the whole $W(c_1)$.
    We recall from Section \ref{S5} that $X \cap \d B(0,1)$ is covered by the spherical caps $S(\xi_i,5c_i)$ with $c_i \geq 10^{-3} \nu c$. Assuming $\varepsilon_0$ small enough (depending on $\nu$ and $c$), it follows that $E_\infty \cap A_r(2)$ is covered by the conical domains $H(8 c_i,\xi_i)$ (this is only a slightly more precise variant of what we did near (\ref{eq_Xcover})).
    For $x \in W(c_1) \cap A_r(2)$, there exists $x' \in E_\infty$ such that $\abs{x - x'} \leq c_1 \abs{x}$ and taking in particular $c_1 \leq 10^{-4} \nu c$, we have
    \begin{equation}
        \abs{\frac{x}{\abs{x}} - \frac{x'}{\abs{x'}}} \leq 2 \abs{x}^{-1} \abs{x - x'} \leq 10^{-3} \nu c.
    \end{equation}
    Since $x'$ belong to one of the $H(8c_i,\xi_i)$ and $c_i \geq 10^{-3} \nu c$, it follows that $x$ belongs to $H(9c_i,\xi_i)$, proving our claim.

    We now prove that $\pi$ is $C$-Lipchitz in $W(c_1)$. 
    For $x, y \in W(c_1)$ such that one of them is zero, say $x = 0$, we clearly have $\abs{\pi(x) - \pi(y)} \leq \abs{x - y}$ because $\pi(0) = 0$ and $\pi$ preserve spheres.
    Next, we focus on the case where $x,y \in W(c_1) \setminus \set{0}$ are such that $\abs{x - y} \leq 10^{-4} \nu c \abs{x}$. One can see that $x,y \in W(c_1) \cap A_r(2)$, where $r = \abs{x} > 0$, and that
    \begin{equation}
        \abs{\frac{x}{\abs{x}} - \frac{y}{\abs{y}}} \leq 2 \abs{x}^{-1} \abs{x - y} \leq 10^{-3} \nu c.
    \end{equation}
    Therefore, there exists a chart $A_r \cap H(100c_i,\xi_i)$ associated to scale $r$ such that $x \in H(9c_i,\xi_i)$ and consequently $y \in H(10c_i,\xi_i)$. Then we can just use the fact that $\pi$ is $C$-Lipschitz in $A_r(2) \cap H(10c_i,\xi_i)$.
    We can deal with the case where $\abs{x - y} \leq 10^{-4} \nu c \abs{y}$ similarly.

    The last case, where $\abs{x - y} \geq 10^{-4} \nu c \max(\abs{x},\abs{y})$ relies on the fact that for all $x \in W(c_1)$,
    \begin{equation}\label{eq_distPIE}
        \abs{\pi(x) - x} \leq C \mathrm{dist}(x,E_\infty). 
    \end{equation}
    This can be checked independently for each Case 1--6, but this also comes from a more generic argument.
    Considering the scale $r = \abs{x}$, we have seen above that $x$ belongs to one of the spherical cap $S_r \cap B(r\xi_i,9c_i r)$ associated to $r$. Moreover, (\ref{eq_distES0}) show that we can find a point $y \in E_\infty \cap S_r \cap B(r\xi_i, 100 c_i r)$ such that $\abs{x - y} \leq C \mathrm{dist}(x,E_\infty) \leq C c_1 r$.
    Taking $c_1$ small enough compared to $n$ and $c$, we have in particular $y \in E_\infty \cap S_r \cap B(r\xi_i,10c_i r)$. Now, $\pi : W(c_1) \cap S_r \cap B(r\xi_i,10rc_i)$ is $C$-Lipschitz and coincides with the identity map on $E_\infty$, so
    \begin{equation}
        \abs{\pi(x) - x} \leq \abs{\pi(x) - y} + \abs{x - y} \leq C \abs{x - y} \leq C \mathrm{dist}(x,E_\infty),
    \end{equation}
    which proves our claim.

    Coming back to proof of Lipschitzness of $\pi$ when $x, y \in W(c_1)$ are such that $\abs{x - y} \geq 10^{-4} \nu c \max(\abs{x},\abs{y})$, we estimate using (\ref{eq_distPIE}) that
    \begin{align}
        \abs{\pi(x) - \pi(y)}   &\leq \abs{\pi(x) - x} + \abs{\pi(y) - y} + \abs{x - y}\\
                                &\leq C c_1 (\abs{x} + \abs{y}) + \abs{x - y} 
        \leq 2 \abs{x - y}, \nonumber 
    \end{align}
    provided that $c_1$ is small enough depending on $\nu$ and $c$.
    This finishes the proof.

Recall that we cannot use $\pi$ directly as a retraction to $E_\infty$, because $W(c_1)$
is not a neighborhood of $0$. Since we do not know where to send the points of
$B(0,r_0) \sm W(c_1)$, we will decide to send them to the origin, and interpolate on part of  
$W(c_1)$. For this a construction of a contraction of $E_\infty \cap B(0,r_0)$ will be useful.

\section{Contraction of \texorpdfstring{$E_\infty$}{Einfty} to the origin} 
\label{S8a}

We keep the assumptions of the previous sections, and in particular we fix 
$x_0 \in E_\infty$, $r_0 > 0$, and $\varepsilon_0$ as in Section \ref{S5}.
Let us say that $x_0 = 0$.
For $x\in E_\infty \cap B_0$, where $B_0 = B(0, r_0)$, we now want to define a path in $E_\infty$ that 
goes from $0$ to $x$, and depends on $x$ in a Lipschitz way.

\begin{pro}\label{l8a1}
    There is a Lipschitz mapping 
    $\sigma : \left[E_\infty \cap B_0\right] \times [0,1] \to E_\infty \cap B_0$, 
    such that
    \begin{equation}\label{8a0}
        |\sigma(x,t) - \sigma(y,s)| \leq C |x - y| + C \min(r(x),r(y)) | s - t |,
    \end{equation}
    where $r(x) = \abs{x}$, and
    \begin{equation} \label{8a1}
        \sigma(x,t) \in E_\infty \cap S_{tr} \ \text{ when }
        x \in E_\infty \cap S_r,
    \end{equation}
    and of course 
    \begin{equation} \label{8a2}
        \sigma(x,0) = 0 \ \text{ and } \sigma(x,1) = x.
    \end{equation}
    The letter $C$ denotes a constant $\geq 1$ that depends only on $n$.
\end{pro}

\ms
Then, by \eqref{8a1}, $|\sigma(x,t)| = t \abs{x}$ and in particular $\sigma(0,t) = 0$.
We will see other properties of $\sigma$ along the way.
It would  have been
nice to make the mapping $x \to \sigma(x,t)$ injective, 
but this cannot be arranged in general because of our Case 2 (see below).

We shall first define $\sigma(x,t)$ for $1/2 \leq t \leq 1$, and we shall see later how to compose and extend $\sigma$ to $0 \leq t <1/2$.
As usual we start by covering $E_\infty \cap B(0,r_0)$ with boxes $A_r(2) \cap H(10c_i,\xi_i)$ (where $0 < r \leq r_0$) such that $E$ has an explicit description in $A_r \cap H(100 c_i,\xi_i)$, and we use each such box to define $\sigma(x,t)$ for $x \in E_\infty \cap A_r(2) \cap H(10c_i,\xi_i)$ and $1/2 \leq t \leq 1$.

Let a radius $0 < r \leq r_0$ and a box $A_r \cap H(100 c_i,\xi_i)$ be given.
To simplify the notations, we write $(c,\xi_0)$ for $(c_i,\xi_i)$.
We only detail the construction of  $\sigma(x,t)$ for $x \in E_\infty \cap S_r \cap B(r \xi_0, 10cr)$
but it can be easily adapted to the thicker domain $E_\infty \cap A_r(2) \cap H(10c,\xi_0)$. Using boxes associated to another radius in $(r/2,2r)$ would work as well.
Although we construct $\sigma$ separately in the regions singled out in Section \ref{S5}, we will make sure that our definitions will be easy to glue.

We focus on the most interesting Case 2 (a sharp $\bV$-set).
The construction will also apply directly to Case 2 Bis (a truncated $\bY$ with a triple junction $\zeta_0 \in \partial B(0,1) \cap B(\xi_0,\nu c) \setminus \set{\xi_0}$). In order, to avoid distinguishing cases, we simply set $\zeta_0 = \xi_0$ in Case 2 and the proof will apply to both cases.
The set $E_{\infty}$ is composed of three relatively closed faces 
$F_i \subset A_r \cap H(100c,\xi_0)$, $i \in I = \set{v, +, -}$, precisely,
\begin{equation}
    E_\infty \cap A_r \cap H(100 c,\xi_0) = F_v \cup F_+ \cup F_-,
\end{equation}
where $F_+$, $F_-$ are two $\tau$-Lipschitz graphs bounded by $G$ and $F_v$ is a vertical wall between $\Gamma$ and $G$.
The set $\Gamma $ coincides in $A_r$ with the graph over the line $L_0$ (the line generated by $\xi_0$) of a $\tau$-Lipschitz function $\psi^0$ such that $\abs{\psi^0} \leq C \tau r$.
Similarly, $G$ is given in $A_r \cap H(100c,\xi_0)$ as a $C \tau$-Lipschitz graph above $L$ (the line generated by $\zeta_0$) of a $\tau$-Lipschitz function $\psi$ such that $\abs{\psi} \leq C \tau r$.
As $\tau$ is small enough, both curves are transverses to spheres. For all $\rho \in (r/2,2r)$, we let $\xi(\rho)$ and $\zeta(r)$ denote the unique points of $\Gamma \cap S_{\rho}$ and $G \cap S_{\rho}$ respectively which lie on the same side as $\xi_0$.
We have $\abs{\xi(\rho) - \rho \xi_0} \leq C \tau \rho$ and $\abs{\zeta(\rho) - \rho \zeta_0} \leq C \tau \rho$.
Besides, the maps $\rho \mapsto \xi(\rho)$ and $\rho \mapsto \zeta(\rho)$ are $(1+C\tau)$-Lipschitz 
on $(r/2,2r)$ (we refer to the footnote at the beginning of Section \ref{S7}).
The curves $\Gamma \cap A_r \cap H(100c,\xi_0)$ and $G$ and the vertical face $F_v$ are entirely contained in $H(\nu c, \xi_0)$.
We thus have $\abs{\zeta(\rho) - \rho \xi_0} \leq \nu c \rho$ and $F_v \cap S_r \subset B(r\xi_0,\nu c \rho)$ and since $\nu \leq 10^{-2}$, we deduce that $F_v \cap S_{\rho} \subset S_{\rho} \cap B(\zeta(r),c\rho/10)$.
We also define
\begin{equation}
    \gamma_i := F_i \cap S_r, 
\end{equation}
which is piece a $C \tau_0$-Lipschitz graph above $\R \xi_i$ (we recall $\xi_{\pm}$ and $\xi_v$ are defined in Case 2 of Section \ref{S5})
This allows to decompose $E_\infty \cap S_r \cap B(r \xi_0, 100 cr)$ as three $C \tau_0$-Lipschitz graphs;
\begin{equation}
    E_\infty \cap S_r \cap B(r \xi_0, 100cr) = \gamma_+ \cup \gamma_- \cup \gamma_v.
\end{equation}
The short vertical curve $\gamma_v$ connects the point $\xi(r) \in \Gamma \cap S_r$ to the $\bY$-point $\zeta(r) \in G \cap S_r$, and it is entirely contained in $S_r \cap B(r \xi_0,\nu c r)$.
On the other hand, the two curves $\gamma_{\pm}$ go all the way from $\zeta(r)$ to $S_r \cap \partial B(r\xi_0,100cr)$.

Because of the sliding condition, we want to choose $\sigma$ such that
\begin{equation} \label{8a3}
    \sigma(x,t) \in \Gamma \ \text{ when } x \in \Gamma,\ 1/2 \leq t \leq 1, 
\end{equation}
which will force us to do something slightly strange with $G$ and the faces.
In the present case, \eqref{8a3} just means that $\sigma(\xi(r),t) = \xi(tr)$ for $1/2 \leq t \leq 1$.

    Along the way, we shall need to know that
    \begin{equation}\label{eq_xirho}
        \text{$\rho \mapsto \frac{\xi(\rho)}{\rho}$ is $C \tau r^{-1}$-Lipschitz for $\rho \in (r/2,2r)$.}
    \end{equation}
    The first step in the proof of (\ref{eq_xirho}) is to show that
    \begin{equation}\label{eq_xi_claim}
        \text{$\rho \mapsto \xi(\rho) - \rho \xi_0$ is $C \tau$-Lipschitz for $\rho \in (r/2,2r)$.}
    \end{equation}
    To simplify the notations, we assume that $\xi_0$ is the first vector of the canonic base of $\R^n$. Then for $\rho \in (r/2,2r)$, there exists a unique $s = s(\rho) \geq 0$ such that $\xi(\rho) = (s,\psi^0(s))$. As $\abs{\xi(\rho)} = \rho$, it satisfies $\abs{s - \rho} \leq \abs{(s,0) - \xi(\rho)} \leq \abs{\psi^0(s)}$ and thus $s \in (r/4,4r)$. Observing that
    \begin{equation}
        \xi(\rho) - \rho \xi_0 = (s(\rho) - \rho, \psi^0(s(\rho))),
    \end{equation}
    it suffices to prove that $\rho \mapsto s(\rho) - \rho$ is $C \tau$-Lipschitz in order to deduce (\ref{eq_xi_claim}). Indeed, this implies that $\rho \mapsto s(\rho)$ is $2$-Lipschitz and then that $\rho \mapsto \psi^0(s(\rho))$ is $C \tau$-Lipschitz. Letting $\rho_1 < \rho_2$ in $(r/2,2r)$ and $s(\rho_1), s(\rho_2) \in (r/4,4r)$ the associated coordinate, one can compute that
    \begin{equation}\label{eq_lip_rho}
        \abs{(\rho_1 - s(\rho_1)) - (\rho_2 - s(\rho_2))} \leq C \tau \abs{s(\rho_1) - s(\rho_2)}.
    \end{equation}
    If $\tau$ is small enough, this implies in particular that, $\abs{s(\rho_1) - s(\rho_2)} \leq 2 \abs{\rho_1 - \rho_2}$ and plugging this in (\ref{eq_lip_rho}) yields
    \begin{equation}
        \abs{(\rho_1 - s(\rho_1)) - (\rho_2 - s(\rho_2))} \leq C \tau \abs{\rho_1 - \rho_2},
    \end{equation}
    justifying our claim. Combining (\ref{eq_xi_claim}) and the fact that $\abs{\xi(\rho) - \rho \xi_0} \leq C \tau \rho$ for all $\rho \in (r/2,2r)$, we deduce that $\rho \mapsto \rho^{-1}(\xi(\rho) - \rho \xi_0)$ is $C \tau$-Lipschitz and (\ref{eq_xirho}) follows easily.

Similarly, the map $\rho \mapsto \rho^{-1} \zeta(\rho)$ is $C \tau$-Lipschitz on $(r/2,2r)$ with the same argument, but replacing $\xi_0$ by $\zeta_0$.

For $x \in E_\infty \cap S_r \cap B(r \xi_0, 10cr)$, set
\begin{equation} \label{8a4}
    v(x) := \dist(x, \xi(r(x)) \in [0,20 cr] 
\end{equation}
where we sometimes write $r(x) = |x|$ instead of just $r$, to insist on the fact that it is a Lipschitz function of $x$.
We also set
\begin{equation}
    v_{max}(r) := \mathrm{dist}(\xi(r), \zeta(r)) \in [0, 2 \nu c r] 
\end{equation}
and we write simply $v_{max}$ when there is no ambiguity.
Let us justify that
\begin{equation}\label{eq_vmax}
    \abs{v_{max}(tr) - tv_{max}(r)} \leq C \tau (1 - t) r.
\end{equation}
We have
\begin{align}
    \abs{v_{max}(tr) - tv_{max}(r)}
    &
    \leq \big| \abs{\xi(tr) - \zeta(tr)} - t \abs{\xi(r) - \zeta(r) \big|
    }
    \nonumber       \\
    &\leq \abs{\xi(tr) - \zeta(tr) - t(\xi(r) - \zeta(r))}
    \nonumber       \\
    &\leq \abs{\xi(tr) - t\xi(r)} + \abs{\zeta(tr) - t \zeta(r)}
    \leq C \tau (1 - t) r,
\end{align}
where the last line comes from (\ref{eq_xirho}).

The curve $\gamma_v$ is transverse to spheres centered %
at $\xi(r)$ since it is a $C \tau_0$-Lipschitz graph starting from $\xi(r)$ (we refer to the footnote for details as the beginning of Section \ref{S7}). Moreover, we have for all $x \in \gamma_v$,
\begin{equation}\label{eq_vlessvmax}
    v(x) \leq v_{max}(r) \leq 2 \nu c r.
\end{equation}
The curve $\gamma_+$ is also transverse to the spheres centred at $\xi(r)$ because it it a $C \tau_0$-Lipschitz graph above $\xi_{\pm}$ starting from $\zeta(r)$ and because, observing that
\begin{equation} \label{8b18}
    \abs{\frac{\zeta(r) - \xi(r)}{\abs{\zeta(r) - \xi(r)}} - \xi_v} \leq C \tau_0 \leq 1/100,
\end{equation}
the vector $\zeta(r) - \xi(r)$ is not too orthogonal to $\R \xi_+$ (we refer again to the footnote).
By concatenation, the curve $\gamma_+ \cup \gamma_v$ is then again transverse to the spheres and goes continuously from $\xi(r)$ all the way to $S_r \cap \partial B(r\xi_0,100cr)$. 
We deduce that for all $0 \leq \rho \leq 20 cr$, 
there exists a unique point $x = x_+(\rho) \in \gamma_{+} \cup \gamma_v$ (resp. $x_-(\rho) \in \gamma_- \cup \gamma_v$) such that $v(x) = \rho$. Moreover, the mapping $\rho \mapsto x_{\pm}(\rho)$ is $(1 + C\tau_0)$-Lipschitz, and we have $\rho \leq v_{max}(r)$ if and only if $x_{\pm}(\rho) \in \gamma_v$.

For $x \in E_\infty \cap S_r \cap B(r\xi_0, 10cr)$ and $1/2 \leq t \leq 1$, we set
\begin{equation} \label{8a5}
    v_t(x) := \max(0, t v(x) - C_0 \tau (1 - t) r(x)), 
\end{equation}
where $C_0 \geq 1$ is a constant that we will fix soon (depending only on $n$).
This formula is chosen so that the following properties hold:
\begin{equation}\label{eq_vt1}
    \text{$v_t(x) = 0$ when $v(x) = 0$,}
\end{equation}
\begin{equation}\label{eq_vt2}
    \abs{v_t(x) - t v(x)} \leq C_0 \tau (1 - t) r  
\end{equation}
and, assuming that $C_0$ is bigger than the constant $C$ in \eqref{eq_vmax},
\begin{equation}\label{eq_vt3}
    \text{$v_t(x) \leq v_{max}(tr)$ when $x \in \gamma_v$}. 
\end{equation}
The two first properties are clear. The third property is requested by the construction of $\sigma$ and is the reason why we couldn't directly set $v_t(x) = tv(x)$ in (\ref{8a5}). Let us check (\ref{eq_vt3}).
For $x \in \gamma_v$, we use $v(x) \leq v_{max}(r)$ (see \eqref{eq_vlessvmax}) and we assume that $C_0$ is bigger than the constant in \eqref{eq_vmax} to estimate
\begin{equation} \label{8b23}
    v_{max}(tr) \geq t v_{max}(r) - C_0 \tau (1 - t) r
    \geq t v(x) - C_0 \tau (1 - t) r.
\end{equation}
Our claim follows.

We are now ready to define $\sigma_1$ (there will be another choice $\sigma_2(x,t)$).
For $x \in \gamma_v$, we choose $\sigma_1(x,t)$ to be the point of $F_v \cap S_{tr(x)}$ that lies exactly at distance $v_t(x)$ from $\xi(tr(x))$. 
This relies on $v_t(x) \leq v_{max}(tr(x))$, see \eqref{eq_vt3}, because if $v_t(x)$ was bigger than 
$v_{max}(t r(x))$, there would be two points of $E_\infty \cap S_{tr(x)}$ at distance $v_t(x)$ from $\xi(tr(x))$ and there is no intrinsic way of choosing between $F_+$ and $F_-$.

Next consider $x \in \gamma_{\pm} \cap B(r\xi_0,10cr)$.
We can proceed similarly, but this time we choose the point of 
$(F_v \cup F_\pm) \cap S_{tr}$ (in the same face $F_\pm$ if we need to choose) 
that lies at distance $v_t(x)$ from $\xi(t r(x))$, and call this point $\sigma_1(x,t)$.
It does not get too far from $tr \xi_0$; we have $\sigma_1(x,t)  \in S_{tr} \cap B(tr\xi_0,20ctr)$ since 
$\abs{\sigma_1(x,t)  - \xi(tr(x))} \leq v_t(x) \leq tv(x) \leq 10ctr$ and $\abs{\xi(tr) - tr\xi_0} \leq C \tau r$.
In general, the part of the face $F_\pm$ near $G$ might be mapped to 
$F_v$ for $t < 1$, but this is all right:
we do not need $\sigma_1(\cdot, t)$ to be injective.
Anyway we essentially have no choice because the face $F_v \cap S_r$ might be reduced to $\set{\xi(r)}$ 
and $F_v \cap S_{t r}$ may be a bit larger but we want that $\sigma_1(\xi(r),t) = \xi(tr)$.

  \begin{figure}[!h]
  \centering
    \includegraphics[width=8cm]{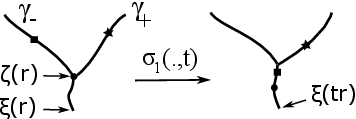}
\vskip-0.2cm
 \caption{Case 2, typical behavior of $\sigma_1 : S_r \cap E \to S_{tr} \cap E$}
\label{fig3}
\end{figure}

It is easy to see that with this choice, $\sigma_1$ has the required properties \eqref{8a1}, \eqref{8a2} when $t=1$ and \eqref{8a3}.
We will also establish that $\sigma_1$ is Lipschitz but we need to show first that for $x \in F_i \cap S_r \cap B(r\xi_0,20cr)$ and $t \in [1/10,10]$,
\begin{equation}\label{eq_txE}
    \mathrm{dist}(tx, F_i \cap S_{tr}) \leq C \tau \abs{1 - t} r.
\end{equation}

    We focus on the case $i \in \set{+,-}$. The reader will see that the proof works with $i = v$; in this case $F_i$ is bounded by two curves but the strategy is the same. 
    It also adapts easily to any other radius in $(r/2,2r)$.
    So we fix a point $x \in F_i \cap S_r \cap B(r\xi_0,20cr)$. We recall that in a suitable choice of coordinate system, $F_i$ can be described as
    \begin{equation}\label{eq_defF}
        \set{z + \varphi_i(z) | z \in P_i \ \text{such that} \ z \cdot \xi_i \geq \psi(z_1) \cdot \xi_i} \cap A_r \cap H(100 c,\xi_0).
    \end{equation}
    where $\varphi : P_i \to P_i^\perp$ is a $\tau$-Lipschitz function such that $\abs{\varphi} \leq C \tau r$, $P_i$ is the vectorial plane generated by $\xi_0 = (1,0,0,0)$ and $\xi_i = (0, \pm \sqrt{3}/2, 1/2,0)$.
    We let $z = (z_1,z_2,z_3,z_4)$ be the orthogonal projection of $x$ onto $P_i$; thus 
    \begin{equation}\label{eq_xz}
        x = z + \varphi_i(z) \quad \text{and} \quad z \cdot \xi_i \geq \psi(z_1) \cdot \xi_i.
    \end{equation}
    Notice that since $\abs{x} = r$ and $\abs{\varphi_i} \leq C \tau r$ with $\tau$ small, we have $r/2 \leq \abs{z} \leq r$.
    For $s \in [1/20,20]$, we recall the definition of
    \begin{equation}
        z(s) = s z_0 \xi_0 + \max(s z \cdot \xi_1, \psi(sz_0) \cdot \xi_1) \xi_1,
    \end{equation}
    which satisfies moreover $\abs{z(s) - sz} \leq C \tau \abs{1 - s} r$ (see Section \ref{S6a})

    Let $x(s)$ denoting $z(s) + \varphi_i(z(s))$ for $s \in [1/20,20]$, it is clear that $\abs{x(s)} \in (r/100,100r)$ and that
    \begin{equation}
        \mathrm{dist}(\abs{x}^{-1}x, \abs{x(s)}^{-1} x(s)) \leq 2 \abs{sx}^{-1} \abs{sx - x(s)} \leq C \tau
    \end{equation}
    so $x(s) \in A_r \cap H(100c,\xi_0)$ (assuming possibly $\tau$ small enough depending on $c$) and in turn that $x(s) \in F_i$ by (\ref{eq_defF}).

    Now, in order to show (\ref{eq_txE}), we want to prove that for all $t \in [1/10,10]$, there exists some $s \in [1/20,20]$ such that $\abs{x(s)} = tr$ and
    \begin{equation}\label{eq_txEpartial}
        \mathrm{dist}(tx, x(s)) \leq C \tau \abs{1 - t} r. 
    \end{equation}
    The map $s \mapsto \abs{x(s)}$ goes continuously from a value $<r/10$ at $s=1/20$ to a value $> 10r$ as $s=20$ so for all $t \in [1/10,10]$, there exists $s \in [1/20,20]$ such that $\abs{x(s)} = tr$. 
    In order to show (\ref{eq_txEpartial}), the main step is to prove that $\abs{t - s} \leq C \tau \abs{1 - t}$. 
    Observe that if two unit vectors $u = (u_1,u_2)$ and $v = (v_1,v_2)$ in $\R^2 \times \R^{n-2}$ are such that $\abs{u_2}, \abs{v_2} \leq 1/2$, then $\abs{\abs{u_1} - \abs{v_1}} \leq C \abs{u_2 - v_2}$.
    Applying this observation to $tx$ and $x(s) = z(s) + \varphi_i(z(s))$ which have the same norm $tr$, we obtain
    \begin{align}
        \big|\abs{z(s)} - \abs{tz}\big| &\leq C \abs{\varphi_i(z(s)) - t \varphi_i(z)}\\
        &\leq C \abs{\varphi_i(z) - t \varphi_i(z)} + C \abs{\varphi_i(z(s)) - \varphi_i(z)}\\
        &\leq C \tau \abs{1 - t} r + C \tau \abs{z(s) - z}.
    \end{align}
    As $\abs{z(s) - sz} \leq C \tau \abs{1 - s}r$, we deduce
    \begin{equation}
        \big|\abs{z(s)} - \abs{tz}\big| \leq C \tau \abs{1 - t} r + C \tau \abs{1 - s} r.
    \end{equation}
    Using the fact that $u \mapsto \big| \abs{tz} - \abs{u} \big|$ is $1$-Lipschitz, 
    we can also bound from below
    \begin{align}
        \big|\abs{z(s)} - \abs{tz}\big| &\geq \big|\abs{sz} - \abs{tz}\big| - \abs{z(s) - sz}\\
                                        &\geq C^{-1} \abs{t - s} r - C \tau \abs{1 - s}r 
    \end{align}
    and it follows that $\abs{t - s} \leq C \tau \abs{1 - t} + C \tau \abs{1 - s}$.
    Noticing that $\abs{1 - s} \leq \abs{t - s} + \abs{1 - t}$, this implies that $\abs{1 - s} \leq C \abs{1 - t}$ and next, $\abs{t - s} \leq C \tau \abs{1 - t}$.
    Together with \eqref{eq_xz}, 
    the fact that $\varphi$ is $\tau$-Lipschitz and $\abs{\varphi} \leq C \tau r$, we conclude that
    \begin{align}
        \mathrm{dist}(tx, z(s) + \varphi_i(z(s))) &\leq \abs{tz - z(s)} + \abs{t\varphi_i(z) - \varphi_i(z(s))}\\
        &\leq \abs{tz - sz} + \abs{sz - z(s)} + \abs{t \varphi_i(z) - \varphi_i(z)} + \abs{\varphi_i(z) - \varphi_i(z(s))}\\
        &\leq C \tau \abs{1 - t} r;
    \end{align}
\eqref{eq_txEpartial} and \eqref{eq_txE} follow. 

\ms 
Our next long-time goal is to 
prove that $\sigma_1$ is Lipschitz in the sense that for all $x,y \in E_{\infty} \cap A_r(2) \cap H(10c,\xi_0)$ 
and $1/2 \leq t,s \leq 1$,
\begin{equation}\label{8a09}
    \abs{\sigma_1(x,t) - \sigma_1(y,s)}   \leq C \abs{x - y} + C \abs{t - s} r.
\end{equation}
and even has the slight ``contraction'' property that
\begin{equation} \label{8a9}
    |\sigma_1(x,1/2)-\sigma_1(y,1/2)| \leq \frac{3}{4} \abs{x - y}.
\end{equation}
This will be useful when we compose mappings. 

Let us first check that 
\begin{equation}\label{eq_lip_sigmat2}
    \abs{\sigma_1(x,t)  - tx} \leq C \tau (1 - t) r(x).      
\end{equation}
Suppose for instance that $x \in F_+ \cup F_v$ (the case when $x \in F_- \cup F_v$ can be done the same way).
Then by \eqref{eq_txE} we can find $z \in (F_+ \cup F_v) \cap S_{tr(x)}$ such that 
\begin{equation} \label{8c42}
  \abs{z -tx} \leq C \tau (1 - t) r(x).
\end{equation}
We want a similar estimate for $ \abs{z -\sigma_1(x)}$, and thanks to the biLipschitz %
behavior of $z \mapsto \dist(z, \xi(t r(x))$ described below \eqref{8b18}, it is enough to show that
\begin{equation} \label{8c43}
    \abs{\dist(z, \xi(t r(x)) - \dist(\sigma_1(x,t)), \xi(t r(x))} \leq C \tau (1 - t) r(x).
\end{equation}
By definition of $\sigma_1(x,t)$, $\dist(\sigma_1(x,t)), \xi(t r(x)) = v_t(x)$, and \eqref{eq_vt2}
says that $\abs{v_t(x) - t v(x)} \leq C \tau (1 - t) r(x)$. Then by \eqref{8c42} it is enough to check that
\begin{equation} \label{8c44}
    \abs{\dist(tx, \xi(t r(x))) - t v(x)} \leq C \tau (1 - t) r(x).
\end{equation}
Furthermore, it follows from \eqref{eq_xirho} and the fact that 
$\abs{\xi(t r(x))} = t r(x) = t \abs{\xi(r(x))}$ that $\abs{\xi(t r(x))-t \xi(r(x))} \leq C \tau (1 - t) r(x)$,
and now \eqref{8c44} and \eqref{eq_lip_sigmat2} follow, because \eqref{8a4} says that
$\dist(x, \xi(r(x)) = v(x)$.

Our next step is to show that
\begin{equation}\label{eq_lip_sigmat}
    \abs{\frac{\sigma_1(x,t)}{tr(x)} - \frac{\sigma_1(y,s)}{sr(y)}} \leq \abs{\frac{x}{r(x)} - \frac{y}{r(y)}} + C \tau_0 r^{-1} \abs{x - y} + C \tau \abs{t - s}.
\end{equation}
We will see later how this implies (\ref{8a09}), (\ref{8a9}).

We start with the case where both $\sigma_1(x,t)$ and $\sigma_1(y,s)$ 
lie $F_+ \cup F_v$ (the case $F_- \cup F_v$ can be done in the same way).
According to (\ref{eq_txE}), applied to $\sigma_1(y,s) \in S_{s r(y)}$ and the factor $\wt t = \frac{t r(x)}{s r(y)}$, 
there exists $z \in (F_+ \cup F_v) \cap S_{tr(x)}$ such that
\begin{equation}
    \abs{\frac{tr(x)}{sr(y)} \sigma_1(y,s) - z} 
    \leq C \tau \abs{tr(x) - sr(y)} \leq C \tau \abs{x - y} + C \tau \abs{t - s} r,
\end{equation}
and thus
\begin{equation}\label{eq_sigmaz}
    \abs{\frac{\sigma_1(y,s)} 
    {sr(y)} - \frac{z}{tr(x)}} \leq C \tau r^{-1} \abs{x - y} + C \tau \abs{t - s}.
\end{equation}
In order to prove (\ref{eq_lip_sigmat}), it then suffices to show 
\begin{equation}\label{eq_sigmaz2}
    \abs{\frac{\sigma_1(x,t)}{tr(x)} 
    - \frac{z}{tr(x)}} \leq \abs{\frac{x}{r(x)} - \frac{y}{r(y)}} + C \tau_0 r^{-1} \abs{x - y} + C \tau \abs{t - s}.
\end{equation}
Apply the discussion below \eqref{8b18} to the radius $t r(x)$. Recall that 
for all $0 \leq \rho \leq 20 c t r(x)$, there exists a unique point 
$x(\rho) \in (F_+ \cup F_v) \cap S_{t r(x)} \cap H(\xi_0,20c)$ such that $v(x(\rho)) = \rho$ 
and the corresponding mapping $\rho \mapsto x(\rho)$ is $(1 + C \tau_0)$-Lipschitz
(here everything, including the mapping 
$v(\cdot) = \dist(\cdot, \xi(t r(x))$ from \eqref{8a4}, 
depends also on $t r(x)$, but we fix it). Now recall that $v(\sigma_1(x,t)) = v_t(x)$
(see below \eqref{8b23}), while by definition of $v$, $v(z) = \abs{z - \xi(tr(x))}$.
So
\begin{equation}
    \abs{\sigma_1(x,t) - z} 
    \leq (1 + C \tau_0) \big| v_t(x) - \abs{z - \xi(tr(x))} \big|.
\end{equation}
Next, by the triangular inequality and since $v_s(y) = \abs{\sigma(y,s) - \xi(sr(y))}$,
\begin{align}
    \abs{\frac{v_t(x)}{tr(x)} - \frac{\abs{z - \xi(tr(x))}}{tr(x)}} &\leq \abs{\frac{v_t(x)}{tr(x)} - \frac{v_s(y)}{sr(y)}} + \abs{\frac{v_s(y)}{sr(y)} - \frac{\abs{z - \xi(tr(x))}}{tr(x)}}\\
                                                           &\leq \abs{\frac{v_t(x)}{tr(x)} - \frac{v_s(y)}{sr(y)}} + \abs{\frac{\sigma(y,s)}{sr(y)} - \frac{z}{tr(x)}} + \abs{\frac{\xi(sr(y))}{sr(y)} - \frac{\xi(tr(x))}{tr(x)}}.
\end{align}
For the first term at the right-hand side, we recall that $\rho \mapsto \xi(\rho)$ is $C \tau r^{-1}$-Lipschitz, and we use the definition of $v_t(x)$ and $v_s(y)$ to estimate
\begin{align}
    \abs{\frac{v_t(x)}{tr(x)} - \frac{v_s(y)}{sr(y)}}  &\leq \abs{\frac{v(x)}{r(x)} - \frac{v(y)}{r(y)}} 
    + C_0  \tau \abs{\frac{1-t}{t} - \frac{1-s}{s}}\\
            &\leq \abs{\frac{x}{r(x)} - \frac{y}{r(y)}} + \abs{\frac{\xi(r(x))}{r(x)} - \frac{\xi(r(y))}{r(y)}} 
            + C_0 \tau \abs{\frac{1}{t} - \frac{1}{s}}\\ 
                         &\leq \abs{\frac{x}{r(x)} - \frac{y}{r(y)}} + C \tau r^{-1} \abs{x - y} + C \tau \abs{t - s}.
\end{align}
We control the second term with (\ref{eq_sigmaz}) and the last term with the Lipschitz property of 
$\rho \mapsto \xi(\rho)/\rho$.
One can easily deduce (\ref{eq_sigmaz2}), assuming $\tau \leq \tau_0$ if necessary, and in turn (\ref{eq_lip_sigmat}).

\ms
Our next case is when $\sigma_1(x,t) \in F_+ \setminus \Gamma$ and $\sigma_2(s,y) \in F_- \setminus \Gamma$ 
(in particular, we necessarily have $x \in F_+$ and $y \in F_-$ in this case). Here, the main point is to show that
\begin{equation}\label{eq_xyr}
    \abs{x - y} \geq C^{-1} (1 - t) r.
\end{equation}
(We could also prove similarly that $\abs{x - y} \geq C^{-1} (1 - s) r $).
The assumption on $\sigma_1(x,t)$ means that $v_t(x) > v_{max}(tr(x))$, i.e., 
$t v(x) - C_0 \tau (1 - t) r(x) \geq v_{max}(tr(x))$, thus
\begin{equation}
    v(x) - \frac{v_{max}(tr(x))}{t} \geq \frac{C_0}{2} \tau (1 - t) r(x).
\end{equation}
Assuming $C_0$ big enough depending on the constant in (\ref{eq_vmax}), it follows that $v(x) - v_{max}(r(x)) \geq (C_0/4) \tau (1 - t) r(x)$ and thus, by the triangular inequality, 
\begin{equation}
    \abs{x - \zeta(r(x))} \geq \frac{C_0}{4} \tau (1 - t)r(x).
\end{equation}
This is the last time we put an assumption on $C_0$, which will now be considered fixed (depending 
only on $n$).
According to (\ref{eq_txE}), there exists $z \in F_- \cap S_{r(x)}$ such that 
\begin{equation}
    \abs{\frac{y}{r(y)} - \frac{z}{r(x)}} \leq C \tau r^{-1} \abs{r(x) - r(y)} \leq C \tau r^{-1} \abs{x - y}.
\end{equation}
From the lower bound on $\abs{x - \zeta(r(x))}$ and the fact that $x \in F_+ \cap S_{r(x)}$ and $z \in F_- \cap S_{r(x)}$, (\ref{eq_distHi}) tell us that
\begin{equation}
    \abs{x - z} \geq \frac{1}{2} \abs{x -\zeta(r(x))} \geq \frac{C_0}{8} (1 - t) r(x).
\end{equation}
We deduce
\begin{equation}
    \abs{\frac{x}{r(x)} - \frac{y}{r(y)}} \geq \frac{C_0}{8} (1 - t) r(x) - C \tau r^{-1} \abs{x - y}.
\end{equation}
but since $\abs{x/r(x) - y/r(y)} \leq 2 r(x)^{-1} \abs{x - y}$, (\ref{eq_xyr}) follows.
Then, (\ref{eq_lip_sigmat2}) and (\ref{eq_xyr}) allow to estimate
\begin{align}
    \abs{\frac{\sigma_1(t,x)}{tr(x)} - \frac{\sigma_1(s,y)}{sr(y)}} &\leq \abs{\frac{\sigma_1(t,x) - tx}{tr(x)}} + \abs{\frac{x}{r(x)} - \frac{y}{r(y)}} + \abs{\frac{\sigma_1(s,y) - sy}{sr(y)}}\\
                                                                    &\leq \abs{\frac{x}{r(x)} - \frac{y}{r(y)}} + C \tau (1 - t) + C \tau (1 - s)\\
                                                                    &\leq \abs{\frac{x}{r(x)} - \frac{y}{r(y)}} + C \tau r^{-1} \abs{x - y}.
\end{align}
This ends the proof of (\ref{eq_lip_sigmat}).

\ms
We still need to deduce (\ref{8a09}) and (\ref{8a9}) from (\ref{eq_lip_sigmat}). 
Observe that since $\abs{\sigma_1(x,t)} = t r(x)$ (and similarly with 
$y$ and $s$), a simple computation (expand the right-hand side, four terms out of six cancel, and we get the
same result from the right-hand side) yields 
\begin{equation}
    \abs{\sigma_1(x,t) - \sigma_1(y,s)}^2 
    - \abs{tx - sy}^2 = st r(x) r(y) \left(\abs{\frac{\sigma_1(x,t)}{tr(x)} - \frac{\sigma_1(y,s)}{sr(y)}}^2 - \abs{\frac{x}{r(x)} - \frac{y}{r(y)}}^2\right).
\end{equation}
Thanks to the usual estimate (\ref{eq_theta_lip}), we also have independently
\begin{equation}
    \abs{\frac{\sigma_1(x,t)}{tr(x)} - \frac{\sigma_1(y,s)}{sr(y)}} + \abs{\frac{x}{r(x)} - \frac{y}{r(y)}} 
    \leq C r^{-1} \left(\abs{\sigma_1(t,x) - \sigma_1(s,y)} + \abs{tx - ty}\right).
\end{equation}
From the two above lines and (\ref{eq_lip_sigmat}), we deduce
\begin{align}
    \abs{\sigma_1(t,x) - \sigma_1(s,y)} - \abs{tx - sy} &\leq C r \left(\abs{\frac{\sigma_1(x,t)}{tr(x)} - \frac{\sigma_1(y,s)}{sr(y)}} - \abs{\frac{x}{r(x)} - \frac{y}{r(y)}}\right)\\
                                                        &\leq C \tau_0 \abs{x - y} + C \tau r \abs{t - s}
\end{align}
and (\ref{8a09}), (\ref{8a9}) follow.

\ms 
As the reader may have guessed, there is a second choice $\sigma_2(x,t)$, that we will rather use when 
$x \in S_r \cap B(r\xi_0,10cr) \setminus B(\zeta(r),2 cr)$ say.
In this case, $E_\infty$ is very flat near $t x$, and we want to use the intrinsic projection {$p_{tx, t r}$ of Section \ref{S6a} associated to $E_\infty \cap S_{tr}$.

Precisely, this projection is well-defined thanks to the following properties: for all $x \in F_{\pm} \cap S_r \cap B(r\xi_0,10cr) \setminus B(\zeta(r), cr)$, for all $1/2 \leq t \leq 1$ and for all $z \in B(tx,10^{-1}ctr)$, we have
\begin{equation}\label{eq_HF1}
    \mathrm{dist}(z, F_{\pm} \cap S_{tr}) \leq 10^{-3} c t r
\end{equation}
and
\begin{equation}\label{eq_HF2}
    E_\infty \cap S_{tr} \cap B(z, 10^{-1} c t r) = F_{\pm} \cap S_{tr} \cap B(z, 10^{-1} c t r) \ \text{with} \ \zeta(r) \notin B(z,10^{-1} ctr).
\end{equation}
    The first point directly comes from (\ref{eq_txE}), as $\tau$ can be chosen small depending on $c$ so that $\mathrm{dist}(tx,F_{\pm} \cap S_{tr}) \leq 10^{-4} ctr$. For the second point, recall that the map $\rho \mapsto \zeta(\rho)/\rho$ is $C \tau r^{-1}$-Lipschitz so, assuming again that $\tau$ small enough depending on $c$, we have $\abs{\zeta(tr) - t \zeta(r)} \leq C \tau (1 - t)r \leq ctr/10$. As $\abs{x - \zeta(r)} \geq cr$, this implies $\abs{tx - \zeta(tr)} \geq 9ctr/10$. Moreover, (\ref{eq_txE}) also shows that there exists a point $y \in F_{\pm} \cap S_{tr}$ such that $\abs{tx - y} \leq ctr/10$ and in particular $\abs{y - \zeta(tr)} \geq 8ctr/10$. Then, we know by (\ref{eq_distHi}) that $y$ is at distance $\geq 4ctr/10$ from the other $F_i \cap S_{tr}$. We conclude that $tx$ is at distance $\geq 3ctr/10$ from the others $F_i \cap S_{tr}$, and thus any $z \in B(tr,10^{-1}ctr)$ is at distance $\geq 2ctr/10$ from them as well. This proves (\ref{eq_HF2}). In particular, this justifies that $E_{\infty} \cap S_{tr}$ coincide with a $\tau_0$-Lipschitz graph in $B(z,10^{-1}ctr)$.

We can thus set for $x \in E_\infty \cap S_r \cap B(r\xi_0,10cr) \setminus B(\zeta(r),2cr)$,
\begin{equation} \label{8a10}
    \sigma_2(x,t) = p_{tx, tr}(t x).
\end{equation}
This one has the advantage of being the same as soon as we are away from the singularities of $E_\infty$. Thanks to (\ref{eq_pLipschitz2}), we still have
\begin{align}
\abs{\frac{\sigma_2(x,t)}{tr(x)} - \frac{\sigma_2(y,s)}{tr(y)}} 
          &\leq \abs{\frac{x}{r(x)} - \frac{y}{r(y)}} + C \tau_0 r^{-1} \abs{tx - sy}\\
                     &\leq \abs{\frac{x}{r(x)} - \frac{y}{r(y)}} + C \tau_0 r^{-1} \abs{x - y} + C \tau_0 \abs{t - s}.
\end{align}
and thus, as we have seen with $\sigma_1$, the map $\sigma_2$ is Lipschitz (\ref{8a09}) with the contraction property (\ref{8a9}).
We also have
\begin{equation}
    \abs{\sigma_2(t,x) - tx} \leq C \tau (1 - t)r
\end{equation}
    using (\ref{eq_pLipschitz3}) and observing that $\abs{\sigma_2(t,x) - tx} = \abs{p_{tx,tr}(tr) - p_{x,r}(x)}$.

We keep $\sigma(x,t) = \sigma_1(x,t)$
when $x \in E_\infty \cap S_r \cap B(\zeta(r),cr)$, and of course we need to interpolate nicely in the remaining region where $x \in E_\infty \cap S_r \cap B(\zeta(r),2cr) \setminus B(\zeta(r),cr)$. We set
\begin{equation}
    a(x) = \frac{\abs{x - \zeta(r(x))} - c r(x)}{c r(x)} \in [0,1]
\end{equation}
and choose 
\begin{equation}\label{8a11}
    \sigma(x,t) = p_{tx, tr}\left(a(x) \sigma_2(x,t) + (1-a(x))\sigma_1(x,t)\right).
\end{equation}
Now, the formula will be easier to glue because of $\sigma_2$.
Here we added the extra projection $p_{tx, t r}$ because the intermediate point may lie slightly away from $E_\infty$ again. The point $z := a(x) \sigma_2(x,t) + (1-a(x))\sigma_1(x,t)$ satisfies $\abs{z - tx} \leq C \tau (1 - t)r \leq 10^{-1} ctr$ so $\sigma(x,t)$ is well-defined as usual thanks to (\ref{eq_HF1}) and (\ref{eq_HF2}).

Let us check that $\sigma$ is Lipschitz (\ref{8a09}) with the contraction property (\ref{8a9}). Letting $z = a(x) \sigma_2(x,t) + (1 - a(x))\sigma_1(x,t)$ and $w = a(y) \sigma_2(y,s) + (1 - a(y)) \sigma_1(y,s)$, we know by (\ref{eq_pLipschitz2}) that
\begin{equation}
    \abs{\frac{\sigma(x,t)}{tr(x)} - \frac{\sigma(y,s)}{sr(y)}} \leq (1 + C \tau) \abs{\frac{z}{tr(x)} - \frac{w}{sr(y)}} + C \tau_0 \abs{x - y} + C \tau_0 \abs{t - s}.
\end{equation}
Then, we show that
    \begin{equation}\label{eq_zwst}
    \abs{\frac{z}{tr(x)} - \frac{w}{sr(y)}} \leq \abs{\frac{x}{r(x)} - \frac{y}{r(y)}} + C \tau_0 \abs{x - y} + C \tau_0 \abs{t - s}.
\end{equation}
We first use the triangular inequality,
\begin{multline}
    \abs{\frac{z}{tr(x)} - \frac{w}{sr(y)}} \leq a(x) \abs{\frac{\sigma_2(x,t)}{tr(x)} - \frac{\sigma_2(y,s)}{sr(y)}} + (1 - a(x)) \abs{\frac{\sigma_1(x,t)}{tr(x)} - \frac{\sigma_1(y,s)}{sr(y)}} \\+ \abs{a(x) - a(y)} \abs{\frac{\sigma_1(y,s) - \sigma_2(y,s)}{sr(y)}}.
\end{multline}
We already know that for $i=1,2$,
\begin{equation}
    \abs{\frac{\sigma_i(x,t)}{tr(x)} - \frac{\sigma_i(y,s)}{sr(y)}} \leq \abs{\frac{x}{r(x)} - \frac{y}{r(y)}} + C \tau_0 \abs{x - y} + C \tau_0 \abs{t - s}.
\end{equation}
As $\rho \mapsto \xi(\rho)/\rho$ is $C \tau r^{-1}$-Lipschitz, we can also estimate
\begin{equation}
    \abs{a(x) - a(y)} \leq \frac{1}{c} \abs{\frac{x}{r(x)} - \frac{y}{r(y)}} + \frac{1}{c} \abs{\frac{\xi(r(x))}{r(x)} - \frac{\xi(r(y)}{r(y)}} \leq \frac{C \abs{x - y}}{cr}.
\end{equation}
And we recall that $\abs{\sigma_i(y,s) - sy} \leq C \tau (1 - s) r \leq C \tau r$, which allows to control
\begin{equation}
    \abs{\frac{\sigma_1(y,s) - \sigma_2(y,s)}{sr(y)}} \leq C \tau.
\end{equation}
Since $\tau$ is allowed to be small depending on $c$ and $\tau_0$, (\ref{eq_zwst}) follow. We deduce
\begin{equation}
    \abs{\frac{\sigma(x,t)}{tr(x)} - \frac{\sigma(y,s)}{sr(y)}} \leq \abs{\frac{x}{r(x)} - \frac{y}{r(y)}} + C \tau_0 \abs{x - y} + C \tau_0 \abs{t - s},
\end{equation}
and consequently (\ref{8a09}), (\ref{8a9}) as usual.

We also have
\begin{equation}
    \abs{\sigma(x,t) - tx} \leq C \tau (1 - t) r.
\end{equation}
    This uses the fact that $z \mapsto p_{tx,tr}(z)$ is $C$-Lipschitz, that the point $z = (a(x) \sigma_2(x,t) + (1-a(x))\sigma_1(x,t))$ satisfies $\abs{z - tx} \leq C \tau (1 - t)r$ and that $\abs{p_{tx,tr}(tx) - tx} \leq C \tau (1 - t)r$.

\ms 
This ends the definition of $\sigma$ for Case 2 and Case 2 bis.
Our next case is near a generic cone of type $\bV$, and we don't need to change 
anything there, except that when the angle $\alpha$ gets closer to $\pi$, 
we rapidly have no vertical face $F_v$ left, and the definition of $v_t$ reduces to $v_t(x) = t v(x)$.
When $\alpha$ becomes close enough to $\pi$, we can brutally 
use $\sigma_2(x,r)$ in the whole region, and as we did near \eqref{7a13}, interpolate nicely 
between the two formulas as follows. Let $\tau_1 \in (0,1)$ be such that previous constructions work for all $\tau \leq \tau_1$ (it only depends on $n$, $\nu$ and $c$).
When $\alpha$ varies between $\pi - \tau_1$ and $\pi - 2 \tau_1$, the map $\sigma$ built for Case 3 is replaced by
\begin{equation}
    p_{tx,tr}((2 - \alpha'/\tau_1) tx + (\alpha'/\tau_1 - 1) \sigma(x,t)).
\end{equation}

Our Case 4, with approximation by half planes $H \in \bH$, is similar to Case 3, 
except that now we have only one face (and this case always stays far from the other ones).

Near a point of type $\bY$ (hence away from $\Gamma$),
we proceed as in the previous two cases with the simple formula $v_t(x) = t v(x)$.
As before, we also need to interpolate smoothly our two formulas, the simple one
$v_t(x) = t v(x)$ for the present case of $\bY$, and the one above for truncated
$\bY$ cones. But the formula can stay the same (as long as the truncated cone stays 
at small distance from a sharp $\bV$; we just use the slightly different projections
discussed above).

In the remaining case when $E_\infty$ is flat near $x$, we use directly 
$\sigma_2(x,t)$ from \eqref{8a10}.

\ms
At this point we have built a function
\begin{equation}
    \sigma : \left[E_\infty \cap B(0,r_0)\right] \times [1/2,1] \to E_\infty \cap B(0,r_0)
\end{equation}
that solves the Lemma but only for $1/2 \leq t \leq 1$.
The map $\sigma$ is Lipschitz in the sense that for all $x, y \in E_\infty \cap B(0,r_0)$ and for all $1/2 \leq t,s \leq 1$, we have
\begin{equation}\label{eq_sigma_lip}
    \abs{\sigma(x,t) - \sigma(y,s)} \leq C \abs{x - y} + C \min(\abs{x},\abs{y}) \abs{t - s}.
\end{equation}
The reasoning is very similar to what we did with $\pi$, at the end of Section \ref{S7}. Let $x, y \in E_{\infty} \cap B(0,r_0)$. If $\abs{x - y} \geq 10^{-4} \nu c \max(\abs{x},\abs{y})$, then
\begin{align}
    \abs{\sigma(x,t) - \sigma(y,s)}
    &\leq \abs{\sigma(x,t) - tx} + \abs{tx - sy} + \abs{\sigma(y,s) - sy}
    \nonumber\\
    &\leq C \tau \abs{x} + \abs{tx - sy} + C \tau \abs{y}
    \leq \abs{tx - sy} + 20C \tau \abs{x - y}.
\end{align}
If on the other hand $\abs{x - y} \leq 10^{-4} \nu c \abs{x}$, one can see as in Section \ref{S7} that $x$, $y$ belong to a chart $A_r(2) \cap H(10c_i,\xi_i)$ and uses the fact that $\sigma$ is Lipschitz in such a box.
We can assume $\tau$ small depending on $\nu$ and $c$ so that in all cases, $\abs{\sigma(x,t) - \sigma(y,s)} \leq \abs{tx - sy} + \abs{x - y}/10$. This concludes the proof of (\ref{eq_sigma_lip}) and shows moreover that $\sigma$ has the contracting property: for $x,y \in E_\infty \cap B(0,r_0)$,
\begin{equation}
    \abs{\sigma(x,1/2) - \sigma(y,1/2)} \leq \frac{3}{4}\abs{x - y}. 
\end{equation}

We come to the last part of this proof which consists in extending $\sigma(x,t)$ to $t \in (0,1]$, and we shall simply compose.
First define $\sigma_k(x)$, $x\in B(0,r_0)$ and $k \geq 0$,  by induction on $k$: set 
$\sigma_0(x) = x$, and then $\sigma_{k+1}(x) = \sigma(\sigma_k(x),1/2)$ for $k \geq 0$.
Now write any $t \in (0,1]$ as $t = s 2^{-k}$, with $k \geq 0$ and $1/2 < s \leq 1$, and set
\begin{equation} \label{8a12}
    \sigma(x,t) = \sigma(\sigma_k(x),s).
\end{equation}
There is no jump across the integers since
\begin{equation}
    \lim_{s \to 1^-} \sigma(x,s 2^{-k}) = \sigma(\sigma_k(x),1) = \sigma_k(x)
\end{equation}
and
\begin{equation}
    \lim_{s \to \tfrac{1}{2}^+} \sigma(x,s 2^{-(k-1)}) = \sigma(\sigma_{k-1}(x),1/2) = \sigma_k(x).
\end{equation}
It is also clear that $\sigma$ preserves the spheres. As a consequence, $\sigma$ naturally extends to $t = 0$ by setting $\sigma(x,0) = 0$.

Finally, we justify that $\sigma$ has the Lipschitz property (\ref{8a0}). This amounts to showing that for all $x,y \in E_\infty \cap B(0,r_0)$ and $0 < t,s \leq 1$,
\begin{equation}\label{eq_sigma_lip1}
    \abs{\sigma(x,t) - \sigma(x,s)} \leq \abs{x} \abs{t - s}
\end{equation}
and
\begin{equation}\label{eq_sigma_lip2}
    \abs{\sigma(x,t) - \sigma(y,t)} \leq \abs{x - y}. 
\end{equation}

We start with (\ref{eq_sigma_lip1}).
If $2^{-k-1} < t_1, t_2 \leq 2^{-k}$ for some $k \geq 0$, we use the fact that $\abs{\sigma_k(x)} = 2^{-k} \abs{x}$ to see that (writing $t_i = s_i 2^{-k}$ where $s_i \in [1/2,1]$),
\begin{align}
    \abs{\sigma(x,t_1) - \sigma(x,t_2)} 
    &= \abs{\sigma(\sigma_k(x),s_1) - \sigma(\sigma_k(x),s_2)}
    \nonumber\\
    &\leq C \abs{\sigma_k(x)} \abs{s_1 - s_2}
    \leq C \abs{x} \abs{t_1 - t_2}.
\end{align}
We deduce the general case $0 < t_1, t_2 \leq 1$ by summing the inequalities in each interval $[2^{-k-1}, 2^{-k}]$ between $t_1$ and $t_2$.

We pass to (\ref{eq_sigma_lip2}). If $2^{-k-1} < t \leq 2^{-k}$ for some $k \geq 0$ (writing $t = s 2^{-k}$ as before), we have
\begin{align}
    \abs{\sigma(x,t) - \sigma(y,t)} = \abs{\sigma(\sigma_k(x),s) - \sigma(\sigma_k(y),s)}.
    \leq C \abs{\sigma_k(x) - \sigma_k(y)} 
\end{align}
The contraction property of $\sigma$ shows that for all $k \geq 1$, $\abs{\sigma_k(x) - \sigma_k(y)} \leq (3/4) \abs{\sigma_{k-1}(x) - \sigma_{k-1}(y)}$ and we deduce by iteration that $\abs{\sigma_k(x) - \sigma_k(y)} \leq \abs{x - y}$. This ends the proof.

\qed

\begin{rem}\label{rem_l8a1}
    In the heat of the construction, it would seem that 
    we forgot that we wanted the sliding condition \eqref{8a3},
    and now we notice that it may fail. This does not happen in Cases 3 
    (generic $\bV$-cones) and  Case 6 ($\bY$-points), or even when we are far from $L_0$, so it only happens in 
    Case 1 (or Case 3 with a very flat cone), or Case 5, which anyway we assimilated to Case 1
    at the beginning of the argument. That is, we are worried when $E_\infty$ is very 
    close to a plane $P$ that contains $L_0$, the point $\zeta(r)$ of $S_r \cap \Gamma$ lies
    in $E_\infty$, and we cannot enforce \eqref{8a3} because anyway $E_\infty$ leaves
    $\Gamma$ gently away from this point. We will have to remember this case and counter it
    with a trick.
\end{rem}

\begin{rem}\label{l8a2} 
    We said that we cannot make $\sigma$ injective because of Case 2 
    (where the various $E_\infty \cap S_r$ do not have the same topology),
    but if Case 2 does not arise, we can use the construction above to obtain a
    biLipschitz parameterization of $E_\infty \cap B(0,r)$ by the cone over $E_\infty \cap S_r$.
    We present this in Section \ref{Sbil} because this was not obvious a priori, but we shall not need this fact for our main theorems.
\end{rem}

\section{A retraction onto \texorpdfstring{$E_\infty$}{Einfty} defined near \texorpdfstring{$E_\infty \cap B(0,r_0)$}{Einfty inter B(0,r0)}} 
\label{S9}

In this section we use the mappings of Sections \ref{S7} and \ref{S8a} 
to construct a Lipschitz retraction on $E_\infty \cap B(0, r_0)$, 
assuming  as usual that $0 \in E_\infty$.

\begin{lem}\label{l9a1}
    Let $E_\infty$, $\varepsilon_0$ and the ball $B_0 = B(0, r_0)$ be as in Section \ref{S5}.
    Then (if $\varepsilon_0$ is small enough in (\ref{5a1})) there is a deformation 
    $\rho : \R^n \times [0,1] \to \R^n$ with the following properties:
    \begin{equation} \label{9a1}
        \rho(x,0) = x \text{ for } x\in \R^n;
    \end{equation}
    \begin{equation} \label{9a2}
        \rho(x,t) = x \text{ for $x\in \R^n \sm B_0$ and for $x\in E_\infty$;}
    \end{equation}
    \begin{equation} \label{9a3}
        \rho(x,t) \text{ is $C$-Lipschitz on } \R^n \times [0,1].
    \end{equation}
    \begin{equation} \label{9a4}
        \rho(x,1) \in E_\infty  \text{ for } x\in B(0,r_0/2).
    \end{equation}
    The constant $C \geq 1$ depends only on $n$.
\end{lem}

\ms
As usual, more properties will appear later on. 
Let $\pi$ and $c_1$ be as in Proposition \ref{l6a1}, and choose a cut-off function 
$\eta : [0,+\infty) \to [0,1]$ such that $\eta(t) = 1$ when $0 \leq t \leq c_1/3$, 
$\eta(t) = 0$ when $t \geq 2c_1/3$, and $\eta$ linear in between, i.e., 
$\eta(t) = (2c_1-3t)/c_1$ when $c_1/3 \leq t \leq 2c_1/3$.
Then define a scaling factor $h : B_0 \setminus \set{0} \to [0,1]$, given by
\begin{equation} \label{9a5}
    h(x) = \eta(\dist(x,E_\infty)/r(x)),\quad x \in B_0 \setminus \set{0}.
\end{equation}
We observe that $h(x) \in [0,1]$ and that, since $\mathrm{dist}(x,E_{\infty}) \leq r(x)$,
\begin{equation}\label{9a55}
    \abs{h(x) - h(y)} \leq \frac{C}{\min(r(x),r(y))} \, \abs{x - y}.
\end{equation}
Then we define a local Lipschitz retraction $g : B_0 \to E_\infty$ by
\begin{equation} \label{9a6}
    g(x) =  \sigma(\pi(x), h(x)) \ \text{when} \ \dist(x,E_\infty) < c_1 r(x),
\end{equation}
and simply
\begin{equation} \label{9a7}
    g(x) = 0 \ \text{when} \ \dist(x,E_\infty) \geq \frac{2 c_1}{3} \,  r(x). 
\end{equation}
Notice that in the overlapping region where $2c_1/3 r(x) \leq \dist(x,E_\infty) < c_1 r(x)$, the definitions coincide because $h(x) = 0$ and hence $\sigma(\pi(x), h(x)) = 0$. 
In all cases, we have $g(x) \in E_\infty$ by definition of $\pi$, $\sigma$ and because $0 \in E_\infty$.
According to the property of $\sigma$ and $\pi$, we have for all $x \in B_0$,
\begin{equation}
    \abs{g(x)} = h(x) r(x) \leq r(x).
\end{equation}
Notice also that when $x \in E_\infty \cap B_0$ (in particular $h(x) = 1$ if $x \ne 0$), we have
\begin{equation} \label{9a9}
    g(x) = \pi(x) = x.
\end{equation}
Next, we check that $g$ is $C$-Lipschitz in $B_0$ by distinguishing three cases.
Let $x, y \in B_0$. Our first case is when $\dist(x,E_\infty) < c_1 r(x)$ and $\dist(y,E_\infty) < c_1 r(y)$. We use the properties of $\sigma$, $\pi$ and $h$ (see in particular (\ref{9a55})), to see that
\begin{align}
    \abs{g(x) - g(y)}  
    \leq C \abs{\pi(x) - \pi(y)} + C \min(\abs{x},\abs{y}) \abs{h(x) - h(y)} 
    \leq C \abs{x - y}.
\end{align}
Next assume that $\dist(x,E_\infty) < c_1 r(x)$ but $\dist(y,E_\infty) \geq c_1 r(y)$.
Actually, we can directly assume $\dist(x,E_\infty) < 2c_1 r(x)/3$ in this case, otherwise $g(x) = g(y) = 0$ 
and there would be nothing to do.
As we have seen before, $\abs{g(x)} \leq r(x)$ and $g(y) = 0$ so $\abs{g(x) - g(y)} \leq r(x)$, and then we show that $r(x) \leq C \abs{x - y}$.
We consider $z \in E_\infty$ such that $\abs{x - z} < 2 c_1/3 r(x)$.
We must have $\abs{y - z} \geq \mathrm{dist}(y,E_{\infty}) \geq c_1 r(y)$ whence
\begin{equation}
    \abs{x - y} \geq  \abs{y - z} - \abs{x - z} \geq c_1 r(y) - 2 c_1 r(x)/3
\end{equation}
but according to the triangular inequality, we have $r(y) \geq r(x) - \abs{x - y}$ so
\begin{equation}
    \abs{x - y} \geq \frac{c_1}{3(1 + c_1)} r(x).
\end{equation}
Then $ \abs{g(x) - g(y)} \leq r(x) \leq C \abs{x - y}$, as needed. The proof is the same when 
$\dist(x,E_\infty) \geq c_1 r(x)$ but $\dist(y,E_\infty) < c_1 r(y)$. Finally, when
$\dist(x,E_\infty) \geq c_1 r(x)$ and $\dist(y,E_\infty) \geq c_1 r(y)$, the Lipschitz estimaze is trivial again.

Before passing to the next step, we extend $g$ as a $C$-Lipschitz function $g : \R^n \to \R^n$.

Finally, we need to localize and make a one parameter family. 
Pick a continuous bump function $\varphi : [0,+\infty) \to [0,1]$,
such that $\varphi(t) = 1$ for $t \leq 1/2$, $\varphi(t) = 0$ for $t \geq 9/10$, and
$\varphi$ is affine in the intermediate interval. Then set
\begin{equation} \label{9a10}
    \rho(x,t) = t \varphi(r(x)/r_0) g(x) + (1 - t \varphi(r(x)/r_0)) \, x
\end{equation}
for $0 \leq t \leq 1$ and $x \in \R^n$ (notice that the definition of $g(x)$ outside $B_0$ does not matter since $\varphi(r(x)/r_0) = 0$ outside $B_0$).
It is straightforward that $\rho(x,0) = x$, which proves (\ref{9a1}).
When $r(x) \leq r_0/2$ and $t = 1$, we have $\rho(x,1) = g(x) \in E_\infty$, which proves (\ref{9a4}).
We have $\rho(x,t) = x$ when $r(x) \geq 9r_0/10$ because $\varphi(\abs{x}/r_0) = 0$ and also when $x \in E_\infty \cap B_0$ because $g(x) = x$. This proves (\ref{9a2}).
Finally, we note that $\rho$ is $C$-Lipschitz on $\R^n \times [0,1]$, the verification is standard and left to the reader.
\qed

\ms
\begin{rem}\label{r9a2}
    It would have been nice to also have that
    \begin{equation} \label{9a11}
        \rho(x,t) \in \Gamma \ \text{ when } x \in \Gamma,
    \end{equation}
    but in general this is not the case. More precisely, let $\Gamma_+$ and $\Gamma_-$ denote
    the two pieces of $\Gamma \cap B(0,r_0) \sm \{ 0 \}$. Only three cases can occur concerning 
    (a given) $\Gamma_\pm$. First, $\Gamma_\pm \cap E_\infty$ can be empty.
    Then \eqref{9a11} probably fails, but we don't care at all, because the sliding condition will be void 
    in $B_0$.
    The second simple case is when $\Gamma_\pm \subset E_\infty$. 
    Near $\Gamma_{\pm}$, we always have Case 4 (a half plane $\bH$), Case 1 (a plane $\bP$), Case 2 (a sharp $\bV$), Case 2 Bis (a truncated $\bY$) or Case 3 (a generic $\bV$); 
    in Cases 2 and 3, our construction already gives $\sigma(x,t) = \xi(t r(x))$ when $1/2 \leq t \leq 1$ 
    and $x = \xi(r(x))$ lies in 
    $\Gamma_\pm$, and by iteration we find that $\sigma(x,t) = \xi(t r(x)) \in \Gamma_\pm$
    for all $t \leq 1$. Hence \eqref{9a11} holds on $\Gamma_\pm$. In Case 1, we decided not to
    care, but when $\Gamma_\pm \subset E_\infty$ we can do better. In this case, since we
    know that all the points $\xi(r)$ lie in $E_\infty$, we can change our mind and treat 
    Case 1 exactly as Case 3, i.e., make sure that $\sigma(\xi(r),t) = \xi(t r(x))$, and 
    then we still have \eqref{9a11} in that case.

    We are left with the most unpleasant case when $\Gamma_\pm \cap E_\infty$
    and $\Gamma_\pm \sm E_\infty$ are both nonempty. This happens only as follows. 
    Let $x = \xi(r)$ be a point that lies in the closure of both sets (a point where $E_\infty$ 
    is leaving $\Gamma$).
    At such a point, the only option is
    Case 1, where we can take an approximation plane $P$ that contains $L_0$. 
    Indeed, at such a point $E_\infty$ can touch and leave $\Gamma$ more or less freely, 
    as long as this happens
    tangentially. We will need to find a way, when this happens, to use one of the pieces of
    $\Gamma_\pm \sm E_\infty$ to neutralize the sliding condition.
\end{rem}

\section{We glue the retractions} 

\label{S10}

So far we considered our limit set $E_\infty$, picked a point $x_0\in E_\infty$
(we rapidly assumed that $x_0 = 0$ to simplify the notation), and worked in any
small enough ball $B_0 = B(x_0,r_0)$ so that the good approximations of $E_\infty$
by cones hold, and eventually we obtained a mapping $\rho$ as in Lemma \ref{l9a1}.
Now we want to compose (a finite number of) these mappings, to obtain a deformation
of a neighborhood of $E_\infty$ onto a subset of $E_\infty$, which will be used to show that
this subset essentially is a sliding competitor for elements of our initial sequence $\{ E_k \}$
(see Section \ref{S2}).

\begin{lem}\label{l10a1}
    Let $E_\infty$ be the (coral) minimal set of the previous sections. 
    Then we can find a small number $\kappa > 0$, that may depend very badly on $E_\infty$,
    and a Lipschitz mapping $\Phi : \R^n \to \R^n$ such that
    \begin{equation}\label{10a2}
        \Phi(x) = x \text{ for } x\in E_\infty.
    \end{equation}
    and
    \begin{equation}
        \Phi(U(\kappa)) \subset E_\infty,
    \end{equation}
    where $U(\kappa) := \set{x\in \R^n | \dist(x,E_\infty) < \kappa}$.
\end{lem}

We need to restrict to a neighborhood of $E_\infty$, because $E_\infty$
could have some topology (with a boundary $\Gamma$ that follows $E_\infty$
sufficiently well, we could arrange that $E_\infty$ is a topological sphere in $\R^3$),
which would prevent the existence of a continuous retraction on $E_\infty$ defined everywhere.
But the reader should be warned that, in order to simplify the proof and avoid playing with 
coverings, we shall take $\kappa$ extremely small and the Lipschitz constant very large, also depending badly on $E_\infty$. We could possibly avoid that but the authors are not sure.

For each $x_0 \in E_\infty$, we can find $r_0 = r_0(x_0) > 0$ such that Lemma \ref{l4a1} 
and the ensuing construction gives a mapping $\rho$ as in Lemma \ref{l9a1}. 
Let us only recall the endpoint $\rho_0$, defined by $\rho_0(x) = \rho(x,1)$.
By compactness, we can find a finite family of points $x_i \in E_\infty$, $i = 1,\ldots,m$,
such that the balls $B_i = B(x_i, r_0(x_i)/4)$ cover $E_\infty$. 
Call $\rho_i : \R^n \to \R^n$ the corresponding mapping and take
\begin{equation} \label{10a3}
    \Phi = \rho_m \circ \ldots \circ \rho_1.
\end{equation}
We claim that if $\kappa$ is small enough, $\Phi$ does the job.
Of course $\Phi$ is Lipschitz (unfortunately, with a constant that depends on $m$ that we don't control), and \eqref{10a2} holds because every mapping $\rho_i$ fixes $E_\infty$. We just need to make sure that $\Phi(x) \in E_\infty$ when $x \in U(\kappa)$, for $\kappa$ sufficienty small.

For all $j=1,\ldots,m$, we define
\begin{equation}
    \phi_j = \rho_j \circ \ldots \circ \rho_1
\end{equation}
and we take the convention $\phi_0 = \mathrm{id}$.
We let $C_j$ denote the Lipschitz constant of $\phi_j$.
The function $(\phi_j - \mathrm{id})$ is $(C_j + 1)$-Lipschitz on $\R^n$ and is $0$ on $E_\infty$ so for all $x \in \R^n$,
\begin{equation}
    |\phi_j(x) - x| \leq (C_j + 1) \dist(x,E_\infty), 
\end{equation}
and in particular $|\phi_j - \mathrm{id}| \leq (C_j + 1) \kappa$ on $U(\kappa)$.

We take $\kappa$ so small that $(C_j +1) \kappa < r_0(x_j)/100$ for all $j = 1,\ldots,m$. 
For each $x\in U(\kappa)$, there is a point $x' \in E_\infty$ such that $|x - x'| \leq \kappa$ and there exists $i = 1,\ldots,m$ such that $x' \in B(x_i,r_0(x_i)/4)$ so $x \in B_i = B(x_i,r_0(x_i)/3)$.
As $\abs{\phi_{i-1}(x) - x} \leq (C_{i-1} + 1) \kappa$, it follows that $y : = \phi_{i-1}(x) \in B(x_i, r_0(x_i)/2)$ and by \eqref{9a4}, $\phi_{i}(x) = \rho_i(y) \in E_\infty$.
From there on, all the successive images are equal to $\Phi_{i}(x) \in E_\infty$ (by \eqref{9a2}).
The Lemma is proved.
\qed

\ms
We are now ready for the last part of the argument.
Recall that we started in the statement of Theorem \ref{t2a1}, by taking a minimizing sequence 
$\{ E_k \}$ for the functional of \eqref{2a5} in the class $\cE(E_0,\Gamma)$ of the early introduction. 
Then we took the measures $\mu_k = J_{\vert E_k}$, extracted a 
subsequence that converges, and proved in Theorem \ref{t2a1} that
the limit $\mu_\infty$ is of the form $\H^{2}_{\vert E_\infty}$ for some sliding
minimizer $E_\infty$. But we observed that this is not enough for Theorem~\ref{t1}:
we now want to show that $E_\infty$ gives a minimizer in the initial class
$\cE = \cE(E_0,\Gamma)$. 

We have to say ``gives'', because $E_\infty$ is probably not in the class $\cE$ itself,
because some of the topology of $E_k$ may have disappeared when we took the limit, so 
what we want to prove is merely that there is another set $E$ such that 
\begin{equation} \label{10a4}
    E \in \cE \text{ and } \H^2(E \sm E_\infty) = 0.
\end{equation}
Thus the true minimizer is $E$, and $E$ will be typically composed of $E_0$, plus
a bunch of wires of dimension $1$ that remember the topology but have no measure.
We seem to allow $E \cap E_\infty$ to be strictly smaller than $E_\infty$, but this won't
happen, because the facts that $J(E_\infty) \leq m_0 \leq J(E)$ (see \eqref{1a5}) 
and $\H^2(E \sm E_\infty) = 0$ imply that
$J(E_\infty) = J(E) = J(E \cap E_\infty)$, hence $\H^2(E_\infty \sm E) = 0$.
Since $E_\infty$ is a coral set and $E$ is closed, this implies that $E_\infty \subset E$.

Let us first prove the existence of $E$ as if the sliding condition did not exist.
Let $\kappa$ be as in Lemma \ref{l10a1} (the thickness of the neighborhood of $E_\infty$ where we define $\Phi$). This parameter depends only on $n$ and $E_\infty$ and is fixed for the rest of the proof.
We use the fact that $\mu_\infty$ is the limit of the $\mu_k$ and is supported by $E_\infty$ to pick $k$ so large that $\H^2(E_k \sm U(\kappa/3)) \leq C \mu_k(\R^n \sm U(\kappa/3)) \leq \varepsilon$,
where $\varepsilon$ will be chosen very soon. 
Note that the index $k$ is now fixed for the rest of the proof.

\newcommand{\phiFF}{p} 

Next we choose a dyadic scale $\tau := 2^{-l}$, where $\tau$ is chosen so small that $10\sqrt n \tau < \kappa$.
The parameter $\tau$ depends on $n$, $\kappa$ but not $\varepsilon$, and we will later take $\varepsilon$ small enough compared to $\tau$.
We let $\Delta$ be the set of all closed dyadic cubes of side length $\tau = 2^{-l}$.
Notice that they all have a diameter $\leq \kappa/10$ by our choice of $\tau$.

We let $\mathcal{S}$ be the set of thoses closed cubes in $\Delta$ that meet $\R^n \setminus U(\kappa/2)$.
Therefore,
\begin{equation}
    \R^n \setminus U(\kappa/2) \subset \bigcup_{Q \in \mathcal{S}} Q \subset \R^n \setminus U(\kappa/3),
\end{equation}
and in particular the cubes $Q \in \mathcal{S}$ cannot meet a cube $Q' \in \Delta$ which contains a point of $E_{\infty}$.
We then perform a Federer-Fleming projection $\phiFF$ of $E_k$ into the $2$-dimensional faces of cubes $Q \in \mathcal{S}$.

There is a quite general formalism to define Federer-Fleming projections in Section 2 of \cite{La} but in summary, it is a map $\phiFF : \bigcup \set{Q | Q \in \mathcal{S}} \to \R^n$ which performs successive projections of the set $E \cap \bigcup_{Q \in \mathcal{S}} Q$ in the faces of the grid until it is projected onto the $1$-dimensional skeleton of the grid.
The map $\phiFF$ preserves all faces of cubes $Q \in \mathcal{S}$ (whether they are $0$-faces, $1$-faces, $2$-faces, etc.) and the projections center can be chosen in such a way that for all cube $Q \in \mathcal{S}$, we have $\H^2(\phiFF(E_k \cap Q)) \leq C \H^2(E_k \cap Q)$.

Let us underline a difference between our construction and the Federer-Fleming projection in (\cite{DS2}).
In \cite{DS2}, the Federer-Fleming projection of a $2$-dimensional set $E$ in a cube $Q_0$ consists in subdividing $Q_0$ into a grid of smaller cubes $Q$ and then define $\phiFF$ by doing projections of $E \cap Q$ in the internal faces of the grid, but not the external faces so that $\phiFF = \mathrm{id}$ on $\partial Q_0$. In our situation, we don't need $\phiFF = \mathrm{id}$ on the boundary of $\bigcup_{Q \in \mathcal{S}} Q$ so we can perform the projections in all faces of cubes, without distinguishing internal and external faces.

Using the fact that the cubes $Q \in \mathcal{S}$ have finite overlap, we can estimate
\begin{equation}
    \H^2\left(\phiFF\left(E_k \cap \bigcup_{Q \in \mathcal{S}} Q\right)\right) \leq C \sum_{Q \in \mathcal{S}} \H^2(E_k \cap Q)
    \leq C \H^2(E_k \setminus U(\kappa/3))
    \leq \varepsilon
\end{equation}
and we can assume $\varepsilon$ small enough depending on $\tau$ so that the image cannot contain any full $2$-face of a cube $Q \in \mathcal{S}$.
This allows to project again on $1$-faces (the resulting projection still denoted by $\phiFF$) so that
\begin{equation}
    \H^2\left(\phiFF\left(E_k \cap \bigcup_{Q \in \mathcal{S}} Q\right)\right) = 0.
\end{equation}

Since $\bigcup \set{Q | Q \in \mathcal{S}}$ is a positive distance from $E_{\infty}$, we can extend $\phiFF$ as a Lipschitz function $\phiFF : \R^N \to \R^N$ in such a way that $\phiFF = \mathrm{id}$ on $E_\infty$. We can even arrange so that $\phiFF$ preserve all face of cubes in $\Delta$ (first set $\phiFF = \mathrm{id}$ at any vertex of $\Delta$ which does not belong to a cube $Q \in \mathcal{S}$, next interpolate linearly on the edges (1-faces) that don't belong to a cube $Q \in \mathcal{S}$, then interpolate on the remaining 2-faces, and so on...). 
As a consequence, since the cubes $Q \in \Delta$ have diameter $\leq \kappa/10$, we have $\abs{\phiFF - \mathrm{id}} \leq \kappa/10$ in $\R^N$.

The restriction of $\phiFF$ to $E_k$ is a Lipschitz deformation of $E_k$. Let us assume that it is even a sliding deformation of $E_k$. In this way, $F_k := \phiFF(E_k)$ lies in $\cE = \cE(E_0,\Gamma)$ too.
By construction, $\phiFF(U(\kappa / 2)) \subset U(\kappa)$, because $\abs{\phiFF - \mathrm{id}} \leq \kappa/10$, whereas $\phiFF(E_k \setminus U(\kappa / 2))$ is $\H^2$-negligible (it is even contained in a one dimensional grid).
Therefore, the set $F_k = \phi_FF(E_k)$ is composed of one piece $F_k^1 \subset U(\kappa)$ and another piece $F_k^2$ which is negligible

Now, if the map $\Phi$ from Lemma \ref{l10a2} were 
also a sliding deformation, the image $E := \Phi(F_k) = \Phi \circ \phiFF(E_k)$ would lie in $\cE$ too.
But
\begin{equation}\label{eq_phi_Fk}
    \Phi(F_k) =  \Phi(F_k^1) \cup \Phi(F_k^2) \subset E_\infty \cup \Phi(F_k^2),
\end{equation}
and $\H^2(\Phi(F_k^2)) = 0$ because $\Phi$ is Lipschitz, so \eqref{10a4} holds, as needed.

\ms
Unfortunately we need to take care of the sliding condition.
The mild constraint on $\Gamma$ can be used to build a Federer-Fleming projection 
$\phiFF$ which preserves $\Gamma$.
This is obvious if $\Gamma$ is a finite union of edges in a grid of cubes, but it can also be done easily when $\Gamma$ is the biLipschitz image of such a set, hence also when $\Gamma$ is as in Theorem \ref{t1}.


We prove a variant of Lemma \ref{l10a1} where $\Phi$ is a (global) sliding deformation.
Let us first define what is a \ub{global sliding deformation}.
Here, this is a Lipschitz map $\Phi : \R^n \to \R^n$ such that there exists a Lipschitz homotopy $F : [0,1] \times \R^n \to \R^n$ satisfying $F_0 = \mathrm{id}$, $F_1 = \Phi$ and $F_t(\Gamma) \subset \Gamma$ for all $t$.
It is usually also required that there exists a compact set $C \subset \R^n$ such that 
$F_t = \mathrm{id}$ in $\R^n \setminus C$, for all $t$.
This last condition would not play a constraining role in the statement of Lemma \ref{l10a2} because, since $\Gamma$ and $E_\infty$ are compact sets,
it is always possible to artificially set $F_t = \mathrm{id}$ away from $\Gamma$ and $E_\infty$.
One recognizes the same definition as in the beginning of Section \ref{S1}, but with $E_0$ replaced by $\R^n$.
The advantage of global sliding deformations is that for any compact set $E$, $\Phi$ 
induces
a sliding deformation of $E$. Hence, in the reasoning above, the image $\Phi(F_k)$ 
would lie in $\cE$ even though the construction of $\Phi$ is independent from $F_k$.

\begin{lem}\label{l10a2}
    Let $E_\infty$ be the (coral) minimal set of the previous sections.
    Then we can find a small number $\kappa > 0$, that may depend very badly on $E_\infty$,
    and a global sliding deformation $\Phi : \R^n \to \R^n$ such that
    \begin{equation}
        \Phi(U(\kappa)) \subset E_\infty,
    \end{equation}
    where $U(\kappa) := \set{x \in \R^n | \dist(x,E_\infty) < \kappa}$.
\end{lem}
Here, $\Phi$ is not a retraction as in Lemma \ref{l10a1} because we lost the property 
$\Phi = \mathrm{id}$ in $E_\infty$.
This is unavoidable. The set $E_\infty$ may leave $\Gamma$ tangentially and in this case, 
there will be a sequence of points $x_k \in \Gamma \setminus E_\infty$ for which the ratio 
$\mathrm{dist}(x_k,E_\infty \cap \Gamma) / \mathrm{dist}(x_k, E_\infty)$ goes to $\infty$.
A retraction $\Phi$ on $E_\infty$ that preserves $\Gamma$ would send these points on $E_\infty \cap \Gamma$.
If the retraction is $C$-Lipschitz and fixes the points of $E_\infty$, we must have $\mathrm{dist}(x_k, \Gamma \cap E_{\infty}) \leq \abs{\Phi(x_k) - x_k} \leq C \mathrm{dist}(x_k,E_\infty)$ and this is incompatible with the above ratio going to $\infty$.
On the other hand, if $E_\infty$ does not leave $\Gamma$ tangentially, Case 1 
does not bother and one can adapt the construction of $\sigma$ in Section \ref{S8a} so that the 
functions $\rho$ in Lemma~\ref{l9a1} preserves the boundary (see Remark~\ref{r9a2}). 
In this case, the construction given by Lemma~\ref{l10a1} is directly a sliding deformation.
We underline that the property $\Phi = \mathrm{id}$ in $E_\infty$ was not needed for the reasoning above; we just needed $\Phi(U(\kappa)) \subset E_\infty$ to establish (\ref{eq_phi_Fk}). Therefore, Lemma \ref{l10a2} will complete the proof of our existence theorem.
\begin{proof}
    Recall that $\Gamma$ is a 
    compact set composed of a finite number of closed smooth loops $\Gamma_j$, 
    $1 \leq j \leq j_{max}$, that lie at positive distances from each other.
    For $\delta > 0$, we set
    \begin{equation}
        \Gamma(\delta) := \set{x \in \R^n | \mathrm{dist}(x,\Gamma) \leq \delta}.
    \end{equation}

    The construction of a sliding deformation will use the following fact.
    There exists a small number $\delta_0 > 0$ (depending only on $\Gamma$) and a $2$-Lipschitz map $\pi \colon \Gamma(\delta_0) \to \Gamma$ such that $\pi = \mathrm{id}$ on $\Gamma$.
    In particular, we have
    \begin{equation}\label{eq_pidist}
        \abs{\pi - \mathrm{id}} \leq 3 \mathrm{dist}(\cdot,\Gamma)
    \end{equation}
    because $(\pi - \mathrm{id})$ is $3$-Lipschitz in $\Gamma(\delta_0)$ and is $0$ on $\Gamma$.
    If the $\Gamma_i$ are $C^2$ or better, we can take the closest point projection,
    but even when the $\Gamma_i$ are just $C^{1+\varepsilon}$, we can define $\pi$ easily. 

    Before building $\Phi$, we make a general remark about how to build global sliding deformations.  Let us use
    $\pi$ to show that any Lipschitz map $\varphi \colon \R^n \to \R^n$ which preserves 
    $\Gamma$ and satisfies $|\varphi - \mathrm{id}| \leq \delta_0$ is a global sliding deformation.
    Indeed, we start by setting
    \begin{equation} \label{10c8}
        F(t,x) =
        \begin{cases}
            x                   &   \text{in} \ \set{0} \times \R^n\\
            \varphi(x)                &   \text{in} \ \set{1} \times \R^n\\
            \pi((1-t) x + t \varphi(x))   &   \text{in} \ [0,1] \times \Gamma.
        \end{cases}
    \end{equation}
    The composition with $\pi$ in the last formula is well-defined because for $x \in \Gamma$ and $t \in [0,1]$,
    \begin{equation}
        \mathrm{dist}((1-t) x + t \varphi(x),\Gamma) \leq t \abs{\varphi(x) - x} \leq \delta_0.
    \end{equation}
    These formulas coincide at the intersections of their domains since $\varphi$ preserves $\Gamma$ and $\pi = \mathrm{id}$ on $\Gamma$. In addition, the function is continuous because each piece is continuous on a closed domain. We can finally make a continuous extension $F \colon [0,1] \times \R^n \to \R^n$ by using the Tietze extension theorem. 
    As $\varphi$ is Lipschitz, then it is easy to check that we can make $F$ Lipschitz too,
    because \eqref{10c8} defines a Lipschitz mapping on $\left(\set{0} \times \R^n\right) \cup \left(\set{1} \times \R^n\right)
    \cup \left([0,1] \times \Gamma\right)$.
    The key point is to estimate that for $x,y \in \R^N$ and $s,t \in [0,1]$,
    \begin{equation}
        \big|\left[(1-t)x + t\varphi(x)\right] - \left[(1-s) x + s \varphi(x)\right]\big| \leq \abs{t - s} \abs{\varphi(x) - x} \leq \delta_0 \abs{t - s}
    \end{equation}
    and
    \begin{equation}
    \big|\left[(1-s)x + s\varphi(x)\right] - \left[(1-s) y + s \varphi(y)\right]\big| 
    \leq \abs{\varphi(x) - \varphi(y)} + \abs{x - y}.
    \end{equation}
    We can then use the Whitney theorem to extend this map in a Lipschitz way.
    This proves that $\varphi$  is in fact a global sliding deformation.

    We let $0 < \delta \leq \delta_0$ be a constant that will be fixed small later. 
    The construction of $\Phi$ relies on the following preparation map.
    We claim that there exist a relative open subset $V \subset \Gamma$ of $\Gamma$ containing $\Gamma \cap E_\infty$ and a Lipschitz map $\psi \colon \Gamma \to \Gamma$ such that $\abs{\psi - \mathrm{id}} \leq \delta/2$ and
    \begin{equation}\label{eq_psi}
        \psi(V) \subset \Gamma \cap E_\infty.
    \end{equation}
    We postpone the details to the end of the proof, but mention that in the easy case 
    when $\Gamma \setminus E_\infty$ is a finite union of disjoint open intervals,
    the idea would simply be to push points near $\Gamma \cap E_\infty$ to 
    $\Gamma \cap E_\infty$. 
   We are  not annoyed by topology, because we don't need to define $\psi$ on the whole open interval so that it lands on $E_\infty$.

    We assume temporarily that such a preparation map $\psi$ exist, and we pass to the construction of $\Phi$.
    We will take Lipschitz extensions repeatedly using, for example, the McShane-Whitney formula (we don't need to preserve the Lipschitz constants).
    Before that, we justify that there exists $0 < \eta \leq \delta_0$ such that for all $x \in \R^n$ satisfying $\mathrm{dist}(x,\Gamma) \leq \eta$ and $\mathrm{dist}(x,E_\infty) \leq \eta$, we have
    \begin{equation}\label{eq_psiV}
        \pi(x) \in V.
    \end{equation}
    The set $\Gamma \setminus V$ is compact and disjoint from the closed set $E_\infty$ so there exists a small $\eta > 0$ such that for all $x \in \Gamma \setminus V$, we have $\mathrm{dist}(x,E_\infty) > 4\eta$.
    By contraposition, for all $x \in \Gamma$, the condition $\mathrm{dist}(x,E_\infty) \leq 4\eta$ implies $x \in V$.
    Thus, for all $x \in \R^n$ satisfying $\mathrm{dist}(x,\Gamma) \leq \eta$ and $\mathrm{dist}(x,E_\infty) \leq \eta$, we have $\pi(x) \in \Gamma$ (because $\eta \leq \delta_0$) and $\mathrm{dist}(\pi(x),E_\infty) \leq 4 \eta$ (by (\ref{eq_pidist})) so $\pi(x) \in V$. The constant $\eta$ depends only on $n$, $\Gamma$, $E_\infty$, $\delta$ and is fixed for the rest of the proof.

    We extend $\psi$ in a neighborhood $\Gamma(\eta)$ of $\Gamma$ by setting 
    \begin{equation}\label{10c11}
        \psi(x) = \psi(\pi(x)) \ \text{ for } x \in \Gamma(\eta).
    \end{equation}
  Remember that $\abs{\psi - \mathrm{id}} \leq \delta/2$ in $\Gamma$. 
  We can assume $\eta \leq \delta/6$ so that we still have 
  $\abs{\psi - \mathrm{id}} \leq \delta$ in $\Gamma(\eta)$; indeed for $x \in \Gamma(\eta)$,
    \begin{align}
        \abs{\psi(x) - x}   = \abs{\psi(\pi(x)) - x}
        \leq \abs{\psi(\pi(x)) - \pi(x)} + \abs{\pi(x) - x}
        \leq \delta/2 + 3 \eta.
    \end{align}
     We finally extend $\psi$ as a Lipschitz function over $\R^n$ such that $\abs{\psi - \mathrm{id}} \leq \delta$. 

    Let us return to the construction of $\Phi$.
    According to Lemma \ref{l10a1}, there exists a small constant $\kappa_0$ and a Lipschitz map 
 $\Phi_0 \colon U(\kappa_0) \to \R^n$ such that $\Phi_0 = \mathrm{id}$ on $E_\infty$ and 
 $\Phi_0(U(\kappa_0)) \subset E_\infty$ where $U(\kappa_0) = \set{x \in \R^n | \dist(x,E_\infty) < \kappa_0}$.
    We can also assume that $\abs{\Phi_0 - \mathrm{id}} \leq \delta_0/2$ by taking $\kappa_0$ a bit smaller if necessary.
     Then can extend $\Phi_0$ in a Lipschitz way over $\R^n$ and such that $\abs{\Phi_0 - \mathrm{id}} \leq \delta_0/2$.

    We define $\Phi := \Phi_0 \circ \psi : \R^n \to \R^n$, which is our 
   candidate for proving the lemma (but we again will to make a few additional changes).
    For the rest of the proof, we fix $\delta \leq \min(\delta_0/2, \kappa_0/2)$. Since $\abs{\psi - \mathrm{id}} \leq \delta$, we have in particular $\abs{\Phi - \mathrm{id}} \leq \delta_0$. Thus we will be able to post-compose $\Phi_{\vert \Gamma}$ with $\pi$ later on.

    We fix a constant $0 < \kappa \leq \min \left(\eta, \kappa_0/2\right)$ and check that $\Phi$ sends $U(\kappa)$ to $E_\infty$.
    We have $\abs{\psi - \mathrm{id}} \leq \delta \leq \kappa_0/2$ and $\kappa \leq \kappa_0/2$ so the map $\psi$ sends $U(\kappa)$ in $U(\kappa_0)$. As $\Phi_0$ sends $U(\kappa_0)$ in $E_\infty$, the map $\Phi$ sends $U(\kappa)$ in $E_\infty$.

    Next, we look at the sliding condition. For $x \in U(\kappa) \cap \Gamma$, we have $x \in V$ 
    (this comes from (\ref{eq_psiV}) since $x \in \Gamma$ and $\mathrm{dist}(x,E_{\infty}) \leq \kappa \leq \eta$) so $\psi(x) \in E_\infty \cap \Gamma$. Since $\Phi_0$ fixes $E_\infty$, the map $\Phi$ sends $U(\kappa) \cap \Gamma$ in $\Gamma$.
    Unfortunately, we may not have $\psi(\Gamma) \subset \Gamma$ because of the lack of control of $\psi$ on $\Gamma \setminus U(\kappa)$, so we need to tweak a bit the definition of $\Phi$.
 For this purpose, we aim to replace $\Phi$ on $\Gamma$ by $\pi \circ \Phi$, but we should check carefully that this still is a well-defined Lipschitz deformation.
   
    Let $L$ be a positive constant. We consider the set
    \begin{equation}\label{eq_W}
        W(\kappa) := \set{x \in U(\kappa) | \mathrm{dist}(\Phi(x),\Gamma) \leq L \mathrm{dist}(x,\Gamma)}
    \end{equation}
    and show that $W(\kappa)$ is a neighborhood of $E_\infty$ when $L$ is big enough.
    For $x \in E_\infty$, we distinguish two cases. 
  If $\mathrm{dist}(x,\Gamma) < \eta$, then $\pi(x) \in V$ (by (\ref{eq_psiV})) so $\psi(x) = \psi(\pi(x)) \in E_\infty \cap \Gamma$ (by (\ref{eq_psi})) and then $\Phi(x) = \psi(x)$ so $\mathrm{dist}(\Phi(x),\Gamma) = 0$.
    Actually, the same reasoning applies to every point in
    \begin{equation}
        \set{y \in \R^n | \mathrm{dist}(y,E_\infty) < \eta,\ \mathrm{dist}(y,\Gamma) < \eta},
    \end{equation}
    which is a neighborhood of $x$, in which $\mathrm{dist}(\Phi(y),\Gamma) = 0$.
    If $\mathrm{dist}(x,\Gamma) > \eta/2$, then we bound directly
    \begin{equation}
        \mathrm{dist}(\Phi(x),\Gamma)   \leq \mathrm{dist}(x,\Gamma) + \abs{\Phi(x) - x}
        \leq \mathrm{dist}(x,\Gamma) + \delta_0
        \leq \mathrm{dist}(x,\Gamma) + 2\eta^{-1} \delta_0 \mathrm{dist}(x,\Gamma)
    \end{equation}
    so we can just take $L = 2\eta^{-1} \delta_0$ in this case. 
    We conclude that every point $x \in E_{\infty}$ has a neighborhood contained in $W(\kappa)$.
    We finally define
    \begin{equation}
        \Phi_1 =
        \begin{cases}
            \Phi &\text{in $W(\kappa)$}\\
            \pi \circ \Phi &\text{in $\Gamma$}.
        \end{cases}  
    \end{equation}
    The composition with $\pi$ in second formula is well-defined because 
    $\abs{\Phi - \mathrm{id}} \leq \delta_0$.
    The formulas coincide at the intersection of the domains because $\Phi$ sends 
   $W(\kappa) \cap \Gamma \subset U(\kappa) \cap \Gamma$ into $\Gamma$.
   Since $W(\kappa)$ is neighborhood of $E_{\infty}$, which is compact, there exists $\kappa_1 > 0$ such that $U(\kappa_1) \subset W(\kappa)$ and thus $\Phi_1(U(\kappa_1)) \subset E_{\infty}$.
    Let us check that $\Phi_1$ is Lipschitz, 
    using the fact that $\Phi$ is Lischitz for some constant $C \geq 1$, the definition (\ref{eq_W}) of $W(\kappa)$ and the fact that $\abs{\pi - \mathrm{id}} \leq 3 \mathrm{dist}(\cdot,\Gamma)$.
    For $x \in W(\kappa)$ and for $y \in \Gamma$, we have indeed
    \begin{align}
           \abs{\Phi(x) - \pi(\Phi(y))} &\leq \abs{\Phi(x) - \Phi(y)} + \abs{\Phi(y) - \pi(\Phi(y))}  \nonumber\\
        &\leq C \abs{x - y} + 3 \mathrm{dist}(\Phi(x),\Gamma)  
        \leq C \abs{x - y} + 3 L \mathrm{dist}(x,\Gamma)  
        \leq (C + 3 L) \abs{x - y}.  
    \end{align}
 Since $\abs{\pi - \mathrm{id}} \leq 3\delta_0$ on $\Gamma(\delta_0)$ and $\abs{\Phi - \mathrm{id}} \leq \delta_0$, we have $\abs{\Phi_1 - \mathrm{id}} \leq 4\delta_0$ on its domain of definition. 
 We can even replace $\delta_0$ by $\delta_0/4$ in the whole proof so that $\abs{\Phi_1 - \mathrm{id}} \leq \delta_0$.
    We extend one last time $\Phi_1$ as a Lipschitz function over $\R^n$ such that $\abs{\Phi - \mathrm{id}} \leq \delta_0$.
    Now, $\Phi_1$ induces a global sliding deformation as we observed at the beginning of the proof (just after the definition of $\pi$).
    This completes the proof of Lemma \ref{l10a2}, modulo the following verification.

    \ms
    As promised, we now
    detail the construction of the preparation map $\psi$. We shall do the construction 
    concerning one of the $\Gamma_j$, but then we shall do
    the same thing with each $\Gamma_j$ (and the constructions will be independent).
    We write $\Gamma$ for $\Gamma_j$ to simplify the notation. 

    The main point will be to reduce to the simple situation where $\Gamma \cap E_\infty$
    and $\Gamma \setminus E_\infty$ have a finite number of connected components. 

    Let $\delta > 0$. First select a (necessarily finite) maximal family $\{ z_i \}$ of points of $E_\infty \cap \Gamma$, with
    $|z_i-z_j| \geq \delta/20$ for $i \neq j$. Call $\{ J_k \}$, $k \in K$, the connected components of $\Gamma \sm \cup_{i} \{ z_i \}$. If $E_\infty \cap \Gamma = \emptyset$, we can take 
    $\psi(x) = x$ and $V = \emptyset$, so let us assume that
    $E_\infty \cap \Gamma \neq \emptyset$; then each $J_k$ is an open interval 
    of $\Gamma$, and the two endpoints of $J_k$ lie in $E_\infty$.
    Let us twist a little the notation and write $J_k = (a_k, b_k)$, with $a_k, b_k \in E_\infty$. 
    We intend to take $\psi(z_i) = z_i$ for all $i$, so we  need to define
    $\psi$ on each $J_k$, so that $\psi(a_k) = a_k$ and $\psi(b_k)=b_k$.

    If $J_k \subset E_\infty$, we keep $\psi(x) = x$ on $J_k$. Next suppose that $J_k$ 
    contains a point of $\Gamma \sm E_\infty$ and that the length of
    $J_k$ is larger than $\delta/4$. Call $\wh J_k$ the set of points of $J_k$
    that lie at distance $\geq \delta/9$ from $a_k$ or $b_k$, $A_k$ the set 
    of points of $J_k$ that lie at distance $\leq \delta/10$ from $a_k$, and 
    $B_k$ the set of points of $J_k$ that lie at distance $\leq \delta/10$ from $b_k$. 
    We take $\psi(x) = x$ on $\wh J_k$, $\psi(x) = a_k$ on $A_k$, and
    $\psi(x) = b_k$ on $B_k$. On the two remaining short intervals of $J_k \sm (\wh J_k \cup A_k \cup B_k)$, we interpolate ``linearly''. Since all the points that move lie at distance $\leq \delta/9$ from $a_k$ or $b_k$,
    we get that $|\psi(x)-x| \leq \delta/9 \leq \delta$ for $x \in J_k$. In addition, by maximality of the $z_i$,
    all the points of $J_k$ that lie at distance less than $\delta/20$ of $\Gamma \cap E_\infty$ must lie in at distance $\leq \delta/10$ from one of the $z_i$ and thus must lie in $A_k \cup B_k$. They are sent to $a_k\in E_\infty$ or $b_k \in E_\infty$, in accordance with \eqref{eq_psi}.

    Finally, when $J_k$ contains some point $c_k \in \Gamma \setminus E_\infty$
    and its length is less than $\delta/4$, we pick a small segment $\wh J_k$ centered at
    $c_k$, so that its ``double'' $2 \wh J_k$ does not meet $E_\infty$, and choose $\psi$
    Lipschitz on $J_k$ so that $\psi(x)=x$ on $\wh J_k$, $\psi$
    takes the values $a_k$ and $b_k$ on the two intervals that compose $J_k \sm 2 \wh J_k$,
    and $\psi$ interpolates on the two remaining intervals that compose $2 \wh J_k \sm J_k$.
    Here, the fact that $|\psi(x)-x| \leq \delta/4 \leq \delta$ for $x \in J_k$ follows from the fact that $\psi$ preserves $J_k$, which is of length less than $\delta/4$.
    If $d_k$ denotes the distance from $2 \wh J_k$ to $\Gamma \cap E_\infty$, then all the points of $J_k$ that lie at distance less than $d_k$ from $\Gamma \cap E_\infty$ are sent to $a_k \in E_\infty$ or $b_k \in E_\infty$.

    It is easy to see that the different pieces that compose $\psi$ can be glued to compose
    $\psi : \Gamma \to \Gamma$, which is Lipschitz (with an estimate that can depend badly on
    the geometry of $E_\infty$ in $\Gamma$, but this is all right). This completes our construction of $\psi$; Lemma \ref{l10a2} follows.

 \end{proof}

\section{Variants of the main theorem}
\label{Slast}
In this section we discuss generalizations of Theorem \ref{t1} that can be deduced from 
Theorem~\ref{t3a1}. The main point is that the proof has some flexibility, due to the fact
that $E_\infty$ is a (sliding) almost minimal set. 

Consider the problem where we minimize a functional $J$ as in \eqref{2a5}-\eqref{2a7}
in the class $\cE(E_0,\Gamma)$. Then we can use the conjunction of Theorems~\ref{t2a1}
and \ref{t3a1} to get the existence of solutions, but this requires finding a compact set
$K$, that contains $\Gamma$, where we know that we can construct a minimizing sequence
(or a limit of a minimizing sequence), and so that $\Gamma$ has a good access to the 
complement of $K$.

If we work with the functional $\H^2$, the best way, and probably the only reasonable way
to get $K$ is to take the convex hull of $\Gamma$, and ask for the good access. 
If we use a slightly different functional $J$, we may have a little more flexibility, 
because we can try to play with the definition of $J$ to show that we can find minimizing sequences
in a slightly more complicated $K$. Yet the best way to arrange this is to have a projection
$\pi : \R^n \to K$ such that $J(\pi(E)) \leq J(E)$ for any set $E$; in the case of $\H^2$,
this would be the shortest distance projection. More complicated arguments can exist 
(for instance, saying that if $E \in \cE(E_0,\Gamma)$
cannot be projected on $K$ as above, then $J(E)$
is too large for some other reason), but the existence of $K$ and $\pi$ is probably our best chance.
Notice that the good access condition seems to allow more
complicated shapes for $K$.

Of course forcing $E$ a priori to lie in $K$ (with the good access property) does not do the  
job, because $E_\infty$ will then only be almost minimal under the constraint that 
$E_\infty \subset K$, and the regularity results of \cite{Dvv} won't hold.

Similarly, we did not mention elliptic integrands $J$ given by a formula like
$J(E) = \int_E f(x,T_E(x)) d\H^d(x)$, with a function $f$ that depends also 
on the direction of the approximate tangent to $E$ at $x$, because they are not so easy 
to treat. If $E_0$ is rectifiable, then $T_E(x)$ exists almost everywhere, but the issue is
whether $E_\infty$ is an almost minimal set with a small enough gauge function. 
Unless $f$ is H\"older-continuous (as in our assumption \eqref{2a7}), or has a small enough
modulus of continuity, there is little chance that we can prove this.
When  $f$ depends also on the direction, we should probably require that this dependence is also
H\"older-continuous, say, but there is no analogue of \cite{Dvv} with that much generality.

In the special, still roughly Euclidean, case where $f(x,T) = |\det(A(x)\circ \pi_T)|$, 
where $x \to A(x)$ is a H\"older continuous mapping with values in $\cL(\R^n, \R^2)$, 
$\pi_T$ denotes the orthogonal  projection on the vector $2$-plane $T \subset \R^n$, 
and $A$ is such that $C^{-1} \leq |\det(A(x)\circ \pi_T)| \leq C$ for all $x$ and $T$ (ellipticity), there is a good chance that the results of this paper go through because those of
\cite{Dvv} could be extended. The point is that when $A$ is constant, 
almost minimizers for $J$ are the same as images under an invertible linear mapping of almost minimizers 
for $H^2$; then the idea would be, when we study an almost minimal set for $J$ at a point $x$, 
to conjugate with a linear mapping to reduce to the case of almost minimal sets. A true
argument would be more complicated than this, as we also need the regularity at the other points
nearby, where we would need to change the conjugation, and this was never written.

The case of truly non euclidean elliptic integrands, even independent of $x$, seems widely open
even far from the boundary.

\section{Local biLipschitz parameterizations of some almost minimal sets}
\label{Sbil}

The methods above allow us to give local biLipschitz parameterizations for coral 
almost minimal sets in some circumstances that we describe now. 
The statements will be cleaner with the following notion of local almost minimal sets. 
Let $\Omega \subset \R^n$ be open,
and let $E \subset \Omega$ and $\Gamma$ be closed in $\Omega$. 
We recall that a gauge is a nondecreasing function $h \colon (0,+\infty) \to [0,+\infty]$ such that $\lim_{r \to 0} h(r) = 0$. We shall assume in addition the condition (\ref{2a1}), that is, there exists $\alpha > 0$, $c_h \geq 0$ and $r_h > 0$ such that 
\begin{equation}
    h(r) \leq c_h r^\alpha \ \text{for $0 < r \leq r_h$}. 
\end{equation}
We say that $E$ is sliding almost minimal in $\Omega$, with boundary $\Gamma$ and gauge  function $h$, when \eqref{2a3} holds for every family $\{ \varphi_t \}$
that satisfies the conditions \eqref{1a1}-\eqref{1a4} with $E_0$ replaced by $E$, 
and in addition such that the deformation happens entirely in some 
$B(y,r) \subset \subset \Omega$,
as in \eqref{2a2}. That is, we keep most of Definition \ref{d2a1} but only require 
\eqref{2a3} for deformations that happen in balls contained in $\Omega$. 
We say that $E$ is (plain) almost minimal in $\Omega$, with gauge function $h$, when
it is sliding almost minimal in $\Omega$, with $\Gamma = \emptyset$. That is, we simply  forget about \eqref{1a2} and the constraint which prevents competitor from being degenerate comes instead from the fact that $\varphi(t,x) = x$ outside a compact subset of $\Omega$.

Finally, recall that the closed set  $E$ in $\Omega$ is called coral when 
$\H^d(E  \cap B(x,r)) > 0$  for all $x\in  E$  and $r > 0$.

We start with a statement for plain almost minimal sets.

\begin{thm}\label{t12a1}
    Suppose $E$ is a (plain) coral almost minimal set of dimension $2$ in an open set 
    $\Omega \subset \R^n$ that contains $B(0,1)$, with a gauge function $h$ that satisfies
    \eqref{2a1} for some $\alpha > 0$, $c_h \geq 0$ and $r_h > 0$.
    Suppose in addition that 
    $0 \in E$. Then there is a radius $r_0 \in (0,1/2)$, a minimal cone $X$, and a 
    biLipschitz mapping $\psi : X\cap B(0,r_0) \to E \cap B(0, r_0)$, such that in addition 
    \begin{equation} \label{12a1}
        |\psi(y)| = |y| \ \text{ for } y \in X\cap B(0,r_0).
    \end{equation}
\end{thm}

\ms
In fact, we will get the following more precise statement. 
There exist constants $C \geq 1$, that depends only on $n$,
and $\varepsilon > 0$, that depends only on $n$ and $\alpha$,
such that if $E$ satisfies the assumptions of the theorem, and if 
in addition $r_0 \in (0,1/2) \cap (0, r_h/2)$ is such that $c_h r_0^\alpha \leq \varepsilon$ and,  for
all $0 < r \leq 2r_0$, we can find a minimal cone $X = X(r)$ such that
\begin{equation} \label{12a2}
    d_{0,r}(E,X) \leq \varepsilon,
\end{equation}
then there is a $C$-biLipschitz mapping 
$\psi : X(r_0) \cap B(0,r_0) \to E \cap B(0, r_0)$ such that \eqref{12a1} holds true.
That is, we may even choose $X$ as one of the $X(r)$.

Theorem \ref{t12a1} will follow from this, and the fact that all the blow-up limits of $E$ at $0$ 
are minimal cones, and we even get that $\psi$ is $C$-biLipschitz and $X$ is a blow-up limit
of $E$ at $0$. The proof is the same as in Lemma \ref{l4a1}: one uses the definition of a blow-up limit (and the existence of convergent subsequences) to show that there exists $r_0 \in (0,1/2)$ such that 
for $0 < r \leq 2r_0$, \eqref{12a2} holds for some blow-up limit $X(r)$ of $E$ at $0$. The fact that
$X(r)$ is a minimal cone is standard.

Recall that with the same assumptions as in the theorem (or its precise form) 
\cite{Dhh} says that $E$ is locally equivalent to $X$ in $B(0,r_0)$ through a bi-H\"older mapping. 
This was better in a way, because we do not say here that $\psi$ extends to a local 
homeomorphism of $\R^n$. That is, we have a parameterization of $E$, but we do not
say that $E$ is nicely embedded in $\R^n$ with a Lipschitz mapping; at least the statement of 
\cite{Dhh} says that $E$ does not make weird knots in $\R^n$ near the origin. 

On the other hand the biLipschitz regularity is better than what we had in \cite{Dhh}, 
and with some extra work our mapping $\psi$ could also be made $C^{1+\varepsilon}$ 
away from $L_0$ and the set of $\bY$-points of $X$. 
This is only new when $n \geq 4$, because otherwise the result of \cite{Ta} is even more precise. 
Finally the additional property \eqref{12a1} may be convenient. 

\ms 
For the sliding almost minimal sets, we will need to forbid some blow-up limits. 
Let $L_0$ be a line through the origin. Denote by $\cX$ the set of sliding minimal cones centred at $0$
(of dimension $2$) with sliding boundary $L_0$. 
Let $\cB_0$ be the set of cones $X \in \cX$ such that for at least one of the two points
$\xi \in L_0 \cap \d B(0,1)$, $X$ coincides in some $B(\xi,\beta)$, $\beta > 0$ with a cone 
of type $\bY$ whose spine contains $L_0$. 
Since we will exclude this type of blow-up limit, we will be in able to apply Lemma \ref{lem_Xstructure}.
Similarly, we will exclude the set $\cB_1$ of cones $X \in \cX$ 
that have a point of  {\bf sharp } 
type $\bV$. This means that near at least one of the two points 
$\xi \in L_0 \cap \d B(0,1)$, $X$ coincides with a cone of type $\bV$, with a $2\pi/3$ angle.
Finally, 
let $\cB_2$ be the set of cones $X \in \cX$ such $X$ coincides, near at least one 
$\xi \in L_0 \cap \d B(0,1)$, with a plane $P$ that contains $L_0$.

\begin{thm}\label{t12a2}
    Suppose $E$ is a sliding coral almost minimal set of dimension $2$ in an open set 
    $\Omega \subset \R^n$ that contains $B(0,1)$, with a sliding boundary $\Gamma$ which is
    a $C^{1+\alpha}$ curve ($\alpha > 0$) through the origin, and  with a gauge function $h$ 
    that satisfies \eqref{2a1} for some $\alpha > 0$ and some $c_h \geq 0$. Let $L_0$ denote the 
    tangent line to $\Gamma$ at $0$.
    Suppose that $0 \in E$, and that no blow-up limit of $E$ at $0$ lies in 
    $\cB_0 \cup \cB_1$.
    Then there is a radius $r_0 \in (0,1/2)$, a sliding minimal cone $X$ associated to the sliding boundary
    $L_0$, and a biLipschitz mapping 
    $\psi : X\cap B(0,r_0) \to E \cap B(0, r_0)$, such that 
    \begin{equation} \label{12a3}
        |\psi(y)| = |y| \ \text{ for } y \in X\cap B(0,r_0).
    \end{equation}
    If in addition to that, no blow-up limit $X$ of $E$ at $0$ lies in $\cB_2$
    then also 
    \begin{equation} \label{12a4}
        \psi(y) \in \Gamma \ \text{ for } y \in X \cap L_0 \cap B(0,r_0).
    \end{equation}
\end{thm}

\ms
As we shall see from the proof, we can take for $X$ one of the blow-up limits of $E$ at $0$, and the
existence of $\psi$ will also show that all the blow-up limits of $E$ at $0$ are bilipschitz images of that
$X$; this does not imply that $E$ has a tangent cone at $0$, but only that all the blow-up limits are
bilipschitz equivalent to each other. Notice however that the biLipschitz constant for $\psi$ will depend on $E$ and the point $0$, through the constant $\beta > 0$ below, which itself depends on the list of blow-up limits of $E$ at $0$.

The theorem will follow from the following more precise result. 
Assume, on top of the hypotheses of the theorem
(we will not try to see to which extent these follow from the following),
that $r_0 \in (0,1/2)$, 
$\beta > 0$, $\varepsilon > 0$ such that 
$c_h r_0^\alpha \leq \varepsilon$,
and for all $0 < r \leq 2r_0$, one can find a sliding minimal cone 
$X = X(r)$ such that $d_{0,r}(E,X) \leq \varepsilon$ holds, 
and also such that 
\begin{equation}\label{12a5}
    \begin{gathered}
        \text{whenever $X(r)$ contains a point of type $\bV$,}\\
        \text{the faces along $L_0$ make an angle in the range $]2\pi/3 + \beta,\pi[$,}
    \end{gathered}
\end{equation}
and
\begin{equation} \label{12a6}
    \dist(y,L_0) \geq \beta \ \text{ for any triple junction } y \in X \cap \d B, 
\end{equation}
Then, if $\varepsilon$ is small enough, depending on $n$, $\alpha$,
and $\beta$,
there is a $C$-biLipschitz mapping $\psi : X(r_0) \cap B(0,r_0) \to E \cap B(0, r_0)$ 
such that \eqref{12a3} holds. Here $C$ depends only on $n$ and $\beta$.
If in addition to \eqref{12a5} and \eqref{12a6}, we also require that 
\begin{equation}\label{12a8}
    \begin{gathered}
        \text{if $X(r)$ coincides near some $\xi \in L_0$ with a cone of type $\bV$, then}\\
        \text{the angle along $L_0$ of the two faces of this cone is in the range $]2\pi/3 + \beta,\pi - \beta[$.}
    \end{gathered}
\end{equation}
and also
\begin{equation} \label{12a9}
    \dist(X(r), \xi) \geq \beta \ \text{ whenever } \xi \in L_0 \cap \d B(0,1) \sm X(r),
\end{equation}
then we also get the sliding condition \eqref{12a4}.

The theorem still does not say anything about the way $E \cap B(0,r_0)$
is embedded in $\R^n$, i.e., whether $\psi$ has a biLipschitz extension to $B(0,r_0)$,
for instance.

Note that $\psi$ is not necessarily smooth along $L_0$, in the sense that the angle of the faces
along $L_0$ could be slightly different at $\xi \in X \cap L_0$ than at $\psi(x) \in E \cap \Gamma$.
On the smooth part of $X$, we can make it $C^{1+\varepsilon}$, and even along the $\bY$-set
of $X$, we can arrange the mapping to preserve the angles. This is not surprising because we have local regularity results that control $E$ far from $0$, so we won't even try to enforce this aspect.

The point of removing the cones of $\cB_0$ is that this way we can use the description 
of \cite{Dvv}, as in the previous sections. 
We remove the cones of $\cB_1$ so that we never get Case 2
in Section \ref{S5} and the topology of $E$ stays the same as the topology of $X$ 
(and we will see that it remains the same for smaller radii). 
Finally, if we want the sliding condition \eqref{12a4}, we also forbid $\cB_2$ to avoid Case 1 above, 
which was the only case when $E$ could leave from $\Gamma$ and make it hard to respect \eqref{12a4}. 

We claim that we only need to prove the precise version of Theorem \ref{t12a2}.
For Theorem~\ref{t12a1} and its precise variant, just remove some cases and all the 
references to $\Gamma$ and $L_0$. For Theorem \ref{t12a2}, let us check now that its assumptions
imply the precise assumptions for small enough $r_0$.

Let $E$ be as in the theorem, and let $\cX_0 \subset \cX$ denote
the set of blow-up limits of $E$ at $0$. We first claim that for $r$ small, we can find 
$X(r) \in \cX_0$ such that $d_{0,r}(X(r), E) \leq \varepsilon(r)$, with 
$\lim_{r \to 0} \varepsilon(r) = 0$. The proof is as in Lemma \ref{l4a1}: 
otherwise, for some $\varepsilon > 0$ 
there are bad radii $r_k$ that tend to $0$ but so that we cannot find $X(r_k)$ such that 
$d_{0,r_k}(X(r_k), E) \leq \varepsilon$; we extract a subsequence
for which the set $r_k^{-1} E$ converge to a limit, this limit lies in $\cX_0$, 
and this contradicts the definition of $r_k$.

Next we want to check that for every blow-up limit $X_0$ of $E$ at $0$, there exists $\varepsilon > 0$ 
and $\beta > 0$ (depending on $E$ and $X_0$) such that whenever $X$ is a sliding minimal cone with $d_{0,2}(X,X_0) \leq \varepsilon$, then $X$
satisfies \eqref{12a5}, \eqref{12a6} and, if needed, \eqref{12a8} and \eqref{12a9}.  We recall that by assumption
$X_0 \notin \cB_0 \cup \cB_1$.
We proceed by contradiction; we first assume that for every $\beta = 2^{-k}$, there exists a sliding minimal cone $X_k$ such that $d_{0,2}(X_k,X_0) \leq 2^{-k}$ but \eqref{12a5} fails.
Each $X_k$ has a point of type $\bV$, with a nearly
sharp angle, and we extract a subsequence so that this is always the same point 
$\xi \in L_0 \cap \d B(0,1)$. The general description of minimal cones (see in particular Lemma \ref{lem_Xstructure}) says that there exists a uniform $c > 0$, that depends only on $n$, such that in $B(\xi,c)$, $X_k$ coincides with the same nearly
sharp cone of type $\bV$. Then, the limit $X_0$ coincides in $B(\xi,c)$ 
with a sharp $\bV$-cone but this contradicts the fact that $X_0 \notin \cB_1$.

For \eqref{12a6}, suppose now that there is a sequence $(X_k)$ that tends to $X_0$ such that
\eqref{12a6} fails for $\beta = 2^{-k}$.  
Since \eqref{12a6} fails, $X_k$ has a $\bY$-point 
$y_k$, such that $\dist(y_k,L_0) \leq 2^{-k}$.
We can take a subsequence so that this point is always on the same side of $L_0$, and so there exists
$\xi \in L_0 \cap \d B(0,1)$ such that $|y_k - \xi| \leq 2^{-k+1}$ 
(we do not exclude the case when $y_k = \xi$). By the general description of $X_k$, 
we then know that for some $c > 0$, that depends only on $n$, $X_k$ coincides in $B(\xi,c)$
with cone of type $\bY$ with a spine through $y_k$, or else with a cone of type $\bY$, truncated
by $L_0$, still with a spine through $y_k$. Indeed, letting $c_*$ denote the constant of Lemma \ref{lem_Xstructure}, we see that either the spine is at distance $< c_*/2$ from $L_0$ and then $X$ coincide with a truncated $\bY$ in $S(y_k,c_*/2)$, or the spine is at distance $> c_*/2$ from $L_0$ and then $X_k$ coincide with a $\bY$ cone in $S(y_k,c_*/20)$. So in all cases, the constant $c := c_*/20$ satisfies our claim.
Then, the limit $X_0$ coincides in $B(\xi,c)$ with a cone of type $\bY$ with a spine through $\xi$,
or else a sharp $\bV$-set. In both case this contradicts the fact that $X_0 \in \cX \sm \cB_0 \cup \cB_1$.

Finally, we check that in the case when we also excluded $\cB_2$, \eqref{12a8} and \eqref{12a9} 
hold too. If \eqref{12a8} fails, as before there is a sequence $(X_k)$ that tends to $X_0$ 
for which the $X_k$ coincide near some $\xi \in L_0 \cap \d B(0,1)$ (and we can suppose that
it is always the same point), with a set $V_k$ of type $\bV$ with an angle that tends to $\pi$. 
By the general description, and since by \eqref{12a6} there is no point of type $\bY$ nearby,
we get that $X(r_k)$ coincides with $V_k$ in some 
$B(\xi,c)$. Then we go to the limit and find that
after a new extraction, $X_0$ tends to a cone of $\cB_2$, a contradiction. 

The case when \eqref{12a9} fails 
(i.e., when some $\xi \in L_0 \cap \d B(0,1)\sm X(r)$ lies very close to $X(r)$)
is similar; now the general description of $X(r)$, plus the fact that we have no point of type 
$\bY$ around $\xi$,
says that $X(r_k)$ coincides in $B(\xi,c)$ with a plane that passes within $\beta_k = 2^{-k}$ of $\xi$, 
and the limit $X_0$ coincides with a plane in $B(\xi,c)$. Thus $X_0 \in \cB_{2}$, a contradiction.

So for every $X_0 \in \cX_0$ (a blow-up limit  of $E$ at $0$), we found $\varepsilon > 0$ and 
$\beta > 0$ (depending on $E$ and $X_0$) such that whenever $X$ is a sliding minimal cone with 
$d_{0,2}(X,X_0) \leq \varepsilon$, then $X$ satisfies the properties required for the precise asssumptions.
We can do a little better with compactnesss: find $\varepsilon > 0$ that does not depend on 
$X_0 \in \cX_0$. Indeed otherwise, there is a sequence $\{ X_{0,k}\}$ for which $\varepsilon = 2^{-k}$
does not work, which we can assume to converge to some $X_{0,\infty} \in \cX$, but then the $\varepsilon$
associated to $X_{0,\infty}$ works for $k$ large, a contradiction.

Because of this, for each $r$ small enough, we can take any cone 
 $X(r) \in \cX_0$ such that 
$d_{0,r} (E, X(r)) \leq \varepsilon(r)$,
with $\lim_{r \to 0} \varepsilon(r) = 0$, and for $r$ small enough,
$X(r)$ satisfies 
 the required assumptions. 

So the more precise assumptions \eqref{12a5} and \eqref{12a6}, 
and when needed \eqref{12a8} and \eqref{12a9} all hold, and Theorem \ref{t12a2} follows from the precise 
version, 
 that we shall prove now.

\ms
So let us assume that $E$ is as in the theorem, and that 
$r_0$ is so small that all our extra assumptions hold. For $0 \leq r \leq 2r_0$,
choose a minimal cone $X = X(r)$ such that \eqref{12a2}, 
and the correct combination of \eqref{12a5}-\eqref{12a9} holds.
Then (by \eqref{12a5}) we can use the description of $E$ 
that we had in Section \ref{S5}, with the additional information (coming from \eqref{12a5}) that
we never have Case 2 and, if we excluded $\cB_2$, we never have Case 1 either.

We need some notation. Set
\begin{equation} \label{12a10}
    r_k = 2^{-k}r_0, \  
    X_k = X(r_k), \ 
    A_k = \ol B(0, r_k) \sm B(0, r_{k+1}), \ 
    S_k = \d B(0,r_k), \  
    \Sigma_k = X_k \cap S_k.
\end{equation}
Our best model for $E \cap \overline B(0,r_0)$ is not really $X_0 \cap \overline  B(0,r_0)$, but a set 
$T = \cup_{k \geq 0} T_k$, where the tube $T_k \subset A_k$ will be build to connect
smoothly $\Sigma_k$ to $\Sigma_{k+1}$. We start with the usual description of $\Sigma_k$.
Since $X_k$ is a sliding minimal cone, we get that $\Sigma_k$ is the disjoint (except for the endpoints)
union of a finite collection of arcs of great circles, the $\cC_{i,j}$, $(i,j) \in I_a(r_k)$,
that go from one vertex $a_i$ of $\Sigma_k$ to another vertex $a_j$, plus a finite union of full
great circles $C_i$, $i \in I_c(r_k)$.
Away from $L_0$, the arcs $\cC_{i,j}$ can only meet by sets of three at their endpoints,
with $\frac{2\pi}{3}$ angles, and otherwise they stay at distances $\geq C^{-1} r_k$ 
from each other unless they have a common endpoint;
they also have a length at least $C^{-1}$. The circles $C_i$ stay at distances
$\geq C^{-1} r_k$ from each other, and from the $\cC_{i,j}$. And for $\xi \in L_0 \cap \Sigma_k$,
there is either a single $\cC_{i,j}$ that ends at $\xi$, or two arcs $\cC_{i,j}$ that ends at $\xi$,
and by \eqref{12a5} they make an angle larger than $\frac{2\pi}{3} + \beta$ at that point.
Because of \eqref{12a6}, we know that we cannot have three arcs leaving
from $\xi$, but also the arcs that leave from $\xi$ have a length at least $\beta r_k$. 
When $\cB_2$ is excluded, \eqref{12a8} says if two arcs $\cC_{i,j}$ leave from $\xi$, 
they make an angle smaller than $\pi-\beta$. 

Also, still when $\cB_2$ is excluded, \eqref{12a9} excludes the case of
any $\cC_{i,j}$ or a circle $C_i$ that comes within $\beta r_k$ from a point $\xi \in L_0 \cap \Sigma_k$ 
if they don't contain $\xi$.

Incidentally, two arcs $\cC_{i,j}$ that leave from some $\xi \in L_0$ in opposite directions 
may simply be parts of a circle $C_i$; we will not need to decide whether this counts as a curve 
or a circle.

All the sets $\Sigma_k$ have this type of description, and we may use the notation $\cC_{i,j}^k$,
or $C_i^k$, to point out the dependence on $k$. 
We are interested in how the descriptions for $\Sigma_k$ and $\Sigma_{k+1}$ fit, 
and let us observe that by \eqref{12a2}, the triangle inequality, and the fact that
$X_k$ and $X_{k+1}$ are cones, 
\begin{equation} \label{12a11}
    d_{0,2}(X_k \cap \d B(0,1), X_{k+1} \cap \d B(0,1)) \leq C \varepsilon.
\end{equation}
Because of this, and modulo a small exception that will be discussed soon, there is a way
to index the vertices $a_i$, and the arcs $\cC_{i,j}$ and circles $C_j$, for $k$ and $k+1$, 
so that 
\begin{equation} \label{12a12}
    |a_i^k - 2a_i^{k+1}| \leq C \varepsilon r_k, \ 
    d_{0,2r_k}(\cC_{i,j}^k, 2 \cC_{i,j}^{k+1}) \leq C \varepsilon, \ 
    d_{0,2r_k}(C_j^k, 2 C_j^{k+1}) \leq C \varepsilon.
\end{equation}
This is easy (but a little tedious): we first associate the $a_i$ by a bijection, then the $\cC_{i,j}$,
then the $C_i$. The exception arises only when we allow $\cB_2$, because for instance
an arc $\cC_{i,j}^k$ may pass very near $L_0$, and the corresponding arc of $\Sigma_{k+1}$
be split in two, with a vertex $\xi \in L_0$, or the other way around, or something similar
with a circle that passes near $\xi \in L_0$ and becomes a piece of $\bV$.
When splitting occurs, (\ref{12a11}) shows that the two arcs make an angle in $[\pi-C\varepsilon,\pi]$.

Let us now say how we define the tube $T_k$, starting with the case when the exception does not arise.
We have to say what is $T_k \cap \d B(0,r)$ when $r_{k+1} \leq r \leq r_k$, given that of course
$T_k \cap S_k = \Sigma_k$ and $T_{k+1} \cap S_{k+1} = \Sigma_{k+1}$. We write
$r = t r_k + (1-t)r_{k+1}$, $0 \leq t \leq 1$, and first place the vertices
$a_{i,r} \in \d B(0,r)$, as close as possible to $ta_i^k + (1-t)a_i^{k+1}$ 
(the latter may not have a norm 
exactly $r$); notice by the way that 
$a_{i,r} \in L_0$ when $a_i^k$ and $a_i^{k+1}$ lie on $L_0$. Then we can define
curves $\cC_{i,j,r} \subset \d B(0,r)$ that correspond to $\cC_{i,j}^k$ and $\cC_{i,j}^{k+1}$.
If $\H^1(\cC_{i,j}^k) \leq 2\pi r_k/3$, say, we just take the geodesic that connects the endpoints
$a_{i,r}$ and $b_{i,r}$ that correspond to the endpoints $a_i^k$ and $b_i^k$ of $\cC_{i,j}^k$
(and similarly for $\cC_{i,j}^{k+1}$). 
For longer arcs $\cC_{i,j}^k$, if they exist, this may be instable or ill-defined (for instance, if we want to connect two antipodal points of $\d B(0,r)$), 
so we cut $\cC_{i,j}^k$ into three pieces of equal length, do the same thing for $\cC_{i,j}^{k+1}$, and interpolate the three geodesics.
The new $\cC_{i,j,r}$ may not be a geodesic, but this does not matter.
For the full circles $C_i$, we can proceed similarly: we use three equally distant points $a, b, c \in C_i^k$
to cut $C_i^k$ into three equal arcs, then we choose three equally distant points 
$a', b', c' \in C_i^{k+1}$, and so that $|2a'-a|+|2b'-b|+|2c'-c|$
is minimal, for instance, and interpolate each of the three geodesics as before to create $C_{i,r}$,
which is thus the union of three geodesics of $\d B(0,r)$. 
The tube $T_k$ is the union of all these geodesics, with $r_{k+1} \leq r \leq r_k$.
It is not as smooth as we implicitly claimed it would be, because we added vertices and small
discontinuities of the tangent planes at some 
points, but this would be easily fixed with minor modifications of the construction.

The cases where there may be an exception above, which can only happen when we allow the bad set
$\cB_2$, are when for one of the two points $\xi$ of $L_0 \cap S_k$, either $\xi$ is a vertex
of $\Sigma_k$ but $\xi/2$ is not a vertex of $\Sigma_{k+1}$, or the other way around. Let us only
discuss the first case; the other one would be treated symmetrically.

By \eqref{12a6}, $\xi/2$ lies very close to $\Sigma_{k+1}$; we choose $\xi' \in \Sigma_{k+1}$ as close to $\xi/2$ as possible, and cut the
arc $\cC_{i,j}^{k+1}$ that contains $\xi'$ with the additional vertex $\xi'$. Recall that the two endpoints
of $\cC_{i,j}^{k+1}$ (call them $a'$ and $b'$) lie at distances $\geq C^{-1} r_k$ from $\xi'$
(because we excluded the proximity to bad $\bY$-cones). Then proceed as above, to interpolate
between the two curves between $\xi$ and $a, b$ (the points of $\Sigma_k$ that correspond 
to $a'$ and $b'$ in $\Sigma_{k+1}$) and the two curves of $\Sigma_{k+1}$ between $\xi'$ and $a',b'$.
In particular, if one of these two curves is too long, or even comes from a circle $C_i$, we also cut it
far from $\xi$ to interpolate in a more stable way.

At this point we have nice tubes $T_k$, and we can glue them to get $T = \cup_{k \geq 0} T_k$.
We now define a natural mapping $f_k : X_k \cap A_k \to T_k$, as follows.
After adding our extra vertices, we have a description of $\Sigma_k$ as a union of arcs of 
geodesics $\gamma_j$, $j\in J$, and for $r_{k+1} \leq r \leq r_k$, we have a similar description
of $T_k \cap \d B(0,r)$ as the union of corresponding arcs $\gamma_j^r$. We let $f_k$ be the only map
such that for $1/2 \leq t \leq 1$,  $f_k(t \gamma_j) = \gamma_j^{t r_k}$, which is run at constant speed.
That is, if $\gamma_j$ is the geodesic from $a_j$ to $b_j$, and its length is $\ell_j$,
and if $x \in \gamma_j$ lies at distance $\ell$ from $a_j$ along $\gamma_j$  (so $0 \leq \ell \leq \ell_j$),
then $\gamma_j^{t r_k}$ is the geodesic between $f_k(t a_j)$ and $f_k(t b_j)$, and  $f_k(tx)$ is the point on that geodesic that lies at distance $(\ell/\ell_j) \ell_j^{t r_k}$ from $f_k(t a_j)$ along $\gamma_j^{t r_k}$, 
where $\ell_j^{t r_k}$ denotes the length of the geodesic $\gamma_j^{t r_k}$. 
This is all more complicated
and specific than really needed, but the point is that with all these specific definitions, it would be easy 
(but very long, and we shall skip) to check that $f_k : X_k \cap A_k \to T_k$ is biLipschitz, 
and even with a biLipschitz constant that can be taken as close to $1$ as we want 
(by taking $\varepsilon$ accordingly small).
Only recall that we have lower bounds on all the lengths of all our geodesics, and the construction
makes them depend on $r$ in a Lipschitz way. 

Notice that $f_k(x) = x$ for $x\in \Sigma_k$; on the other side the restriction of $f_k$ to 
$X_k \cap \d B(0,r_{k+1}) = \frac12 \Sigma_k$
is a (biLipschitz) mapping that we shall call $g_{k+1} :  \frac12 \Sigma_k \to \Sigma_{k+1}$.
Let us extend $g_{k+1}$ by homogeneity to the whole cone $X_k$, 
so that $g_k(t x) = t g_k(x)$ for $x\in X_k$ and $t  >0$; 
note that $|g_k(x)| = |x|$ and $g_k$ maps $X_k$ to $X_{k+1}$. We also define
the $h_k : X_0 \to X_k$ by $h_0(x) = x$ and the induction relation
\begin{equation} \label{12a13}
    h_{k+1} = g_{k+1} \circ h_k = g_{k+1} \circ \ldots g_1.
\end{equation}

For each $k \geq 0$, Proposition \ref{l6a1} gives us a projection 
$\pi : W(c_1) \to E$, where $W(c_1)$ is a small conic neighborhood of $E$. 
For $x \in T =  \cup_k T_k$, $\mathrm{dist}(x,E) \leq C \varepsilon \abs{x}$ so, if $\varepsilon$ is small enough, $T$ is contained in $W(c_1)$,
and we shall concentrate on the restriction of $\pi$ to $T$. We claim that
\begin{equation} \label{12a14}
    \pi: T \to E \cap \ol B(0,r_0) \ \text{ is a biLipschitz bijection.}
\end{equation}
Let $x, y \in T$ be given; we want to estimate $|\pi(x)-\pi(y)|$, and
the most interesting case is 
when $\big||x|-|y|\big| \leq c |x|$,
with $c$ so small that we can use the same chart to define $\pi(x)$ and $\pi(y)$.
Then we can use our nice local description of $E$, $\pi$,
and of the $T_k$, to prove that $\big||\pi(x) - \pi(y)| - |x-y|\big| \leq C \varepsilon |x-y|$; we skip the details.
We get a similar estimate when 
when $\big||x|-|y|\big| \leq c |y|$.
When instead $\big||x|-|y|\big|$ is larger than $c |x|$ and $c|y|$,
we can use the fact that
$|\pi(x)-x| \leq C \dist(x, E) \leq C \varepsilon |x|$ and $|\pi(y)-y| \leq C \varepsilon |y|$
to get that $|\pi(x)-\pi(y)| \leq C \varepsilon (|x|+|y|) \leq C \varepsilon |x-y|$, as needed
for the Lipschitz part of \eqref{12a14}, 
while the lower bound holds because $|x-y| \leq |x|+|y| \leq  C \big||x|-|y|\big| 
= C\big||\pi(x)|-|\pi(y)|\big| \leq C |\pi(x)-\pi(y)|$
because $\pi$ acts separately on spheres.

Set $X = X_0$.  As the reader probably guessed, our biLipschitz parameterization will be of the form
$\psi = \pi \circ F$, where $F : X \cap \ol B(0,r_0) \to T$ is itself biLipschitz, and 
constructed with the help of the mappings above. We will define $F$ on each $X \cap A_k$ separately, 
so that $F: X \cap A_k \to T_k$. It is reasonable to take $F = f_0$ on $X \cap A_0 = X_0 \cap A_0$, 
and then we get that $F = f_0 = g_1$ on $\frac12\Sigma_0$. We shall take 
\begin{equation} \label{12a15}
    F(x) = f_k\circ h_k(x) 
    \in T_k 
    \ \text{ for } x \in X \cap A_k, \ k \geq 1,
\end{equation}
but let us verify a few things. Recall that $F = f_0$ in $A_0$.
Notice that $h_k : X \to X_k$, and also $h_k(X\cap A_k) \subset X_k \cap A_k$ because the 
$g_k$ and $h_k$ preserve the distance to the origin. Then recall that $f_k: X_k \cap A_k \to T_k$,
so $F(x)$ is defined and lies in $T_k$. Next we check that $F$ is continuous across each sphere 
$S_{k+1}$, $k \geq 0$. The definition from the $A_k$ side gives $F(x) = f_k\circ h_k(x) \in T_k$ on $S_{k+1}$,
but in fact $F(x) \in \Sigma_{k+1} = T_k \cap S_{k+1}$. In addition $h_k(x) \in X_k \cap S_{k+1}$,
and then its image by $f_k$ is also called $g_{k+1}(h_k(x)) = h_{k+1}(x)$. So $F(x) = h_{k+1}(x)$.
But now the definition coming from $A_{k+1}$ is also that $F(x) = f_{k+1}(h_{k+1}(x)) = h_{k+1}(x)$,
this time because $h_{k+1}(x) \in X_{k+1} \cap S_{k+1} = \Sigma_{k+1}$.
So our definition of $F: X \cap \ol B(0,r_0) \to T$ is coherent, $F$ is continuous, and we only need to
check that $F$ is biLipschitz. The main ingredient is the following.

\begin{lem}\label{l12a3}
    The mapping $h_k : X_0 \cap \d B(0,1) \to X_k \cap \d B(0,1)$ is $C$-biLipschitz, for some $C$ that depends only on $n$ and $\beta$.
\end{lem}

There is an issue here, because $h_k$ is the composition of an unbounded number
of biLipschitz mappings, but we have some rigidity, coming from the description of the $X_k$
and the fact that in the definition of $f_k$ we tried to parameterize the geodesics with 
constant speed.

We shall use the description of $X_k \cap \d B(0,1)$ as a union of geodesic arcs
$\cC_{i,j}$, $(i,j) \in I_a(r_k)$ and circles $C_i$, $i \in I_c(r_k)$, and start with
the simpler case where we also exclude $\cB_2$. 
In this case each function $g_{k+1}$, $k \geq 0$
(first extended by homogeneity and then restricted to the unit sphere), is a mapping
$g_{k+1} : X_k \cap \d B(0,1) \to X_{k+1} \cap \d B(0,1)$, which acts in a simple way
on the decomposition above: there are bijections from $I_a(r_k)$ to $I_a(r_{k+1})$
and from $I_c(r_k)$ to $I_c(r_{k+1})$, and then $g_{k+1}$ maps each geodesic $\cC_{i,j}$
or $C_i$ to the corresponding one for $k+1$, with constant speed (given by the ratio of lengths 
between the arc and its image). That is, even when we cut the
geodesics into three pieces because they were too long, we managed to map the whole geodesic
at constant speed. Then of course the same can be said about the composed mappings $h_k$.
Now there is a lower bound on the length of  each piece, which depends on $n$ and $\beta$
(because we also want a bound on the length of the short legs that go from $L_0$ to a
point of type $\bY$), so we have bounds on the speed of $h_k$ on each geodesic. 
The fact that $h_k$ is biLipschitz now follows, because it is easy to deduce from our description of
the $X_k \cap \d B(0,1)$ that the geodesic distance on $X_k \cap \d B(0,1)$ is equivalent to
the Euclidean distance.
Here again the bound depends on $\beta$, because we need a lower
bound on the distance of two arcs that do not share an endpoint.

\ms 
Let us now take care of the case when we allow $\cB_2$. Then the only case when 
$g_{k+1} : X_k \cap \d B(0,1) \to X_{k+1} \cap \d B(0,1)$ does not come from a bijection on
the arcs, with parameterizations at constant speeds, is when for some $\xi \in L_0 \cap \d B(0,1)$,
either $\xi$ is a vertex of $X_k \cap \d B(0,1)$ and not of $X_{k+1} \cap \d B(0,1)$, or the other way around. Notice that when this happens, the number of branches leaving from $\xi$ goes from
$2$ to $0$, and is never $1$. In fact, the local description of $E$ near a point of type $\bH$
(where $E$ is approximated by a half plane) shows that if  for some $r \in (0,r_0)$,
the point $\xi \in L_0 \cap \d B(0,r)$ corresponds to such a point (or in other words 
$X(r)$ has a single branch leaving from $\xi$), 
then this happens for every $r \in (0,r_0)$,
and $\xi/|\xi|$ never shows up in the discussion above.
For short, we will say in the first case mentioned above that the two arcs of $X_k \cap \d B(0,1)$
that leave from $\xi$ merge at $k$, and in the second case that the arc through $\xi$ splits at $k$. 

We would not want the same arc to be split a large number of times with no counterpart, because 
our estimate for the running speed comes from the fact that we always go from a collection of arcs of comparable
length to another one, at constant speed on each arc. What we fear is for instance that near a point, the arc
that contains this point is split (hence sent to two arcs, hence at potentially at roughly twice the speed) 
many times, while the opposite happens, but somewhere else and hence without compensating.
So we need to understand a little how splitting and merging occurs. For this we introduce auxiliary graphs,
$\Xi_k$, that in fact do not really depend on $k$. 

Recall that each $K_k = X_k \cap \d B(0,1)$ has a representation 
as a union of geodesics, which connect a set of vertices $V_k$; we choose a minimal representation, where the
vertices $x$ are either $Y$-point where $3$ geodesics leave, or else are points of $K_k \cap L$, and then exactly one or two geodesics leave from $x$, because we excluded $Y$-points on $K_k \cap L$. 
Set $V'_k = V_k \sm L$ (the triple points), and let $\Xi_k$ denote the graph with
$V'_k$ as the set of vertices, and edges chosen as follows. 
We say that the (non oriented) edge $(x,y)$ lies in the graph when 
the geodesic from $x$ to $y$ is one of the geodesics of the description of $K_k$, 
but also when there is a $\xi \in K_k \cap L$ such that both the geodesic from $x$ to $\xi$ and from $\xi$ to $y$ 
lie in $K_k$ (are members of its description). 
In this case we think of the union of these two geodesic as the representation of the edge $(x,y)$ in the graph.
We even add the edge $(x,y)$ when the three geodesics from $x$ to $\xi$, from $\xi$ to $-\xi$, and from $-\xi$ to $y$, belong also in the minimal representation of $K_k$, and then the representative of $(x,y)$ is 
the union of the three geodesics. 
Finally, we also add to our graph loops with no vertices, that correspond to full great circles that would lie in
$K_k$, and also pairs of geodesics from $\xi$ to $-\xi$ (the two points of $K_k \cap L$), when they lie in $K_k$. 

The main point of $\Xi_k$ is that it stays the same, in the sense that $\Xi_{k+1}$ is always isomorphic to
$\Xi_k$. Indeed, by the discussion above (the slow variation of $K_k$), the only difference between 
$\Xi_k$ and $\Xi_{k+1}$ could only occur when our geodesics split and merge. But when this happens,
we still keep the same set of vertices $V_k'$ (in reality, their representatives move a little but there is an obvious bijection), and maybe add or remove a point $\xi \in K_k \cap L$ 
from the representatives of an edge (exceptionally,
we may even add or remove two at the same time, when there is a double splitting or merging that concerns the same
geodesic). The same thing happens to the loops, as one could go from a circle to a union of two half circles, and back
(but this happens at the same time, as great circles either contain no point $\xi$ or both of them).

This is good for our uniform bilipschitz estimate, because
as before, we find that all our compositions $h_k$ of mapping $g_\ell$ are obtained
by parameterization with constant speed of the pieces that compose the representative of $K_k$; as before, we use the fact that the total number of pieces stays roughly the same (it may move by $1$ or $2$), and the lengths are bounded from below. Lemma \ref{l12a3} follows.
\qed

\ms
Now we have to deduce from the lemma that $F : X \cap \overline B(0,r_0) \to T$ 
is biLipschitz.
First we check that $h_k : X \cap  A_k \to X_k \cap  A_k$ is $C$-biLipschitz. It is Lipschitz 
because for $x, y \in X \cap  A_k$, 
\begin{eqnarray} \label{12a16}
    |h_k(x) - h_k(y)| &=& ||x| h_k(x/|x|) - |y| h_k(y/|y|)|
    \cr
    &\leq& |x| |h_k(x/|x|) - h_k(y/|y|)| + |(|x|-|y|) h_k(y/|y|)|
    \cr
    &\leq& |x| \Big|\frac{x}{|x|} - \frac{y}{|y|} \Big| + |x-y| \leq 3|x-y|
\end{eqnarray}
because $x \to x/|x|$ is $r_{k-1}^{-1}$-Lipschitz on $A_k$. And it has a Lipschitz inverse:
$X_k \cap  A_k \to X \cap A_k$, which is given by a similar formula $x \to |x| h_k^{-1}(x/|x|)$.

Now $f_k: X_k \cap  A_k \to T_k$ is also bilipschitz; this was checked a little above \eqref{12a13},
so by the formula \eqref{12a15} $F : X \cap A_k \to T_k$ is $C$-biLipschitz too. 
We still need to glue the pieces. 

Let $x \in X \cap A_k$ and $y\in X \cap A_\ell$ be given; we want to
estimate $|\psi(x)-\psi(y)|$. 
When $k=\ell$ we have the desired estimates, by \eqref{12a14} 
and since $\psi = \pi \circ F$. 
When $\ell \geq k+2$, we can conclude easily,
because $r_{k+1}/2 \leq |x-y| \leq 2r_k$, and since \eqref{12a15} and \eqref{7a2} say that 
$F(x) \in A_k$ and $F(y) \in A_\ell$, we also get that $\frac12 r_{k+1} \leq |\psi(x)-\psi(y)| \leq 2r_k$,
as needed. 
Modulo exchanging the names of $x$ and $y$ if needed, we are left with the case when
$x \in X \cap A_k$ and $y\in X \cap A_{k+1}$.
Set $\xi = r_{k+1} \frac{y}{|y|} \in X \cap S_k = X \cap A_k \cap A_{k+1}$.
Then $|\psi(x)-\psi(y)| \leq |\psi(x)-\psi(\xi)| + |\psi(\xi)-\psi(y)| \leq C |x-\xi|+ C |\xi-y|
\leq C |x-y|+ 2C |\xi-y| \leq 3C |x-y|$
because $F$ and $\psi$ are continuous across $S_{k+1}$, by the estimates on $A_k$ and
$A_{k+1}$, and because $|\xi-y| = r_{+1}-|y| \leq |x|-|y|$. For the lower bound,
$|\psi(x)-\psi(y)| \geq |\psi(x)-\psi(\xi)| - |\psi(\xi)-\psi(y)| \geq C^{-1} |x-\xi| - C |\xi-y|
\geq C^{-1} |x-y| - (C+C^{-1}) |\xi-y| \geq (2C)^{-1}|x-y|$ if
$2 C (C+C^{-1}) |\xi-y| \leq |x-y|$. Otherwise, 
$|\psi(x)-\psi(y)| \geq |\psi(x)|-|\psi(y)| = |x|-|y| \geq r_{k+1} - |y| =  |\xi-y| \geq 
[2 C (C+C^{-1})]^{-1}|x-y|$. This completes our biLipschitz estimate for $\psi$.

We still need to know that, when we excluded $\cB_2$, we have the sliding property \eqref{12a4}. 
Suppose that we can find $y\in X \cap L_0 \cap B(0,r_0)$. Then $X = X(r_0)$ meets 
$L_0$ at $y/|y|$, and at this point $X$ is either of type $\bH$ (coincides with a half plane near $y/|y|$),
or of type $\bV$ (by \eqref{12a6}), and with an angle which is far from flat and from sharp
(by \eqref{12a6} and \eqref{12a8}). In the first case, by the local description of $E$ near
a point of type $\bH$, we know that since this set is open, all the points of $\Gamma \cap B(0,r_0)$
are of type $\bH$, which implies that all the cones $X(r)$, $0 < r \leq r_0$, are of type $\bH$ at the
point $y/|y|$. In the second case, we have seen that there is no splitting or merging in the
process above, which means that all the cones $X(r)$, $0 < r \leq r_0$, contain $y/|y|$,
and even are of type $\bV$ there, with an angle which is far from sharp or flat. Because of
this, every point of $\Gamma$ on the same side of $0$ as $y/|y|$ lies in $E$
(and is a point of type $\bH$ or generic $\bV$).

Now we follow the construction of $\psi$. Let $k$ be such that $y \in A_k$.
First we map $y$ to $h_k(y)$, and since $y$ is a vertex of $X$ (because $y\in L_0$),
it is sent to $L_0$, in fact, to itself because $h_k$ preserves the norm. 
Then (according to \eqref{12a15}) we send $y$ to $f_k(y) \in T_k$. But again,
$f_k$ is constructed to keep the vertices in $L_0$, so $f_k(y) \in L_0$. In fact, since
$f_k$ also preserves the norm (it was constructed, above \eqref{12a13}, to act on spheres),
$f_k(y) = y$. Finally, $\psi(y) = \pi(y)$, where $\pi$ is the mapping of Proposition \ref{l6a1}.

Set $r = |y|$. With a good construction it should have happened that $y$ is the point $\zeta(r)$
of $S_r \cap \Gamma$ that lies close to $y$, but this is not what we did. This is unfortunate because
this way maybe $\pi(y)$ lies in one of the two branches of $E$ near $\zeta(r)$, and not precisely
on $\Gamma$. However this is not hard to fix, because we only need to manipulate the values
of our mappings near $L_0$. The least dirty way to proceed consists in defining a new bilipschitz
mapping $\tau$, that maps our tube $T = \cup_k T_k$ to a similar tube $\wt T$, with the property that
$|\tau(y)-y| \leq C \varepsilon |y|$ on $T$, and $\tau(y) = \zeta(|y|)$ when $y \in T \cap L_0$.
On the rest of $T$, just translate by a smooth extension of $\zeta(|y|) -y$. The fact that 
$|\tau(y)-y| \leq C \varepsilon |y|$ allows us to still use $\pi$ on $\wt T$, and finally
$\psi = \pi \circ \tau \circ F$ has all the desired properties.

This completes our proof of the precise variant of Theorem \ref{t12a2}, and then 
we saw earlier that Theorem \ref{t12a2}, Theorem \ref{t12a1}, and its precise variant follow.
\qed

\bigskip
\vfill \vfill \vfill \vfill
\noindent Guy David,  
\par\noindent 
Math\'ematiques, B\^atiment 307,
\par\noindent 
Universit\'e Paris-Saclay, CNRS, Laboratoire de math\'ematiques d'Orsay, 
\par\noindent 
91405, Orsay Cedex, France
\par\noindent 
guy.david@universite-paris-saclay.fr 
\par\noindent 
http://www.math.u-psud.fr/$\sim$gdavid/
\ms

\noindent Camille Labourie,
\par\noindent 
Université de Lorraine, CNRS, IECL,
\par\noindent
F-54000 Nancy, France
\par\noindent 
camille.labourie@univ-lorraine.fr
\par\noindent 
https://iecl.univ-lorraine.fr/membre-iecl/labourie-camille/
\ms


\begin{thebibliography}{AAA}


    \bibitem [A3]{AlmgrenMemoir} F. J. Almgren, 
        Existence and regularity almost everywhere 
        of solutions to elliptic variational problems with constraints.
        Memoirs of the Amer. Math. Soc. 165, volume 4 (1976), i-199.

    \bibitem [Da1]{Limits} 
        G. David, Limits of Almgren-quasiminimal sets. 
        Proceedings of the conference on Harmonic Analysis, 
        Mount Holyoke, A.M.S. Contemporary Mathematics series, Vol. 320 (2003), 119--145.
    \bibitem [Da2]{Dhh} 
        G. David, H\"older regularity of two-dimensional almost-minimal sets in $\R^n$. 
        Annales de la Facult\'e des Sciences de Toulouse, 
        Vol 18, 1 (2009), 65--246.
    \bibitem [Da3]{Dcc} 
        G. David, $C^{1+\alpha}$-regularity
        for two-dimensional almost-minimal sets in $\R^n$.
        J. Geom. Anal. 20 (2010), no. 4, 837--954.
    \bibitem [Da4]{DStein} 
        G. David, 
        Should we solve Plateau's problem again? 
        Advances in analysis: the legacy of Elias M. Stein, 108--145, Princeton Math. Ser., 50, 
        Princeton Univ. Press, Princeton, NJ, 2014.
    \bibitem [Da5]{Dss} 
        G. David, 
        Local regularity properties of almost- and quasiminimal sets with a sliding boundary condition. 
        Ast\'erisque No. 411 (2019).
    \bibitem [Da6]{Dvv} 
        G. David, 
        A local description of $2$-dimensional almost minimal sets bounded by a curve.
        Annales de la Facult\'e des Sciences de Toulouse,
        Vol XXXI, no. 1 (2022), 1--382.
    \bibitem [Da7]{Dpc} 
        G. David, 
        Sliding Almost Minimal Sets and the Plateau Problem.
        Harmonic analysis and applications,
        199--255, IAS/Park City Math. Ser., 27, Amer. Math. Soc., Providence, RI, 2020.

    \bibitem [DS1]{DS1}  G. David and S. Semmes, 
        Quasiminimal surfaces of codimension 1 and John domains. 
        Pacific J. Math. 183 (1998), no. 2, 213--277.

    \bibitem [DS2]{DS2} G. David and S. Semmes, 
        Uniform rectifiability and quasiminimizing sets of arbitrary codimension. 
        Mem. Amer. Math. Soc. 144 (2000), no. 687, viii+132 pp.

    \bibitem [Do]{Do} J. Douglas, Solutions of the Plateau problem. 
        Trans. Amer. Math. Soc. 33 (1931), no. 1, 263--321. 


    \bibitem[Fe]{Federer}  H. Federer, 
        \underbar{Geometric measure theory}. 
        Die Grundlehren der mathematischen Wissenschaften, 
        Band 153 Springer-Verlag New York Inc., New York 1969 xiv+676 pp.

    \bibitem [Fv1]{Fv1} V. Feuvrier, 
        Remplissage de l'espace Euclidien par des complexes poly\'edriques 
        d'orientation impos\'ee et de rotondit\'e uniforme. 
        Bull. Soc. Math. France 140 (2012), no. 2, 163--235.
    \bibitem [Fv2]{Fv2} 
        V. Feuvrier, 
        Condensation of polyhedric structures onto soap films. preprint, 2009,
        arXiv:0906.3505.

    \bibitem[Gi]{Giusti} E. Giusti,
        \underbar{Minimal surfaces and functions of bounded variation}.
        Monographs in Mathematics, 80. Birkh\"auser Verlag, Basel-Boston, 
        Mass., 1984.

    \bibitem [DGM]{Ita1} %
        C. De Lellis, F. Ghiraldin, and F. Maggi. 
        A direct approach to Plateau's problem. 
        J. Eur. Math. Soc. (JEMS) 19 (2017), no. 8, 2219--2240, DOI 10.4171/JEMS/716. 

    \bibitem [DDG]{Ita2} 
        G. De Philippis, A. De Rosa, and F. Ghiraldin, 
        A direct approach to Plateau's problem in any codimension. 
        Adv. Math. 288 (2016), 59--80.

    \bibitem[La]{La} C. Labourie,
        Labourie, Camille (F-UPS7-LM)
        Weak limits of quasiminimizing sequences.
        J. Geom. Anal. 31 (2021), no. 10, 10024–10135.

    \bibitem [Li1]{Li} 
        X. Liang, 
        Topological minimal sets and their applications.
        Calc. Var. Partial Differential Equations 47 (2013), no. 3-4, 523--546.

    \bibitem [Li2]{Li2} X. Liang, 
        On the topological minimality of unions of planes of arbitrary dimension. 
        Int. Math. Res. Not. IMRN 2015, no. 23, 12490--12539.

    \bibitem[Ma]{Mattila} 
        P. Mattila, 
        \underbar{Geometry of sets and measures in Euclidean space}.
        Cambridge Studies in Advanced Mathematics 44, Cambridge University Press l995.

    \bibitem [Mo1]{Mo1} F. Morgan, 
        Size-minimizing rectifiable currents.
        Invent. Math. 96 (1989), no. 2, 333--348.

    \bibitem [Mo2]{Mo} 
        F. Morgan, 
        \underbar{Geometric measure theory. A beginner's guide.}
        Fourth edition. Elsevier/Academic Press, Amsterdam, 2009. viii+249 pp. 

    \bibitem [Ra1]{Ra1} T. Rad\'o, 
        The problem of least area and the problem of Plateau. 
        Math. Z. 32 (1930), 763--796. 

    \bibitem [Ra2]{Ra2} T. Rad\'o, 
        On the problem of Plateau. Ergebnisse der 
        Mathematik und ihrer Grenzgebiete, Vol. 2, Springer, Berlin 1933; 
        Reprinted Chelsea, New York 1951.

    \bibitem [Re]{Re} E. R. Reifenberg, 
        Solution of the Plateau Problem for $m$-dimensional surfaces 
        of varying topological type.
        Acta Math. 104, 1960, 1--92.



    \bibitem[St]{Stein} E. M. Stein, 
        \underbar{Singular integrals and differentiability properties of functions.} 
        Princeton university press 1970.


    \bibitem [Ta]{Ta} J. Taylor, The structure of singularities in 
        soap-bubble-like  and soap-film-like minimal surfaces.
        Ann. of Math. (2) 103 (1976), no. 3, 489--539.

\end{thebibliography}
\end{document}